\documentclass[12pt]{article}

\usepackage[english]{babel}
\usepackage[T1]{fontenc}
\usepackage[applemac]{inputenc}
\usepackage{lmodern}
\usepackage{microtype}
\usepackage{bbm}
\usepackage{tikz}
\usepackage{graphicx}
\usepackage{subfloat}
\usepackage{amsmath}
\usepackage{amssymb}
\usepackage{amsthm}
\usepackage{latexsym}
\usepackage{color}
\usepackage{dsfont}
\usepackage{appendix}

\usepackage{hyperref}
\usepackage{float}
\usepackage{enumitem}
\usepackage{mathabx}
\usepackage{wrapfig}
\usepackage{graphicx} 
\graphicspath{{figures/}}

\DeclareSymbolFont{calletters}{OMS}{cmsy}{m}{n}
\DeclareSymbolFontAlphabet{\mathcal}{calletters}

\def\be{\begin{eqnarray}}
\def\ee{\end{eqnarray}}

\def\b*{\begin{eqnarray*}}
\def\e*{\end{eqnarray*}}

\def \E{\mathbb{E}}

\def \L{\mathbb{L}}

\def \N{\mathbb{N}}
\def \P{{\mathbb P}}
\def \Q{{\mathbb Q}}
\def \R{\mathbb{R}}

\def \eps{\varepsilon}

\def\Bc{{\cal B}}
\def\Cc{{\cal C}}
\def\Dc{{\cal D}}

\def\Fc{{\cal F}}

\def\Hc{{\cal H}}

\def\Jc{{\cal J}}
\def\Kc{{\cal K}}

\def\Mc{{\cal M}}
\def\Nc{{\cal N}}

\def\Pc{{\cal P}}

\def\Rc{{\cal R}}
\def\Sc{{\cal S}}
\def\Tc{{\cal T}}

\def\Xc{{\cal X}}
\def\Yc{{\cal Y}}

\def\Cfrak{{\mathfrak C}}

\def\Hb{\overline{H}}

\def\yb{{\bar y}}

\def\Ht{{\widetilde H}}

\def\Jt{{\widetilde J}}

\def \ri{{\rm ri\hspace{0cm}}}
\def \rf{{\rm rf\hspace{0cm}}}
\def \cl{{\rm cl\hspace{0.05cm}}}
\def \interior{{\rm int\hspace{0cm}}}
\def \dom{{\rm dom}}
\def \supp{{\rm supp}}

\def \csupp{{\wideparen{\supp}}}

\def \conv{{\rm conv}}

\def \aff{{\rm aff}}
\def \Aff{{\rm Aff}}

\def \ker{{\rm ker}}

\def \interf{{\rm interf}}

\def \Jo{J^\circ}
\def \Jco{\Jc^\circ}

\def \Leb{{\L}}
\def \Ctn{{\rm{C}}}

\def \Ibf{{\mathbf I}}
\def \Sbf{{\mathbf S}}
\def \Tbf{{\mathbf T}}

\def \med{{\mathfrak m}}

\def \dist{{\rm dist}}

\def \muxpw{\mu{\otimes}{\rm pw}}

\def \Kcirc{{\wideparen{\Kc}}}

\def \proj{{\rm proj}}

\def\no{\noindent}
\def\x{\times}

\def\={\;=\;}
\def\.{\;.}

\def\eps{\varepsilon}

\def \1{{\bf 1}}
\def \ep{\hbox{ }\hfill{ ${\cal t}$~\hspace{-5.1mm}~${\cal u}$   } }

\def \proof{{\noindent \bf Proof. }}
\def \ep{\hbox{ }\hfill$\Box$}

 \def\normeL2#1{\left\|{#1}\right\|_{L^2}}

\title{Quasi-sure duality for multi-dimensional martingale optimal transport\thanks{The author gratefully acknowledges the financial support of the ERC 321111 Rofirm, and the Chairs Financial Risks (Risk Foundation, sponsored by Soci\'et\'e G\'en\'erale) and Finance and Sustainable Development (IEF sponsored by EDF and CA).}}

\author{Hadrien De March\thanks{CMAP, \'Ecole Polytechnique, hadrien.de-march@polytechnique.org.}}

\date{\today}

\begin{document}

\newtheorem{Theorem}{Theorem}[section]
\newtheorem{Lemma}[Theorem]{Lemma}
\newtheorem{Corollary}[Theorem]{Corollary}
\newtheorem{Proposition}[Theorem]{Proposition}
\newtheorem{Remark}[Theorem]{Remark}
\newtheorem{Example}[Theorem]{Example}
\newtheorem{Definition}[Theorem]{Definition}
\newtheorem{Assumption}[Theorem]{Assumption}

\maketitle


\abstract{Based on the multidimensional irreducible paving of De March \& Touzi \cite{de2017irreducible}, we provide a multi-dimensional version of the quasi sure duality for the martingale optimal transport problem, thus extending the result of Beiglb\"ock, Nutz \& Touzi \cite{beiglbock2015complete}. Similar to \cite{beiglbock2015complete}, we also prove a disintegration result which states a natural decomposition of the martingale optimal transport problem on the irreducible components, with pointwise duality verified on each component. As another contribution, we extend the martingale monotonicity principle to the present multi-dimensional setting. Our results hold in dimensions 1, 2, and 3 provided that the target measure is dominated by the Lebesgue measure. More generally, our results hold in any dimension under an assumption which is implied by the Continuum Hypothesis. Finally, in contrast with the one-dimensional setting of \cite{beiglbock2017dual}, we provide an example which illustrates that the smoothness of the coupling function does not imply that pointwise duality holds for compactly supported measures.

\vspace{1cm}

}

\section{Introduction}

The problem of martingale optimal transport was introduced as the dual of the problem of robust (model-free) superhedging of exotic derivatives in financial mathematics, see Beiglb\"ock, Henry-Labord\`ere \& Penkner \cite{beiglbock2013model} in discrete time, and Galichon, Henry-Labord\`ere \& Touzi \cite{galichon2014stochastic} in continuous-time. This robust superhedging problem was introduced  by Hobson \cite{hobson1998robust}, and was addressing specific examples of exotic derivatives by means of corresponding solutions of the Skorokhod embedding problem, see \cite{cox2011robust,hobson2015robust,hobson2012robust}, and the survey \cite{hobson2011skorokhod}. 

Given two probability measures $\mu,\nu$ on $\R^d$, with finite first order moment, martingale optimal transport differs from standard optimal transport in that the set of all interpolating probability measures $\Pc(\mu,\nu)$ on the product space is reduced to the subset $\Mc(\mu,\nu)$ restricted by the martingale condition. We recall from Strassen \cite{strassen1965existence} that $\Mc(\mu,\nu)\neq\emptyset$ if and only if $\mu\preceq\nu$ in the convex order, i.e. $\mu(f)\le\nu(f)$ for all convex functions $f$. Notice that the inequality $\mu(f)\le\nu(f)$ is a direct consequence of the Jensen inequality, the reverse implication follows from the Hahn-Banach theorem.

This paper focuses on proving that quasi-sure duality holds in higher dimension, thus extending the results by Beiglb\"ock, Nutz and Touzi \cite{beiglbock2015complete} who prove that quasi-sure duality holds by identifying the polar sets. The structure of these polar sets is given by the critical observation by Beiglb\"ock \& Juillet \cite{beiglboeck2016problem} that, in the one-dimensional setting $d=1$, any such martingale interpolating probability measure $\P$ has a canonical decomposition $\P=\sum_{k\ge 0}\P_k$, where $\P_k\in\Mc(\mu_k,\nu_k)$ and $\mu_k$ is the restriction of $\mu$ to the so-called irreducible components $I_k$, and $\nu_k := \int_{x\in I_k}\P(dx,\cdot)$, supported in $J_k$ for $k\ge 0$, is independent of the choice of $\P\in\Mc(\mu,\nu)$. Here, $(I_k)_{k\ge 1}$ are open intervals, $I_0:=\R\setminus(\cup_{k\ge 1} I_k)$, and $J_k$ is an augmentation of $I_k$ by the inclusion of either one of the endpoints of $I_k$, depending on whether they are charged by the distribution $\P_k$.

In \cite{beiglbock2015complete}, this irreducible decomposition gives a form of compactness of the convex functions on each components, and plays a crucial role for the quasi-sure formulation, and represents an important difference between martingale transport and standard transport. Indeed, while the martingale transport problem is affected by the quasi-sure formulation, the standard optimal transport problem is not changed. We also refer to Ekren \& Soner \cite{ekren2016constrained} for further functional analytic aspects of this duality.

Our objective in this paper is to extend the quasi-sure duality, find a disintegration on the components, and a monotonicity principle for an arbitrary $d-$dimensional setting, $d\ge 1$. The main difficulty is that convex functions may lose information when converging. A first attempt to find such duality results was achieved by Ghoussoub, Kim \& Lim \cite{ghoussoub2015structure}. Their strategy consists in finding the largest sets on which pointwise monotonicity holds, and prove that it implies a pointwise existence of dual optimisers.

The paper is organized as follows. Section \ref{sect:relaxed} collects the main technical ingredients needed for the definition of the relaxed dual problem in view of the statement of our main results. Section \ref{sect:mainresults} contains the main results of the paper, namely the duality for the relaxed dual problem, the disintegration of the problem in the irreducible components identified in \cite{de2017irreducible}, and a monotonicity principle. In all the cases there are some claims that hold without any need of assumption, and a second part using Assumption \ref{ass:domination} defined in the beginning of the section. Section \ref{sect:examples} shows the identity with the Beiglb\"ock, Nutz \& Touzi \cite{beiglboeck2016problem} duality theorems in the one-dimensional setting, and provides non-intuitive examples, in particular Example \ref{expl:nopwdualitydim2} showing that there is no hope of having pointwise duality. The remaining sections contain the proofs of these results. In particular, Section \ref{sect:proofs} contains the proofs of the main results, and Section \ref{sect:assumption} checks the situations in which Assumption \ref{ass:domination} holds.
\\

\no {\bf Notation}\quad We denote by $\bar\R$ the completed real line $\R\cup\{-\infty,\infty\}$, and similarly denote $\overline{\R}_+:=\R_+\cup \{\infty\}$. We fix an integer $d\ge 1$. If $x\in\Xc$, and $A\subset \Xc$, where $(\Xc,{\rm d})$ is a metric space, $\dist(x,A):=\inf_{a\in A}{\rm d}(x,a)$. In all this paper, $\R^d$ is endowed with the Euclidean distance.

If $V$ is a topological affine space and $A\subset V$ is a subset of $V$, $\interior A$ is the interior of $A$, $\cl A$ is the closure of $A$, $\aff A$ is the smallest affine subspace of $V$ containing $A$, $\conv A$ is the convex hull of $A$, $\dim(A):=\dim(\aff A)$, and $\ri A$ is the relative interior of $A$, which is the interior of $A$ in the topology of $\aff A$ induced by the topology of $V$. We also denote by $\partial A:=\cl A\setminus\ri A$ the relative boundary of $A$. If $A$ is an affine subspace of $\R^d$, we denote by $\proj_A$ the orthogonal projection on $A$, and $\nabla A$ is the vector space associated to $A$ (i.e. $A-a$ for $a\in A$, independent of the choice of $a$). We finally denote $\Aff(V,\R)$ the collection of affine maps from $V$ to $\R$.

The set $\Kc$ of all closed subsets of $\R^d$ is a Polish space when endowed with the Wijsman topology\footnote{The Wijsman topology on the collection of all closed subsets of a metric space $(\Xc,{\rm d})$ is the weak topology generated by $\{\dist(x,\cdot):x\in\Xc\}$.} (see Beer \cite{beer1991polish}). As $\R^d$ is separable, it follows from a theorem of Hess \cite{hess1986contribution} that a function $F:\R^d\longrightarrow\Kc$ is Borel measurable with respect to the Wijsman topology if and only if
\b*
F^-(V):=\{x\in\R^d:F(x)\cap V\neq\emptyset\}&\mbox{is Borel for each open subset}~V\subset\R^d.
\e*
The subset $\Kcirc\subset\Kc$ of all the convex closed subsets of $\R^d$ is closed in $\Kc$ for the Wijsman topology, and therefore inherits its Polish structure. Clearly, $\Kcirc$ is isomorphic to $\ri\,\Kcirc := \{\ri K : K\in\Kcirc\}$ (with reciprocal isomorphism $\rm{cl}$). We shall identify these two isomorphic sets in the rest of this text, when there is no possible confusion.

We denote $\Omega:=\R^d\times\R^d$ and define the two canonical maps
\b*
 X :(x,y)\in\Omega
 \longmapsto x\in\R^d
 &\mbox{and}&
 Y :(x,y)\in\Omega
 \longmapsto y\in\R^d.
 \e*
For $\varphi,\psi:\R^d\longrightarrow\bar\R$, and $h:\R^d\longrightarrow\R^d$, we denote 
\b*
\varphi\oplus\psi
:=
\varphi(X)+\psi(Y),
&\mbox{and}&
h^\otimes := h(X)\cdot(Y-X),
\e*
with the convention $\infty-\infty = \infty$. Finally, for $A\subset \Omega$, and $x\in\R^d$, we denote $A_x:=\{y\in\R^d:(x,y)\in A\}$, and $A_x^c:=\{y\in\R^d:(x,y)\notin A\}$.

For a Polish space $\Xc$, we denote by $\Bc(\Xc)$ the collection of Borel subsets of $\Xc$, and $\Pc(\Xc)$ the set of all probability measures on $\big(\Xc,\Bc(\Xc)\big)$. For $\P\in\Pc(\Xc)$, we denote by $\Nc_\P$ the collection of all $\P-$null sets, $\supp\,\P$ the smallest closed support of $\P$, and $\csupp\,\P:=\cl\conv\,\supp\,\P$ the smallest convex closed support of $\P$. For a measurable function $f:\Xc\to\R$,  we use again the convention $\infty-\infty = \infty$ to define its integral, and we denote
 \b*
 \P[f]:=\E^\P[f] 
 = \int_\Xc f d\P 
 = \int_\Xc f(x) \P(dx)
 &\mbox{for all}&
 \P\in\Pc(\Xc).
 \e*
Let $\Yc$ be another Polish space, and $\P\in\Pc(\Xc\x\Yc)$. The corresponding conditional kernel $\P_x$ is defined by:
$$\P(dx,dy) = \mu(dx)\otimes \P_x(dy),\text{ where }\mu:=\P\circ X^{-1}.
$$
We denote by $\Leb^0(\Xc,\Yc)$ the set of Borel measurable maps from $\Xc$ to $\Yc$. We denote for simplicity $\Leb^0(\Xc):=\Leb^0(\Xc,\bar\R)$ and $\Leb^0_+(\Xc):=\Leb^0(\Xc,\bar\R_+)$. For a measure $m$ on $\Xc$, we denote $\Leb^1(\Xc,m):=\{f\in\Leb^0(\Xc):m[|f|]<\infty\}$. We also denote simply $\Leb^1(m):=\Leb^1(\bar\R,m)$ and $\Leb^1_+(m):=\Leb^1_+(\bar\R_+,m)$.

We denote by $\Cfrak$ the collection of all finite convex functions $f:\R^d\longrightarrow\R$. We denote by $\partial f(x)$ the corresponding subgradient at any point $x\in\R^d$. We also introduce the collection of all measurable selections in the subgradient, which is nonempty (see e.g. Lemma 9.2 in \cite{de2017irreducible}),
$$
\partial f:=\big\{p\in\L^0(\R^d,\R^d): p(x)\in\partial f(x)\text{ for all }x\in\R^d\big\}.
$$
Let $f:\R^d\longrightarrow\overline{\R}$, $f_{conv}(x):= \sup\{g(x)\mbox{ such that }g:\R^d\longrightarrow\overline{\R}\mbox{ is convex and }g \le f\}$ denotes the lower convex envelop of $f$. We also denote $\underline{f}_\infty := \liminf_{n\to\infty}f_n$, for any sequence $(f_n)_{n\ge 1}$ of real number, or of real-valued functions.

Let $I:\R^d\longmapsto \Kcirc$ be the irreducible components mapping defined in \cite{de2017irreducible}, which is the $\mu-$a.s. unique mapping such that for some $\widehat\P\in\Mc(\mu,\nu)$, $\ri\,\conv\,\supp\,\widehat\P_X = I(X)\supset \ri\,\conv\,\supp\,\P_X$, $\mu-$a.s. for all $\P\in\Mc(\mu,\nu)$.

\section{The relaxed dual problem}\label{sect:relaxed}
\setcounter{equation}{0}

\subsection{Preliminaries}

Throughout this paper, we consider two probability measures $\mu$ and $\nu$ on $\R^d$ with finite first order moment, and $\mu \preceq \nu$ in the convex order, i.e. $\nu(f)\ge \mu(f)$ for all $f\in\Cfrak$. Using the convention $\infty-\infty=\infty$, we may then define $(\nu-\mu)(f)\in[0,\infty]$ for all $f\in\Cfrak$. 

We denote by $\Mc(\mu,\nu)$ the collection of all probability measures on $\R^d\times\R^d$ with marginals $\P\circ X^{-1}=\mu$ and $\P\circ Y^{-1}=\nu$. Notice that $\Mc(\mu,\nu)\neq\emptyset$ by Strassen \cite{strassen1965existence}.

An $\Mc(\mu,\nu)-$polar set is an element of $\Nc_{\mu,\nu}:=\cap_{\P\in\Mc(\mu,\nu)}\Nc_\P$. A property is said to hold $\Mc(\mu,\nu)-$quasi surely (abbreviated as q.s.) if it holds on the complement of an $\Mc(\mu,\nu)-$polar set.

For a derivative contract defined by a non-negative cost function $c:\R^d\times\R^d\longrightarrow\R_+$, the martingale optimal transport problem is defined by:
 \be
 \Sbf_{\mu,\nu}(c)
 &:=&
 \sup_{\P\in\Mc(\mu,\nu)}
 \P[c].
 \ee

The corresponding robust superhedging problem is
 \be
 \Ibf_{\mu,\nu}(c)
 &:=&
 \inf_{(\varphi,\psi,h)\in\Dc_{\mu,\nu}(c)} \mu(\varphi)+\nu(\psi),
 \ee
where
 \be
 \Dc_{\mu,\nu}(c)
 &:=&
 \big\{ (\varphi,\psi,h)\in\L^1(\mu)\x\L^1(\nu)\x\L^0(\mu,\R^d):
                                                           ~\varphi\oplus\psi+h^\otimes\ge c
 \big\}.~~~~~
 \ee
 The following inequality is immediate:
  \be\label{eq:weakduality}
\Sbf_{\mu,\nu}(c) \le \Ibf_{\mu,\nu}(c).
\ee
This inequality is the so-called weak duality. For upper semi-continuous cost function, Beiglb\"ock, Henry-Labord\`ere, and Penckner \cite{beiglbock2013model} proved that there is no duality gap, i.e. $\Sbf_{\mu,\nu}(c)= \Ibf_{\mu,\nu}(c)$. See also Zaev \cite{zaev2015monge}. The objective of this paper is to establish a similar duality result for general measurable positive cost functions, thus extending the findings of Beiglb\"ock, Nutz, and Touzi \cite{beiglbock2015complete}.

For a probability $\P\in\Pc(\Omega)$, we say that $\P'\in\Pc(\Omega)$ is a competitor to $\P$ if $\P\circ X^{-1}=\P'\circ X^{-1}$, $\P\circ Y^{-1}=\P'\circ Y^{-1}$, and $\P[Y|X] = \P'[Y|X]$. Let $f:\Omega\longrightarrow \bar\R$, we say that a set $A\subset\Omega$ is $f-$martingale monotone if for all probability $\P$ having a finite support in $A$, and for all competitor $\P'$ to $\P$, we have $\P[f]\geq \P'[f]$.

\subsection{Tangent convex functions}

\begin{Definition}\label{def:Thetamunu}
Let $\theta:\Omega\to\overline{\R}_+$ be a universally measurable function, and a Borel set $N\in\Nc_{\mu,\nu}$ with $\{X= Y\}\subset N^c$. We say that $\theta$ is a $N-$tangent convex function if

\no{\rm (i)} $\theta(x,x)=0$, and $\theta(x,\cdot)$ is partially convex in $y$ on $N^c_x$;

\no{\rm (ii)} $N^c$ is $\theta-$martingale monotone;

\no{\rm (iii)} for all $\P$ with finite support in $N^c$, and any competitor $\P'$ to $\P$ such that $\supp\, \P'\cap N$ is a singleton, we have $\P'[N] = 0$;

\no{\rm (iv)} $A:=\{X\notin N_\mu\}\cap\{Y\in I(X)\}\subset N^c$, and $\mathbf{1}_A\theta$ is finite Borel measurable, for some $N_\mu\in\Nc_\mu$.
\end{Definition}

We denote by $\Theta_{\mu,\nu}$ the collection of all functions $\theta$ which are $N-$tangent convex for some $N$ as above. Clearly, $\Theta_{\mu,\nu}\supset \{\Tbf_p f:f\in\Cfrak,p\in\partial f\}$, where
\b*
\Tbf_pf(x,y)
:=
f(y)-f(x)-p^\otimes(x,y),
&\mbox{for all}&
f:\R^d\longmapsto\bar\R,~\mbox{and}~p:\R^d\longmapsto\R^d.
\e*

Indeed, for $f\in\Cfrak$, and $p\in\partial f$, $\Tbf_pf$ is convex in the second variable, thus satisfying (i) with $N=\emptyset$. For all $\P_0$ with finite support in $N^c=\Omega$, and $\P'$ competitor to $\P_0$, $\P_0[f(X)]=\P'[f(X)]$, $\P_0[f(Y)]=\P'[f(Y)]$, and $\P_0[p(X)\cdot(Y-X)]=\P_0[p(X)\cdot(\P_0[Y|X]-X)]=\P'[p(X)\cdot(Y-X)]$, and therefore $\P_0[\Tbf_pf]=\P'[\Tbf_pf]$.

\begin{Definition}
We say that a sequence $(\theta_n)_{n\ge 1}\subset\Theta_{\mu,\nu}$ generates some $\theta\in\Theta_{\mu,\nu}$ (and we denote $\theta_n\rightsquigarrow\theta$) if
\b*
\underline\theta_\infty\le\theta,&\mbox{and}&\P[\theta]\le\limsup_{n\to\infty}\P[\theta_n],~\mbox{for all }~\P\in\Pc(\Omega).
\e*
\end{Definition}

Notice that some sequences in $\Theta_{\mu,\nu}$ may generate infinitely many elements of $\Theta_{\mu,\nu}$. For example, for any nonzero $\theta\in\Theta_{\mu,\nu}$, the sequence $(\theta_n)_{n\in\N}:=(0,\theta,0,\theta,...)$ generates any $\theta'\in\Theta_{\mu,\nu}$ which is smaller than $\theta$. In particular $\theta_n\rightsquigarrow x\theta$, as $n$ goes to infinity, for all $0\le x \le 1$, which are uncountably many.

\begin{Definition}
{\rm (i)} A subset $\Tc\subset\Theta_{\mu,\nu}$ is semi-closed if $\theta\in\Tc$ for all $(\theta_n)_{n\ge 1}\subset\Tc$ generating $\theta$ (in particular, $\Theta_{\mu,\nu}$ is semi-closed).
\\
{\rm (ii)} The semi-closure of a subset $A\subset\Theta_{\mu,\nu}$ is the smallest semi-closed set containing $A$:
 $$
 \widetilde A
 :=
 \bigcap\big\{\Tc\subset\Theta_{\mu,\nu}:~
              A\subset\Tc,
              ~\text{and}~
              \Tc~\mbox{semi-closed}\,
         \big\}.
 $$
\end{Definition}

We next introduce for $a\ge 0$ the set $\Cfrak_a:=\big\{f\in\Cfrak:(\nu-\mu)(f)\le a\big\}$, and
$$\widetilde{\Tc}(\mu,\nu):=\underset{a\ge 0}{\bigcup}\,\widetilde{\Tc}_a,
\text{ where }
\Tc_a:=
 \big\{ \Tbf_p f: f\in\Cfrak_a,p\in\partial f\big\}.
 $$
 
\begin{Remark}
Notice that even though the construction of $\widetilde{\Tc}(\mu,\nu)$ is very similar to the construction of $\widehat{\Tc}(\mu,\nu)$ in \cite{de2017irreducible}, these objects may be different, see Lemma \ref{lemma:existencehat} below.
\end{Remark}

\begin{Proposition}\label{prop:convexonpolar}
$\widetilde{\Tc}(\mu,\nu)$ is a convex cone.

\end{Proposition}
\proof
The proof is similar to the proof of Proposition 2.9 in \cite{de2017irreducible}, using the fact that for $\theta,\theta_n,\theta_\infty\in\Theta_{\mu,\nu}$, the generation $\theta_n\rightsquigarrow \theta_\infty$ implies the generation $\theta_n+\theta\rightsquigarrow \theta_\infty+\theta$.
\ep

\subsection{Structure of polar sets}\label{subsect:polarmainresults}

The main results of this paper require the following assumption.

\begin{Assumption}\label{ass:domination}
\no{\rm (i)} For all $(\theta_n)_{n\ge 1}\subset \widetilde\Tc_1$, we may find $\theta\in\widetilde\Tc_1$ such that $\theta_n\rightsquigarrow \theta$.

\no{\rm (ii)} $I(X)\in\Cc\cup\Dc\cup\Rc$, $\mu-$a.s. for some subsets $\Cc,\Dc,\Rc\subset \Kcirc$ with $\Cc$ well ordered, $\dim(\Dc)\subset\{0,1\}$, and $\cup_{K\neq K'\in\Rc}\big[K\x (\cl K\cap\cl K')\big]\in\Nc_{\mu,\nu}$.
\end{Assumption}

\begin{figure}[h]
\centering

 \includegraphics[width=0.8\linewidth]{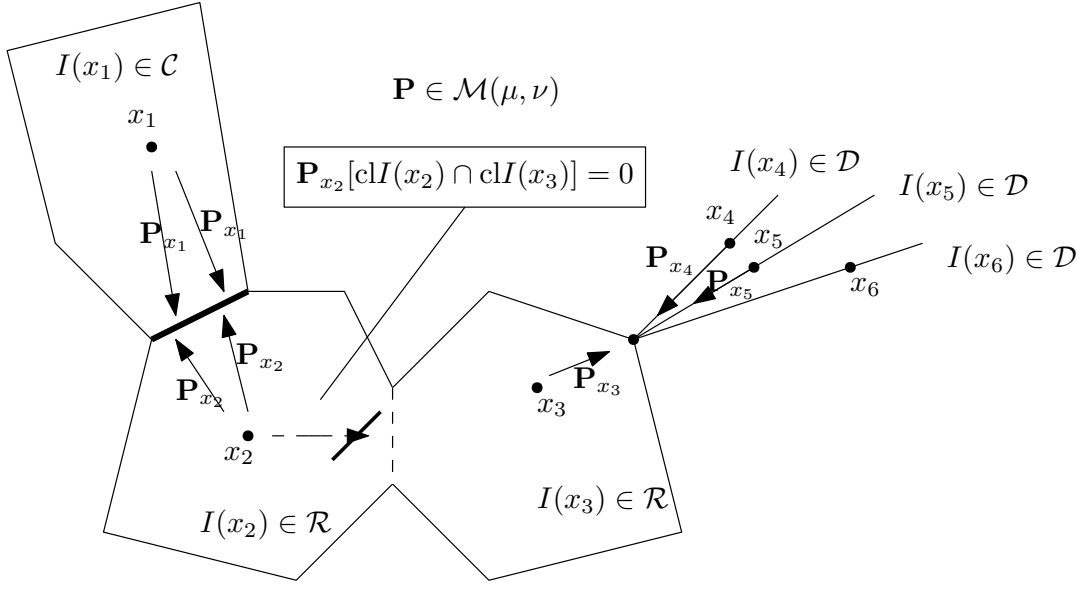}
    \caption{\label{fig:nofrontier} No communication between frontiers of elements in $\Rc$.}
\end{figure}

The condition $\cup_{K\neq K'\in\Rc}\big[K\x (\cl K\cap\cl K')\big]\in\Nc_{\mu,\nu}$ means that the probabilities in $\Mc(\mu,\nu)$ do not charge the intersections between frontiers of elements in $\Rc$, see Figure \ref{fig:nofrontier}.

We provide in Section \ref{sec:assumptions} some simple sufficient conditions for the last assumption to hold true. In particular, Assumption \ref{ass:domination} holds true in dimensions $d=1,2$, in dimension 3 with $\nu$ dominated by the Lebesgue measure, and in arbitrary dimension under the continuum hypothesis.

Recall that by Theorem 3.7 in \cite{de2017irreducible}, a Borel set $N\in\Bc(\Omega)$ is $\Mc(\mu,\nu)-$polar if and only if
\be\label{eq:polar}
N \subset
 \{X\in N_\mu\}\cup\{Y\in N_\nu\}\cup\{Y\notin J_\theta(X)\},
 \mbox{ for some }(N_\mu,N_\nu,\theta)\in\Nc_\mu\x\Nc_\nu\x\widehat{\Tc}(\mu,\nu),& 
 \ee
with $J_\theta:=\dom\theta(X,\cdot)\cap \bar{J}$, for some $I\subset \bar{J}\subset\cl I$, characterized $\mu-$a.s. by $\csupp\P_{X|\partial I(X)}\subset \bar{J}(X)\setminus I(X) = \csupp\widehat\P_{X|\partial I(X)}$, $\mu-$a.s., for some $\widehat\P\in\Mc(\mu,\nu)$, for all $\P\in\Mc(\mu,\nu)$. The definition of $\widehat{\Tc}(\mu,\nu)\subset \Leb^0_+(\Omega)$ is reported to Subsection \ref{subsect:reminders}. By Remark 3.5 in \cite{de2017irreducible}, $J_\theta$ is constant on $I(x)$ for all $x\in\R^d$. Then the random variable $J_\theta$ is $I-$measurable. Notice as well that by this remark we have
\b*
I\subset\underline{J} \subset J_\theta\subset \bar{J}\subset \cl I,&\mu-\mbox{a.s.}
\e*
Where $\underline{J}$ is characterized in Proposition 2.4 in \cite{de2017irreducible}. These sets $J_\theta$ are very important for characterising the polar sets. However they are not satisfactory as they may not be convex. We extend the notion in next proposition. Let $A\subset \Omega$, we say that $A$ is martingale monotone if for all finitely supported $\P\in\Pc(\Omega)$, and all competitor $\P'$ to $\P$, $\P[A]=1$ if and only if $\P'[A] = 1$. Notice that $A$ is martingale monotone if and only if $A$ is $\mathbf{1}_{A^c}-$martingale monotone.

\begin{Proposition}\label{prop:J}
Under Assumption \ref{ass:domination}, for any $N-$tangent convex $\theta\in\widetilde\Tc(\mu,\nu)$, we may find $\theta\le\theta'\in\widetilde\Tc(\mu,\nu)$ and $(N_\mu^0, N_\nu^0)\in\Nc_\mu\x\Nc_\nu$ such that for all $(N_\mu^0,N_\nu^0)\subset (N_\mu,N_\nu)\in\Nc_{\mu}\x\Nc_{\nu}$, the maps $I$, $\underline{J}$, and $\bar J$ from \cite{de2017irreducible} may be chosen so that $J(X):=\conv(\dom\theta'(X,\cdot)\setminus N_\nu)\cap \aff\,I(X)$ satisfies, up to a modification on $N_\mu$:

\no{\rm (i)} $J(X)=\conv\big(J(X)\setminus N_\nu\big)$, and on $N_\mu^c$, we have $J(X)\subset\dom\theta(X,\cdot)$;

\no{\rm (ii)} $N\subset N':=\{X\in N_\mu\}\cup\{Y\in N_\nu\}\cup\{Y\notin J(X)\}\in\Nc_{\mu,\nu}$ and $N'^c$ is martingale monotone;

\no{\rm (iii)} the set-valued map $\Jo(x):=\cup_{x'\in J(x)\setminus N_\mu}I(x')\cup \left(J(x)\setminus N_\nu\right)\cup\{x\}$ satisfies $\underline{J}\subset \Jo\subset J\subset \bar{J}$, furthermore $J$ and $\Jo$ are constant on I(x), for all $x\in\R^d$.
\end{Proposition}

 
The proof of Proposition \ref{prop:J} is reported in Subsection \ref{subsect:polarproofs}. We denote by $\Jc(\mu,\nu)$ (resp. $\Jco(\mu,\nu)$) the set of these modified set-valued mappings $J$ \big(resp. $\Jo$\big) from Proposition \ref{prop:J}.

\begin{Remark}\label{rmk:J}
Let $J\in\Jc(\mu,\nu)$, $N_\nu\in\Nc_\nu$, and $\Jo\in\Jco(\mu,\nu)$ from Proposition \ref{prop:J}. The following holds for $\Jt \in \{J,\Jo,J\setminus N_\nu\}$. Let $x, x'\in\R^d$,

\no{\rm (i)}  $Y\in \Jt(X)$, $\Mc(\mu,\nu)-$q.s.;

\no{\rm (ii)} $\Jt(x)\cap \Jt(x') = \aff\left(\Jt(x)\cap \Jt(x')\right)\cap \Jt(x)$;

\no{\rm (iii)} $J(x)\cap J(x') = \conv\left(\Jt(x)\cap \Jt(x')\right)$; 

\no{\rm (iv)} if $I(x')\cap \Jt(x)\neq \emptyset$, then $\Jt(x')\subset \Jt(x)$.
 \end{Remark}

Remark \ref{rmk:J} will be justified in Subsection \ref{subsect:polarproofs}. We next introduce a subset of polar sets which play an important role.

\begin{Definition}
We say that $N\in\Nc_{\mu,\nu}$ is canonical if $N=\{X\in N_\mu\}\cup\{Y\in N_\nu\}\cup\{Y\notin J(X)\}$, for some $(N_\mu,N_\nu,J)\in\Nc_\mu\x\Nc_\nu\x\Jc(\mu,\nu)$ from Proposition \ref{prop:J} for some $\theta\in\widetilde{\Tc}(\mu,\nu)$.
\end{Definition}
 
 \begin{Theorem}\label{thm:polar}
Under Assumption \ref{ass:domination}, an analytic set $N\subset\Omega$ is $\Mc(\mu,\nu)-$polar if and only if it is contained in a canonical $\Mc(\mu,\nu)-$polar set.
 \end{Theorem}
 
 The proof of Theorem \ref{thm:polar} is reported in Subsection \ref{subsect:polarproofs}. 
 \begin{Remark}\label{rmk:Jmeasurable}
 For a fixed $x\in\R^d$, even though $J(x)$ is convex for $J\in\Jc(\mu,\nu)$, it may not be Borel anymore, unlike $J_\theta(x)$ when $\theta\in\widehat{\Tc}(\mu,\nu)$. The same holds for $\Jo(x)$, with $\Jo\in\Jco(\mu,\nu)$ or for a canonical polar sets, they may not be Borel but only universally measurable (i.e. $\P-$measurable\footnote{A set $A$ is said to be $\P-$measurable if $\P\big[(A\cup B)\setminus (A\cap B)\big] = 0$ for some Borel set $B\subset\Omega$.} for all $\P\in\Pc(\Omega)$). Similar to $J_\theta$ for $\theta\in\widehat{\Tc}(\mu,\nu)$, the invariance of $J\in\Jc(\mu,\nu)$ and $\Jo\in\Jco(\mu,\nu)$ on $I(x)$ for each $x\in\R^d$ proves that $J$ is $I-$measurable.
 \end{Remark}

\subsection{Weakly convex functions}

We see from \cite{beiglbock2015complete} 4.2 that the integral of the dual functions needs to be compensated by a convex (concave in \cite{beiglbock2015complete}) moderator to deal with the case $\mu[\varphi]+\nu[\psi] = -\infty+\infty$. However, they need to define a new concave moderator for each irreducible component before summing them up on the countable components. In higher dimension, as the components may not be countable there may be measurability issues arising. We need to store all these convex moderators in one single moderator which is convex on each component, but that may not be globally convex (see Example \ref{expl:Mmunuconvexfunction}).

\begin{Definition}
A function $f:\R^d\longrightarrow\R$ is said to be $\Mc(\mu,\nu)$-convex or weakly convex if there exists a tangent convex function $\theta\in\widetilde\Tc(\mu,\nu)$ such that
\b*
\Tbf_pf = \theta,& \mbox{on }\{Y\in \Jo(X), X\notin N_\mu\},&\mbox{for some}~p:\R^d\to \R^d,~\mbox{and}~(N_\mu,\Jo)\in\Nc_\mu\x\Jco(\mu,\nu).
\e*
\end{Definition}

Under these conditions, we write that $\theta \approx \Tbf_pf$. Notice that by Remark \ref{rmk:J}, $Y\in\Jo(X)$, $\Mc(\mu,\nu)-$q.s., whence $\theta\approx\Tbf_pf$ implies that $\theta=\Tbf_pf$, $\Mc(\mu,\nu)-$q.s. We denote by $\Cfrak_{\mu,\nu}$ the collection of all $\Mc(\mu,\nu)$-convex functions. Similarly to convex functions, we introduce a convenient notion of subgradient:
 \b*
 \partial^{\mu,\nu} f
 &:=&
 \big\{ p:\R^d\longmapsto\R^d : \Tbf_p f \approx \theta\in\widetilde\Tc(\mu,\nu)
 \big\},
 \e*
which is by definition non-empty. A key ingredient for all the results of this paper is that the sets $\Theta_{\mu,\nu}$ and $\Cfrak_{\mu,\nu}$ turn out to be in one-to-one relationship. 

\begin{Proposition}\label{prop:CctoTc}
Under Assumption \ref{ass:domination},
$$\widetilde\Tc(\mu,\nu)=\{\theta \approx \Tbf_pf,\mbox{ for some }f\in\Cfrak_{\mu,\nu},\mbox{ and }p\in\partial^{\mu,\nu}f\}.$$
\end{Proposition}

The proof of this proposition is reported in Subsection \ref{subsect:CctoTc}.

\begin{figure}[h]
\centering

 \includegraphics[width=0.8\linewidth]{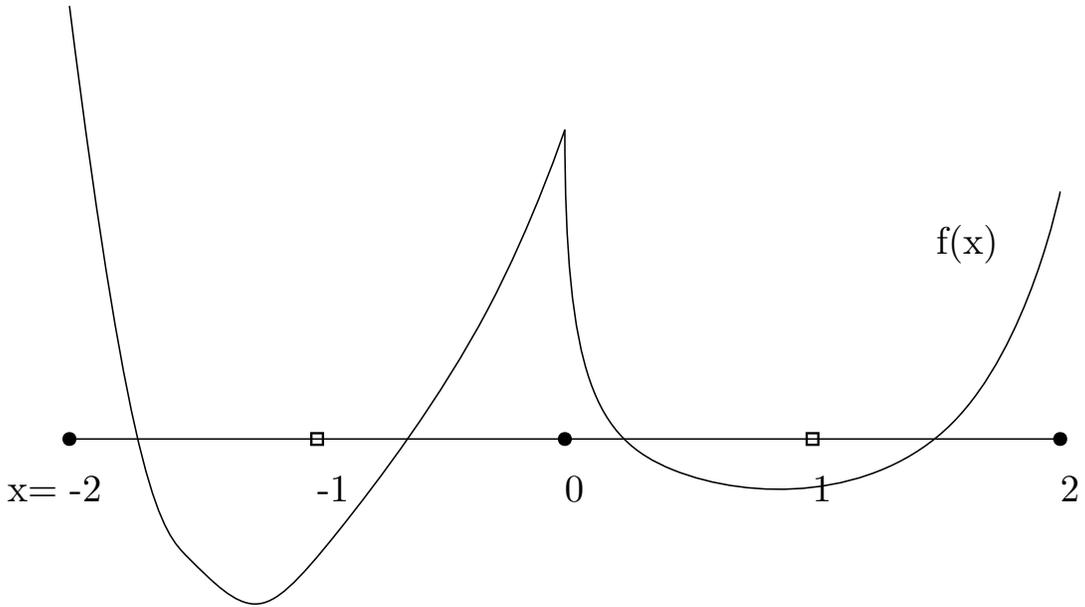}
    \caption{\label{fig:qsconvex} Example of a $\Mc(\mu,\nu)-$convex function.}
\end{figure}

\begin{Example}\label{expl:Mmunuconvexfunction}[$\Mc(\mu,\nu)-$convex function in dimension one]
Let $\mu:=\frac12(\delta_{-1}+\delta_{1})$, and $\nu(dy):=\frac18\big(\mathbf{1}_{[-2,2]}(y)dy+\delta_{-2}(dy)+2\delta_{0}(dy)+\delta_{2}(dy)\big)$. For these measures, one can easily check that the irreducible components from \cite{beiglboeck2016problem}, \cite{beiglbock2015complete}, and \cite{de2017irreducible} are given by $I(-1)=(-2,0)$, and $I(1)=(0,2)$, and the associated $\bar{J}$ mapping is given by $\bar{J}(-1)=[-2,0]$, and $\bar{J}(1)=[0,2]$. By Example \ref{expl:Setsindimone} in this paper, $f:\R\longrightarrow\R$ is $\Mc(\mu,\nu)-$convex if it is convex on each irreducible components. See Figure \ref{fig:qsconvex}.
\end{Example}

The next result shows that the weakly convex functions are convex on each irreducible component. Let $\eta:=\mu\circ I^{-1}$, and recall that any $J\in\Jc(\mu,\nu)$ is $I-$measurable by Remark \ref{rmk:Jmeasurable}.

\begin{Proposition}\label{prop:Mmunuconvex}
Let $f\in\Cfrak_{\mu,\nu}$ and $p\in\partial^{\mu,\nu}f$. Then $f$ is convex on $\Jo$, and $\proj_{\nabla\aff \Jo}(p)(X)\in \partial f|_{\Jo}(X)$, $\mu-$a.s. Furthermore, we may find $\widetilde{f}\in\Cfrak_{\mu,\nu}$ and $\widetilde{p}\in\partial^{\mu,\nu}\widetilde{f}$ such that $f=\widetilde{f}$, $\mu+\nu-$a.s., $\widetilde{p} = \proj_{\nabla\aff \Jo}(p)$, $\mu-$a.s., and $\widetilde{f}$ is convex on $J$ with $\widetilde{p}\in\partial\widetilde{f}|_I$, $\eta$-a.s. for some $J\in{\Jc}(\mu,\nu)$.
\end{Proposition}

The proof of this proposition is reported in Subsection \ref{subsect:CctoTc}.

\subsection{Extended integrals}

The following integral is clearly well-defined:
\be\label{nu-mu:0}
(\nu-\mu)[f]
=
\P[\Tbf_pf] 
&\mbox{for all}&
\P\in\Mc(\mu,\nu),~f\in\Cfrak\cap\L^1(\nu),\, p\in\partial f.
\ee
Similar to Beiglb\"ock, Nutz \& Touzi \cite{beiglbock2015complete}, we need to introduce a convenient extension of this integral. For $f\in\Cfrak_{\mu,\nu}$, define:
 \be\label{eq:nu-mu}
 \nu\overline{\ominus}\mu[f]
 &:=&
 \inf\big\{a\ge 0 :\Tbf_p f\approx\theta\in \widetilde{\Tc}_a,~\mbox{for some}~
 p\in\partial^{\mu,\nu}f\big\} 
 \ee
 \be
 \nu\underline{\ominus}\mu[f]
 :=
 \Sbf_{\mu,\nu}(\Tbf_pf),
 &\mbox{for}&
 p\in\partial^{\mu,\nu}f,
 \ee
where the last value is not impacted by the choice of $p\in\partial^{\mu,\nu}f$, whenever $\nu\underline{\ominus}\mu[f]<\infty$. Indeed, if $p_1,p_2\in \partial^{\mu,\nu}f$ such that $\P[\Tbf_{p_1}f]<\infty$ and $\P[\Tbf_{p_2}f]<\infty$, then $\Tbf_{p_1}f-\Tbf_{p_2}f=(p_2-p_1)^\otimes\in\L^1(\P)$, and it follows from the Fubini theorem that $\P[\Tbf_{p_1}f-\Tbf_{p_2}f] = \P[(p_2-p_1)^\otimes] = \P[\P[(p_2-p_1)^\otimes|X]] = 0$.

We also abuse notation and define for $\theta\in\widetilde{\Tc}(\mu,\nu)$, $\nu\overline{\ominus}\mu{}[\theta]:=\inf\big\{a\ge 0 :\theta\in \widetilde{\Tc}_a\big\}$.

\begin{Proposition}\label{prop:numoinsmuqs}
For $f\in\Cfrak_{\mu,\nu}$ and $\theta\in\widetilde{\Tc}(\mu,\nu)$, we have
\\
{\rm (i)} $\nu\overline{\ominus}\mu[f]\ge\nu\underline{\ominus}\mu[f]\ge 0$, and $\nu\overline{\ominus}\mu{}[\theta]\ge\Sbf_{\mu,\nu}(\theta)\ge 0$;
\\
{\rm (ii)} if $f\in\Cfrak\cap\L^1(\nu)$, then $\nu\overline{\ominus}\mu[f] =\nu\underline{\ominus}\mu[f] = \nu\overline{\ominus}\mu{}[\Tbf_pf] = (\nu-\mu)[f]$, for all $p\in\partial f$;
\\
{\rm (iii)} $\nu\underline{\ominus}\mu$ and $\nu\overline{\ominus}\mu$ are homogeneous and convex.
\end{Proposition}
\proof
The proof is similar to the proof of Proposition 2.11 in \cite{de2017irreducible}.
\ep

We can prove the next simple characterization of $\widetilde{\Tc}(\mu,\nu)$, $\Cc(\mu,\nu)$ and $\widehat\Tc(\mu,\nu)$ in the one-dimensional setting. In dimension $1$, by Beiglb\"ock, Nutz \& Touzi \cite{beiglbock2015complete}, there are only countably many irreducible components of full dimension. The other components are points. Then we can write these components $I_k$ for $k\in\N$ like in \cite{beiglbock2015complete} Proposition 2.3. We also have uniqueness of the $J(x)$ from Theorem 3.7 in \cite{de2017irreducible}, that is equivalent in dimension $1$ to Theorem 3.2. We denote them $J_k$ as well. We also take another notation from the paper, $\mu_k$ and $\nu_k$ the restrictions of $\mu$ and $\nu$ to $I_k$ and $J_k$, and $(\nu_k-\mu_k)$ extending their Definition 4.2 to non integrable convex functions, which corresponds to the operator $\nu{\ominus}\mu$ in this paper. 

\begin{Example}\label{expl:Setsindimone}
If $d=1$,
\b*
\Cfrak_{\mu,\nu} &=& \Big\{f:\R^d\to\R : f_{|J_k}\text{ is convex for all }k\Big\},\\
\widetilde{\Tc}(\mu,\nu) &=& \Big\{\theta = \sum_k\mathbf{1}_{X\in I_k}\Tbf_{p_k} f_k:f_k \text{ convex finite on }J_k,\,p_k\in\partial f_k,
\text{ and }\underset{k}{\sum}(\nu_k-\mu_k)(f_k)< \infty\Big\},\\
&&\text{ and }\\
\nu\overline{\ominus}\mu[f]&=&\nu\underline{\ominus}\mu[f]=\underset{k}{\sum}(\nu_k-\mu_k)(f_{|J_k}),\quad\mbox{for all }f\in\Cfrak_{\mu,\nu}.
\e*
This characterization follows from the same argument than the proof of Proposition 3.11 in \cite{de2017irreducible}.
\end{Example}

\subsection{Problem formulation}

\begin{Definition}
Let $\varphi,\psi:\R^d\longrightarrow\R$ and $f\in\Cfrak_{\mu,\nu}$. We say that $f$  is a convex moderator for $(\varphi,\psi)$ if
 \b*
 \varphi+f \in \L_+^1(\mu),~~
 \psi-f \in \L_+^1(\nu),
 &\mbox{and}&
 \nu{\ominus}\mu[f]  
 :=
 \nu\overline{\ominus}\mu[f]=\nu\underline{\ominus}\mu[f]<\infty.
 \e*
We denote by $\widehat\L(\mu,\nu)$ the collection of triplets $(\varphi,\psi,h)$ such that $(\varphi,\psi)$ has some convex moderator $f$ with $h+p\in\Leb^0(\R^d,\R^d)$ for some $p\in\partial^{\mu,\nu}f$.
\end{Definition}

We now introduce the objective function of the robust superhedging problem for a pair $(\varphi,\psi)\in\widehat\L(\mu,\nu)$ with convex moderator $f$:
 \be
 \mu[\varphi]{\oplus}\nu[\psi]
 &:=&
 \mu[\varphi+f]+\nu[\psi-f]+\nu{\ominus}\mu[f].
 \ee
We observe immediately that this definition does not depend on the choice of the convex moderator. Indeed, if $f_1,f_2$ are two convex moderators for $(\varphi,\psi)$, it follows that $f_1-f_2\in\L^1(\mu)\cap\L^1(\nu)$, and consequently $\mu\underline{\ominus}\nu[f_1]=\mu\underline{\ominus}\nu[f_2]+(\nu-\mu)[f_1-f_2]$ by Proposition \ref{prop:numoinsmuqs}. This implies that
 \b*
 \mu[\varphi+f_1]+\nu[\psi-f_1]+\nu{\ominus}\mu[f_1]
 &=&
 \mu[\varphi+f_2]+\nu[\psi-f_2]+\nu{\ominus}\mu[f_2].
 \e* 

For a cost function $c:\R^d\times\R^d\longrightarrow\R_+$, the relaxed robust superhedging problem is
 \be
 \Ibf^{qs}_{\mu,\nu}(c)
 &:=&
 \inf_{(\varphi,\psi,h)\in\Dc^{qs}_{\mu,\nu}(c)} \mu[\varphi]{\oplus}\nu[\psi],
 \ee
where
 \be
 \Dc^{qs}_{\mu,\nu}(c)
 &:=&
 \big\{ (\varphi,\psi,h)\in\widehat\L(\mu,\nu):
 ~\varphi\oplus\psi+h^\otimes\ge c,~\Mc(\mu,\nu)-\mbox{q.s.}
 \big\}.~~~~~
 \ee

\begin{Remark}
This dual problem depends on the primal variables $\Mc(\mu,\nu)$. However this issue is solved by the fact that Theorem 3.7 in \cite{de2017irreducible} gives an intrinsic description of the polar sets. See also Theorem \ref{thm:polar}.
\end{Remark}
 
We also introduce the pointwise version of the robust superhedging problem:
 \be
 \Ibf^{\rm pw}_{\mu,\nu}(c)
 &:=&
 \inf_{(\varphi,\psi,h)\in\Dc^{\rm pw}_{\mu,\nu}(c)} \mu[\varphi]{\oplus}\nu[\psi],
 \ee
where
 \be
 \Dc^{\rm pw}_{\mu,\nu}(c)
 &:=&
 \big\{ (\varphi,\psi,h)\in\widehat\L(\mu,\nu):
                                                           ~\varphi\oplus\psi+h^\otimes\ge c
 \big\}.~~~~~
 \ee

The following inequalities extending the classical weak duality \eqref{eq:weakduality} are immediate,
\be 
\Sbf_{\mu,\nu}(c)\le \Ibf^{qs}_{\mu,\nu}(c)\le\Ibf^{\rm pw}_{\mu,\nu}(c).
\ee

\section{Main results}\label{sect:mainresults}
\setcounter{equation}{0}

\begin{Remark}
All the results in this section are given for $c\ge 0$. The extension to the case $c\geq \varphi_0\oplus\psi_0+h_0^\otimes$ with $(\varphi_0,\psi_0,h_0)\in \L^1(\mu)\times\L^1(\nu)\times\L^1(\mu,\R^d)$, is immediate by applying all results to $c-\varphi_0\oplus\psi_0-h_0^\otimes\ge 0$.
\end{Remark}

\subsection{Duality and attainability}

We recall that an upper semianalytic function is a function $f:\R^d\to\R$ such that $\{f\ge a\}$ is an analytic set for any $a\in\R$. In particular, a Borel function is upper semianalytic.

\begin{Theorem}\label{thm:duality}
Let $c:\Omega\to\overline{\R}_+$ be upper semianalytic. Then, under Assumption \ref{ass:domination}, we have
\\
{\rm (i)} $\Sbf_{\mu,\nu}(c) = \Ibf^{qs}_{\mu,\nu}(c)$;
\\
{\rm (ii)} If in addition $\Sbf_{\mu,\nu}(c)<\infty$, then existence holds for the quasi-sure dual problem $\Ibf^{qs}_{\mu,\nu}(c)$.
\end{Theorem}

This Theorem will be proved in Subsection \ref{subsect:proofduality}.

\begin{Remark}
{\rm For an upper-semicontinuous coupling function $c$, we observe that the duality result $\Sbf_{\mu,\nu}(c) = \Ibf^{qs}_{\mu,\nu}(c) = \Ibf^{pw}_{\mu,\nu}(c)$ holds true, together with the existence of an optimal martingale interpolating measure for the martingale optimal transport problem $\Sbf_{\mu,\nu}(c)$, without any need to Assumptions \ref{ass:domination}. This is an immediate extension of the result of Beiglb\"ock, Henry-Labord\`ere \& Penckner \cite{beiglbock2013model}, see also Zaev \cite{zaev2015monge}. However, dual optimizers may not exist in general, see the counterexamples in Beiglb\"ock, Henry-Labord\`ere \& Penckner and in Beiglb\"ock, Nutz \& Touzi \cite{beiglbock2015complete}. Observe that in the one-dimensional case, Beiglb\"ock, Lim \& Ob{\l}{\'o}j \cite{beiglbock2017dual} proved that pointwise duality, and integrability hold for $\Ctn^2$ cost functions together with compactly supported $\mu$, and $\nu$. We show in Example \ref{expl:nopwdualitydim2} below that this result does not extend to higher dimension.
}
\end{Remark}

\begin{Remark}
{\rm
An existence result for the robust superhedging problem was proved by Ghoussoub, Kim \& Lim \cite{ghoussoub2015structure}. We emphasize that their existence result requires strong regularity conditions on the coupling function $c$ and duality, and is specific to each component of the decomposition in irreducible convex pavings, see Subsection \ref{sect:irreducible pavings} below. In particular, their construction does not allow for a global existence result because of non-trivial measurability issues. Our existence result in Theorem \ref{thm:duality} (ii) by-passes these technical problems, provides global existence of a dual optimizer, and does not require any regularity of the cost function $c$.
}
\end{Remark}

\subsection{Decomposition on the irreducible convex paving}
\label{sect:irreducible pavings}

The measurability of the map $I$ stated in Theorem 2.1 (i) in \cite{de2017irreducible}, induces a decomposition of any function on the irreducible paving by conditioning on $I$. We shall denote $\eta:= \mu\circ I^{-1}$, and set $\mu_I:= \mu\circ (X|X\in I)^{-1}$. Then for any measurable $f:\R^d \longrightarrow\R$, non-negative or $\mu-$integrable, we have
 $
 \int_{\R^d} f(x) \mu(dx) 
 =
 \int_{I(\R^d)}\left(\int_I f(x)\mu_I(dx)\right)\eta(dI).
 $

Similar to the one-dimensional context of Beiglb\"ock, Nutz \& Touzi \cite{beiglbock2015complete}, it turns out that the martingale transport problem reduces to componentwise irreducible martingale transport problems for which the quasi-sure formulation and the pointwise one are equivalent. For $\P\in\Mc(\mu,\nu)$, we shall denote $\nu^\P_I:=\P\circ (Y|X\in I)^{-1}$ and $\P_I:=\P\circ ((X,Y)|X\in I)^{-1}$.

\begin{Theorem}\label{thm:desintegration}
Let $c:\Omega\to\overline{\R}_+$ be upper semianalytic with $\Sbf_{\mu,\nu}(c)<\infty$. Then we have:
\be\label{eq:supinteg}
 \Sbf_{\mu,\nu}(c)=\underset{\P\in\Mc(\mu,\nu)}{\sup}\int_{I(\R^d)} \Sbf_{\mu_I,\nu^\P_{I}}(c)\eta(dI).
\ee

\no Furthermore, we may find functions $(\varphi,h)\in\Leb^0(\R^d)\x\Leb^0(\R^d,\R^d)$, and $(\psi_K)_{K\in I(\R^d)}\subset\Leb_+^0(\R^d)$ with $\psi_{I(X)}(Y)\in\Leb^0_+(\Omega)$, and $\dom\,\psi_I = J_\theta$, $\eta-$a.s. for some $\theta\in\widehat{\Tc}(\mu,\nu)$, such that

\no{\rm (i)} $c
 \le \bar c :=
\varphi(X)+\psi_{I(X)}(Y)+h^\otimes$, and $
\Sbf_{\mu,\nu}(c)
 =
\Sbf_{\mu,\nu}\big(\bar c\big).
 $
 
\no{\rm (ii)} If the supremum \eqref{eq:supinteg} has an optimizer $\P^*\in\Mc(\mu,\nu)$, then we may chose $\big(\varphi,h,(\psi_K)_K\big)$ so that $(\varphi,\psi_I,h)\in\Dc^{pw}_{\mu_I,\nu^{\P^*}_I}(c_{|I\x J_\theta})$, and
$\Sbf_{\mu_I,\nu^{\P^*}_{I}}(c)
=
\Ibf^{\rm pw}_{\mu_I,\nu^{\P^*}_{I}}(c)
=
\mu_I[\varphi]{\oplus}\nu_I^{\P^*}[\psi_I], ~\eta-\mbox{a.s.}$

\no{\rm (iii)} If Assumption \ref{ass:domination} holds, we may find $J\in\Jc(\mu,\nu)$, and $(\varphi',\psi',h')\in\Dc^{qs}_{\mu,\nu}(c)$ optimizer for $\Ibf_{\mu,\nu}^{qs}(c)$ such that $c\le \varphi'\oplus\psi'+h'^\otimes$, on $\{Y\in J(X)\}$.

\no{\rm (iv)} Under the conditions of {\rm (ii)} and {\rm (iii)}, we may find $(\varphi',\psi',h')\in\Dc^{pw}_{\mu_I,\nu^{\P^*}_I}(c_{|I\x J})$, such that $\Sbf_{\mu_I,\nu^{\P^*}_{I}}(c)
=
\mu_I[\varphi']{\oplus}\nu_I^{\P^*}[\psi'], ~\eta-\mbox{a.s.}
$
\end{Theorem}

Theorem \ref{thm:desintegration} will be proved in Subsection \ref{subsect:proofdisintegretion}

\begin{Remark}
Notice that $(\mu_I,\nu_I^{\P^*})$ may not be irreducible. See Example \ref{expl:decompositionnotirreducible}. This is an important departure from the one-dimensional case.
\end{Remark}

\begin{Remark}
Existence holds for the maximization problem \eqref{eq:supinteg} (and therefore (ii) in Theorem \ref{thm:desintegration} holds) under any of the following assumptions:

\no{\rm (i)} $\nu_I:=\nu_I^\P$ is independent of $\P\in\Mc(\mu,\nu)$ (see Remark \ref{rmk:invariancenu} for some sufficient conditions);

\no{\rm (ii)} There exists a primal optimizer for the problem $\Sbf_{\mu,\nu}(c)$.
\end{Remark}

\subsection{Martingale monotonicity principle}

As a consequence of the last duality result, we now provide the martingale version of the monotonicity principle which extends the corresponding result in standard optimal transport theory, see Theorem 5.10 in Villani \cite{villani2008optimal}. The following monotonicity principle states that the optimality of a martingale measure reduces to a property of the corresponding support.

The one-dimensional martingale monotonicity principle was introduced by Beiglb\"ock \& Juillet \cite{beiglboeck2016problem}, see also Zaev \cite{zaev2015monge}, and Beiglb\"ock, Nutz \& Touzi \cite{beiglbock2015complete}.

\begin{Theorem}\label{thm:monotonicity}
Let $c:\Omega\to\overline{\R}_+$ be upper semianalytic with $\Sbf_{\mu,\nu}(c)<\infty$.
\\
{\rm (i)} Then we may find a Borel set $\Gamma\subset\Omega$ such that:

{\rm (a)} Any solution $\P$ of ${\bf S}_{\mu,\nu}(c)$, is concentrated on $\Gamma$;

{\rm (b)} we may find $\theta\in\widehat{\Tc}(\mu,\nu)$ and $(\Gamma_K)_{K\in I(\R^d)}$ such that $\Gamma=\cup_{K\in I(\R^d)}\Gamma_{K}$ with $\Gamma_{I}\subset I\x J_\theta$, $\Gamma_{I}$ is $c$-martingale monotone, and for any optimizer $\P^*$ of ${\bf S}_{\mu,\nu}(c)$, we have that any optimizer $\P\in\Mc(\mu_{I},\nu_{I}^{\P^*})$ of ${\bf S}_{\mu_{I},\nu_{I}^{\P^*}}(c)$, is concentrated on $\Gamma_{I}$.
\\
\no{\rm (ii)} if Assumption \ref{ass:domination} holds, we may find a universally measurable $\Gamma'\subset N^c$, for some canonical $N\in\Nc_{\mu,\nu}$, satisfying {\rm (a)} and {\rm (b)}, such that $\Gamma'$ is $c$-martingale monotone.
\end{Theorem}
\proof
Let functions $(\varphi,h)\in\Leb^0(\R^d)\x\Leb^0(\R^d,\R^d)$ and functions $(\psi_K)_{K\in I(\R^d)}\subset\Leb_+^0(\R^d)$ with $\psi_{I(X)}(Y)\in\Leb^0_+(\Omega)$ from Theorem \ref{thm:desintegration}. Recall that pointwise we have $c \le \varphi(X)+\psi_{I(X)}(Y)+h^\otimes$. We set $\Gamma:=\{c = \varphi(X)+\psi_{I(X)}(Y)+h^\otimes<\infty\}$.

(i) If $\P^*$ is optimal for the primal problem then,
\b*
\infty>\P^*[c]=\P^*[\varphi(X)+\psi_{I(X)}(Y)+h^\otimes] = \Sbf_{\mu,\nu}(c)
&and&
\P^*[\varphi(X)+\psi_{I(X)}(Y)+h^\otimes-c] = 0
\e*
As $\varphi(X)+\psi_{I(X)}(Y)+h^\otimes-c\geq 0$, and the expectation of $c$ is finite, and therefore $\P^*[c<\infty]=1$, it follows that $\P^*$ is concentrated on $\Gamma$.

(ii) Let $\theta\in\widehat{\Tc}(\mu,\nu)$ such that $J_\theta = \dom \psi_{I(X)}$ from Theorem \ref{thm:desintegration}. For $K\in I(\R^d)$, let $\Gamma_K:=\Gamma\cap K\x\R^d$. Then we have $\Gamma_{I(x)}\subset I(x)\x J_\theta(x)$ for all $x\in\R^d$, $\Gamma_{I(x)}$ is $c$-martingale monotone because of the pointwise duality on each component, and $\Gamma=\cup_{x\in \R^d}\Gamma_{I(x)}$ by definition because $I(\R^d)$ is a partition of $\R^d$.

If Assumption \ref{ass:domination} holds, we consider $(\varphi',\psi', h')\in\Dc^{qs}_{\mu,\nu}(c)$ from the second part of Theorem \ref{thm:desintegration}. Let a canonical $N\in\Nc_{\mu,\nu}$ be such that $c=\varphi'\oplus\psi'+\bar h'^\otimes$ on $N^c$. $\Gamma := N^c\cap\{c=\varphi'\oplus\psi'+\bar h'^\otimes\}$. Similarly, (i) and (ii) hold.

(iii) By definition of $\Theta_{\mu,\nu}$, for $\P_0$ with finite support, supported on $\Gamma\subset N^c$, and $\P'$ competitor to $\P_0$. As $N^c$ is canonical, it is martingale monotone by definition. Then $\P'[N^c]=1$, and therefore $\P'[c]\le\P'[\varphi'\oplus\psi'+h'^\otimes]=\P[\varphi'\oplus\psi'+h'^\otimes]= \P_0[c]$.

Finally, by definition we have $\Gamma\subset N^c$.
\ep

\begin{Remark}\label{rmk:probainvariant}
Let $(\varphi,\psi,h)\in\Dc^{qs}_{\mu,\nu}(c)$ be a minimizer of $\Ibf_{\mu,\nu}^{q.s.}(c)$. Assume that $\P[\varphi\oplus\psi+h^\otimes]$ does not depend on the choice of $\P\in\Mc(\mu,\nu)$ (e.g. if $(\varphi,\psi)\in\L^1(\mu)\x\L^1(\nu)$, or if $d=1$). Then we may chose $\Gamma$ such that a measure $\P\in\Mc(\mu,\nu)$ is optimal for ${\bf S}_{\mu,\nu}(c)$ if and only if it is concentrated on $\Gamma$. Indeed, with the notations from the previous proof, if $\P\in\Mc(\mu,\nu)$ is concentrated on $\Gamma$, $\P[\varphi\oplus\psi+h^\otimes-c] = 0$ and as $\P[\varphi\oplus\psi+h^\otimes]=\mu[\varphi]{\oplus}\nu[\psi]$ because of the invariance,
$$\P(c) = \P[\varphi\oplus\psi+h^\otimes] = \Ibf^{qs}_{\mu,\nu}(c) = \Sbf_{\mu,\nu}(c).$$
\end{Remark}

\subsection{On Assumption \ref{ass:domination}}
\label{sec:assumptions}

\begin{Proposition}\label{prop:emptyboundaries}
Assumption \ref{ass:domination} holds true under either one of the following conditions:

\no{\rm (i)} $Y\notin \partial I(X)$, $\Mc(\mu,\nu)$-q.s. or equivalently $\mu\circ I^{-1} = \nu\circ I^{-1}$.

\no{\rm (ii)} $\dim I(X)\in\{0,1,d\}$, $\mu-$a.s.

\no{\rm (iii)} $\nu$ is dominated by the Lebesgue measure and $\dim I(X)\in\{0,1,d-1,d\}$, $\mu-$a.s.

\no{\rm (iv)}  $I(X)\in\Cc\cup\Dc\cup\Rc$, $\mu-$a.s. for some subsets $\Cc,\Dc,\Rc\subset \Kcirc$ with $\Cc$ countable, $\dim(\Dc)\subset\{0,1\}$, and $\cup_{K\in\Rc}K\x \partial K\in\Nc_{\mu,\nu}$.
\\
Furthermore, {\rm (iv)} is implied by either one of {\rm (i)}, {\rm (ii)}, and {\rm (iii)}.
\end{Proposition}

This proposition is proved in Subsection \ref{subsect:verifass}.

\begin{Remark}
Assumption \ref{ass:domination} holds in dimension $1$ by Proposition \ref{prop:emptyboundaries}. Theorem \ref{thm:duality} is equivalent to \cite{beiglbock2015complete} Theorem 7.4 and the monotonicity principle Theorem \ref{thm:monotonicity} is equivalent to \cite{beiglbock2015complete} Corollary 7.8.
\end{Remark}

\begin{Remark}\label{rmk:invariancenu}
Notice that under either one of (i) or (iii) of Proposition \ref{prop:emptyboundaries}, or in dimension one, the disintegration $\nu^\P_I:=\P\circ (Y|X\in I)^{-1}$ is independent of the choice of $\P\in\Mc(\mu,\nu)$. See Subsection \ref{subsect:verifass} for a justification of this claim.
\end{Remark}

\begin{Remark}\label{rmk:densityofirr}
Proposition \ref{prop:emptyboundaries} may be applied in particular in the trivial case in which there is a unique irreducible component. We state here that any pair of measures $\mu,\nu\in\P(\R^d)$ in convex order may be approximated by pairs of measures that have a unique irreducible component, and therefore satisfy Assumption \ref{ass:domination}. We may then use a stability result like in Guo \& Ob{\l}{\'o}j \cite{guo2017computational} to use the approximation $(\mu_\epsilon,\nu_\epsilon)$ of $(\mu,\nu)$ in practice.

Let $\mu'\preceq\nu'$ in convex order with $(\mu',\nu')$ irreducible, and $\supp\,\nu\subset \ri\,\conv\,\supp\,\nu'$. Then $(\mu_\epsilon,\nu_\epsilon):=\frac{1}{1+\epsilon}(\mu+\epsilon\mu',\nu+\epsilon\nu')$ is irreducible for all $\epsilon>0$. Indeed by Proposition 3.4 in \cite{de2017irreducible}, we may find $\hat\P\in\Mc(\mu',\nu')$ such that $\conv\,\supp\,\hat\P_X=\ri\,\conv\,\supp\,\nu'$, $\mu'-$a.s. Then, $\frac{1}{1+\epsilon}(\P+\eps\hat\P)\in \Mc(\mu_\epsilon,\nu_\epsilon)$ for all $\P\in\Mc(\mu,\nu)$, and $\ri\,\conv\,\supp\,\nu'\subset I(X)$ on a set charged by $\mu_\epsilon$, which proves that $I=\ri\,\conv\,\supp\,\nu'\supset \supp\,\nu$, preventing other components from appearing on the boundary. Thus $(\mu_\epsilon,\nu_\epsilon)$ is irreducible.

Convenient measures to consider are for example $\mu':=\delta_0$ or $\mu':=\Nc(0,1)$, and $\nu':=\Nc(0,2)$. For finitely supported $\mu$ and $\nu$ we may consider $y_1,...,y_k\in \R^d$ for some $k\ge 1$ such that $\supp\,\nu\subset \interior\,\conv(y_1,...,y_k)$, $\nu':=\frac{\delta_{y_1}+...+\delta_{y_k}}{n}$, and $\mu':=\delta_{\frac{y_1+...+y_n}{n}}$.
\end{Remark}

\begin{Proposition}\label{prop:axioms}
Assumption \ref{ass:domination} holds if we assume existence of medial limits and Axiom of choice for $\R$.
\end{Proposition}

We prove this Proposition in Subsection \ref{subsect:verifassaxiom}.

\begin{Remark}
Notice that existence of medial limits and Axiom of choice for $\R$ is implied by Martin's axiom and Axiom of choice for $\R$, which is implied by the continuum hypothesis. Furthermore, all these axiom groups are undecidable under either the Theory ZF nor the Theory ZFC. See Subsection \ref{subsect:verifassaxiom}.
\end{Remark}

\subsection{Measurability and regularity of the dual functions}

In the main theorem, only $\varphi\oplus\psi+h^\otimes$ has some measurability. However, we may get some measurability on the separated dual optimizers.

\begin{Proposition}\label{prop:measurability}
For all $(\varphi,\psi,h)\in\widehat\Leb(\mu,\nu)$,
\\
{\rm (i)} $\big(\varphi,\psi,proj_{\nabla\aff I}(h)\big)\in \Leb^0(I)\x\Leb^0(I)\x\Leb^0(I,\nabla\aff I)$;
\\
{\rm (ii)} under any one of the conditions of Proposition \ref{prop:emptyboundaries}, we may find $(\varphi',\psi',h')\in\widehat\Leb(\mu,\nu)$ such that $\varphi\oplus\psi + h^\otimes=\varphi'\oplus\psi' + h'^\otimes$, q.s. and $(\varphi',\psi',h')\in \Leb^0\big(\R^d\big)^2\x\Leb^0\big(\R^d,\R^d\big)$. Furthermore, the canonical set from Theorem \ref{thm:polar}, and the set $\Gamma'$ from Theorem \ref{thm:monotonicity} may be chosen to be Borel measurable, and $\{Y\in J(X)\}$ (resp. $\{Y\in \Jo(X)\}$) for $J\in{\Jc}(\mu,\nu)$ (resp. $\Jo\in{\Jco}(\mu,\nu)$) may be chosen to be analytically measurable.
\end{Proposition}

The proof of this proposition is reported to Subsection \ref{subsect:CctoTc}. We may as well prove some regularity of the dual functions, provided that the cost function has some appropriate regularity. This Lemma is very close to Theorem 2.3 (1) in \cite{ghoussoub2015structure}.

\begin{Lemma}\label{lemma:regularity}
Let $c:\Omega\to\overline{\R}_+$ be upper semi-analytic. We assume that $x\longmapsto c(x,y)$ is locally Lipschitz in $x$, uniformly in $y$, and that $\Sbf_{\mu,\nu}(c)=\Sbf_{\mu,\nu}(\varphi\oplus\psi +h^\otimes)<\infty$, with $\varphi:\R^{d}\longmapsto \R\cup\{\infty\}$, $\psi:\R^{d}\longmapsto \R\cup\{\infty\}$, and $h:\R^{d}\longmapsto \R^{d}$ such that $c\le \varphi\oplus\psi +h^\otimes$, pointwise. Then, we may find $(\varphi',h') = (\varphi,h)$, $\mu-a.e.$ such that $c\le \varphi'\oplus\psi+h'^\otimes\le \varphi\oplus\psi+h'^\otimes$, $\varphi'$ is locally Lipschitz, and $h'$ is locally bounded on $\ri\,\conv\,\dom\,\psi$.
\end{Lemma}

The proof of Lemma \ref{lemma:regularity} is reported in Subsection \ref{subsect:regularity}.

\section{Examples}\label{sect:examples}
\setcounter{equation}{0}

\subsection{Pointwise duality failing in higher dimension}

In the one-dimensional case, Beiglb\"ock, Lim \& Ob{\l}{\'o}j \cite{beiglbock2017dual} proved that pointwise duality, and integrability hold for $\Ctn^2$ cost functions together with compactly supported $\mu$, and $\nu$. We believe that integrability may hold in higher dimension, and strong monotonicity holds. However the following example shows that dual attainability does not hold with such generality for a dimension higher than 2.

\begin{Example}\label{expl:nopwdualitydim2}
Let $y_{--}:=(-1,-1)$, $y_{-+}:=(-1,1)$, $y_{+-}:=(1,-1)$, $y_{++}:=(1,1)$, $y_{0-}:=(0,-1)$, $y_{0+}:=(0,1)$, $y_{00}:=(0,0)$, $y_{+0}:=(1,0)$, $C:=\conv(y_{--},y_{-+},y_{+-},y_{++})$, $x_1:= (-\frac12,0)$, $x_2:= (\frac12,\frac12)$, $x_3:= (\frac12,-\frac12)$, $\mu := \frac12\delta_{x_1}+\frac14\delta_{x_2}+\frac14\delta_{x_3}$, and $\nu:=\frac14 \mathbf{1}_{C} \mbox{Vol}$. We can prove that for these marginals, the irreducible components are given by
\b*
I(x_1):=\ri\,\conv(y_{--},y_{-+},y_{0+},y_{0-}),\quad I(x_2):=\ri\,\conv(y_{0+},y_{++},y_{+0},y_{00}),\\
\mbox{and}\quad I(x_3):=\ri\,\conv(y_{00},y_{+0},y_{+-},y_{0-}),
\e*
and $\Mc(\mu,\nu)$ is a singleton $\{\P\}$, with
$$\P(dx,dy):=\frac14\big(2\delta_{x_1}(dx)\mathbf{1}_{y\in I(x_1)}+\delta_{x_2}(dx)\mathbf{1}_{y\in I(x_2)}+\delta_{x_3}(dx)\mathbf{1}_{y\in I(x_3)}\big)\otimes \mbox{Vol}(dy).$$
 Now we define a cost function $c$ such that $c\big(x_1,\cdot\big)$ is $0$ on $\cl I(x_1)$, $c\big(x_2,\cdot\big)$ is $0$ on $\cl I(x_2)$, and $c\big(x_3,\cdot\big)$ is $0$ on $\cl I(x_3)$. However we also require $c(x_2,y_{+-})=1$. We may have these conditions satisfied with $c\ge 0$, and $\Ctn^\infty$. Let $(\varphi,\psi,h)$ be pointwise dual optimizers, then $\varphi\oplus\psi+h^\otimes=c$, $\P-$a.s. then $\psi$ is affine on each irreducible components: $\psi(y) = c(x_i,y)-\varphi(x_i)-h(x_1)\cdot(y-x_i)=-\varphi(x_i)-h(x_1)\cdot(y-x_i)$, Lebesgue-a.e. on $I(x_i)$, for $i=1,2,3$. By the last equality, we deduce that $\varphi(x_i)=-\psi(x_i)$, and $h(x_i)=-\nabla\psi(x_i)$. Now by the superhedging inequality, $\psi(y)-\psi(x_i)-\nabla\psi(x_i)\cdot(y-x_i) \ge c(x_i,y) \ge 0$. Therefore $\psi$ is a.e. equal to a convex function, piecewise affine on the components. However a convex function that is affine on $I(x_1)$, $I(x_2)$, and $I(x_3)$ is affine on $\cl I(x_2)\cup \cl I(x_3)$ (it follows from the verification at the angles between the regions where $\psi$ has nonzero curvature). Then $c(x_2,y)\le \psi(y)-\psi(x_2)-\nabla\psi(x_2) = 0$ for a.e. $y\in\cl I(x_3) \subset \cl I(x_2)\cup \cl I(x_3)$. This is the required contradiction as $c(x_2,y_{+-})=1$ and $c$ is continuous, and therefore nonzero on a non-negligible neighborhood of $(x_2,y_{+-})$.

Notice that in this example, $\mu$ is not dominated by the Lebesgue measure for simplicity, however this example also holds when $\delta_{x_i}$ is replaced by $\frac{1}{\pi\epsilon^2}\mathbf{1}_{B_\epsilon(x_i)}\mbox{Vol}$ for $\epsilon>0$ small enough.
\end{Example}

\subsection{Disintegration on an irreducible component is not irreducible}

\begin{figure}[h]
\centering

 \includegraphics[width=0.8\linewidth]{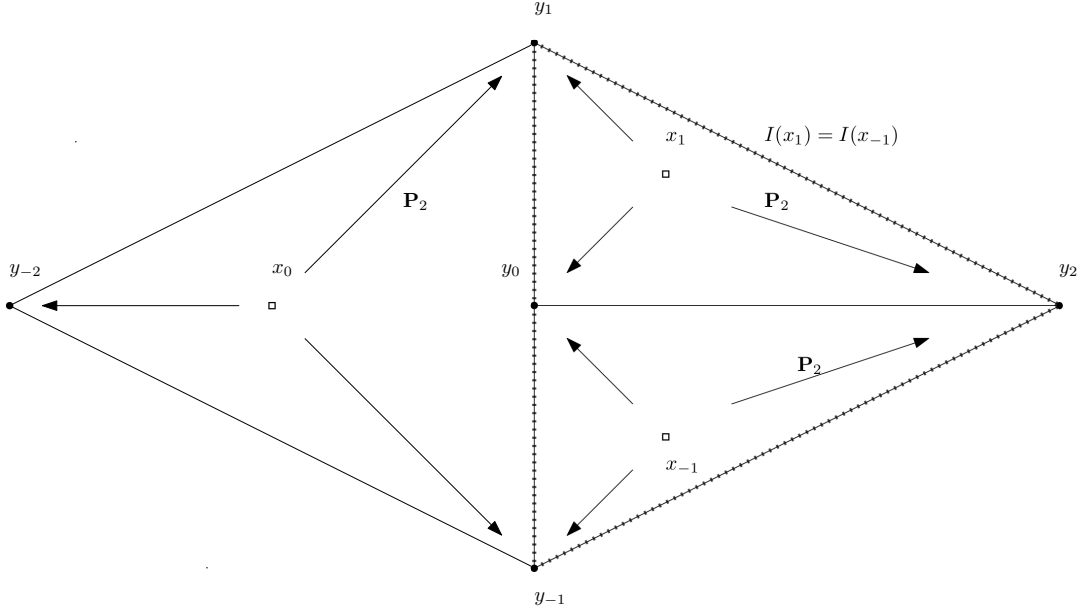}
    \caption{\label{fig:disinnotirr} Disintegration on an irreducible component is not irreducible.}
\end{figure}

\begin{Example}\label{expl:decompositionnotirreducible}
Let $x_{0}:= (-1,0)$, $x_1:= (\frac12,\frac12)$, $x_{-1}:= (\frac12,-\frac12)$, $y_{1}=(0,1)$, $y_2:=(2,0)$, $y_{-1}:= - y_1$, $y_{-2}:= - y_2$, and $y_0:=0$. Let the probabilities
\b*
\mu := \frac13(\delta_{x_0}+\delta_{x_{1}}+\delta_{x_{-1}}),&\mbox{and}&\nu:=\frac16(\delta_{y_{-2}}+ \delta_{y_{2}}+ \delta_{y_{0}})+\frac14(\delta_{y_1}+ \delta_{y_{-1}}).
\e*
We can prove that for these marginals, the irreducible components are given by
\b*
I(x_0)=\ri\,\conv(y_{-2},y_{1},y_{-1}),& \mbox{and }I(x_1)=I(x_{-1})=\ri\,\conv(y_{2},y_{1},y_{-1}),
\e*
indeed, $\Mc(\mu,\nu) = \conv(\P_1,\P_2)$, with
\b*
\P_1:= \frac16\delta_{(x_0,y_{-2})}+\frac16\delta_{(x_0,y_{0})}
&+&\frac{1}{12}\delta_{(x_1,y_2)}+\frac{3}{16}\delta_{(x_1,y_1)}+\frac{1}{16}\delta_{(x_1,y_{-1})}\\
&+&\frac{1}{12}\delta_{(x_{-1},y_2)}+\frac{3}{16}\delta_{(x_{-1},y_{-1})}+\frac{1}{16}\delta_{(x_{-1},y_1)},
\e*
and
\b*
\P_2:=\frac16\delta_{(x_0,y_{-2})}+\frac{1}{12}\delta_{(x_0,y_1)}+\frac{1}{12}\delta_{(x_0,y_{-1})}&+&\frac{1}{6}\delta_{(x_1,y_1)}+\frac{1}{12}\delta_{(x_1,y_{0})}+\frac{1}{12}\delta_{(x_1,y_2)}\\
&+&\frac{1}{6}\delta_{(x_{-1},y_{-1})}+\frac{1}{12}\delta_{(x_{-1},y_{0})}+\frac{1}{12}\delta_{(x_{-1},y_2)}.
\e*
(See Figure \ref{fig:disinnotirr}). Let $c$ be smooth, equal to $1$ in the neighborhood of $(x_0,y_1)$ and $0$ at a distance higher than $\frac12$ from this point, $\P_2$ is the only optimizer for the martingale optimal transport problem $\Sbf_{\mu,\nu}(c)$. However, $\mu_{I(x_1)}= \frac12(\delta_{x_1}+\delta_{x_{-1}})$, and $\nu^{\P_2}_{I(x_1)} = \frac14(\delta_{y_2}+\delta_{y_0}+\delta_{y_1}+\delta_{y_{-1}})$, and the associated irreducible components are
\b*
I_{\mu_{I(x_1)},\nu^{\P_2}_{I(x_1)}}(x_1) = \ri\,\conv(y_0,y_1,y_2),&\mbox{and}& I_{\mu_{I(x_1)},\nu^{\P_2}_{I(x_1)}}(x_{-1}) = \ri\,\conv(y_0,y_{-1},y_2),
\e*
and therefore, the couple $\left(\mu_{I(x_1)},\nu^{\P_2}_{I(x_1)}\right)$ obtained from the disintegration of the optimal probability $\P_2$ in the irreducible component $I(x_1)=I_1$ can be reduced again in two irreducible sub-components.
\end{Example}

 \subsection{Coupling by elliptic diffusion}
 
Assumption \ref{ass:domination} holds when $\nu$ is obtained from an Elliptic diffusion from $\mu$.

\begin{Remark}\label{rmk:diffusion}
Notice that {\rm (iii)} in Proposition \ref{prop:emptyboundaries} holds if $\nu$ is the law of $X_\tau:=X_0+\int_0^t \sigma_sdW_s$, where $X_0\sim \mu$, $W$ a $d-$dimensional Brownian motion independent of $X_0$, $\tau$ is a positive bounded stopping time, and $(\sigma_t)_{t\ge 0}$ is a bounded cadlag process with values in $\Mc_{d}(\R)$ adapted to the $W-$filtration with $\sigma_0$ invertible. We observe that the strict positivity of the stopping time is essential, see Example \ref{expl:zerostoppingtime}.
\end{Remark}

We justify Remark \ref{rmk:diffusion} in Subsection \ref{subsect:verifass}.

 \begin{Example}\label{expl:zerostoppingtime} Let $C:=[-1,1]\x[0,2]\x[-1,1]$, $F:=\{0\}\x[-1,1]\x[-1,1]$, $x_0:=(0,0,0)$, $x_1:=(0,1,0)$ $\mu:=\frac12\delta_{x_0} + \frac12\delta_{x_1}$, a $\Fc-$Brownian motion $W$, and $X$ a random variable $\Fc_0-$measurable with $X_0\sim\mu$. Consider the bounded stopping time $\tau:=1\wedge\inf\{t\ge 0:W_t\in\partial C\}$, and $\nu$, the law of $X_0+W_\tau$. We have $\mu\preceq\nu$ in convex order, as the law $\P$ of $(X,Y):=(X_0,X_0+W_\tau)$ is clearly a martingale coupling. However, observe that $p:=\P[X=x_1,Y\in C]>0$, and that by symmetry $\P[Y|X=x_1,Y\in C] = x_0$. Let $\nu_C$ be the law of $Y$, conditioned on $\{X=x_1,Y\in C\}$. Then $\P':= \P+p\big((\delta_{x_0}-\delta{x_1})\otimes\nu_C-(\delta_{x_0}-\delta{x_1})\otimes\delta_{x_0}\big)$ is also in $\Mc(\mu,\nu)$. We may prove that the irreducible components are $\ri C$, and $\ri F$, and therefore (iii) of Proposition \ref{prop:emptyboundaries} does not hold. This proves the importance of the strict positivity of the stopping time $\tau$ in Remark \ref{rmk:diffusion}. In dimension $4$, we may find an example in which (v) of Proposition \ref{prop:emptyboundaries} does not hold either, by replacing $F$ by a continuum of translated $F$ in the fourth variable, thus introducing an orthogonal curvature in the lower face of $C$ to avoid the copies of $F$ to communicate with each other. 
 \end{Example}

\section{Proof of the main results}\label{sect:proofs}
\setcounter{equation}{0}

\subsection{Moderated duality}

Let $c\ge 0$, we define the moderated dual set of $c$ by
\b*
\Dc_{\mu,\nu}^{\widetilde{mod}}(c)&:=\Big\{&(\bar\varphi,\bar\psi,\bar h,\theta)\in\Leb^1_+(\mu)\x\Leb^1_+(\nu)\x\Leb^0(\R^d,\R^d)\x\widetilde{\Tc}(\mu,\nu):\\
&&c\le\bar\varphi\oplus\bar\psi+\bar h^\otimes+\theta, \mbox{ on }\{Y\in\aff\,\rf_X\conv\,\dom(\theta+\bar\psi)\}\Big\}.
\e*
We then define for $(\bar\varphi,\bar\psi,\bar h,\theta)\in\Dc_{\mu,\nu}^{\widetilde{mod}}(c)$, $Val(\bar\varphi,\bar\psi,\bar h,\theta):=\mu[\bar\varphi]+\nu[\bar\psi]+\nu\widetilde{\ominus}\nu[\theta]$, and the moderated dual problem $\Ibf_{\mu,\nu}^{\widetilde{mod}}(c):=\inf_{\xi\in\Dc_{\mu,\nu}^{\widetilde{mod}}(c)}Val(\xi)$.

\begin{Theorem}\label{thm:dualitymod}
Let $c:\Omega\to\overline{\R}_+$ be upper semianalytic. Then, under Assumption \ref{ass:domination}, we have
\\
{\rm (i)} $\Sbf_{\mu,\nu}(c) = \Ibf^{\widetilde{mod}}_{\mu,\nu}(c)$;
\\
{\rm (ii)} If in addition $\Sbf_{\mu,\nu}(c)<\infty$, then existence holds for the moderated dual problem $\Ibf^{\widetilde{mod}}_{\mu,\nu}(c)$.
\end{Theorem}

This Theorem will be proved in Subsection \ref{subsect:proofduality}.

\subsection{Definitions}\label{subsect:reminders}

We first need to recall some concepts from \cite{de2017irreducible}. For a subset $A\subset\R^d$ and $a\in\R^d$, we introduce the face of $A$ relative to $a$ (also denoted $a-$relative face of $A$): $
 \rf_a A  
 := 
 \big\{y\in A: (a-\eps(y-a),y+\eps(y-a)) \subset A,\text{ for some }\eps>0\big\}
 $. Now denote for all $\theta:\Omega\to\bar\R$: 
 \b*
 \dom_x\theta
 &:=&
 \rf_x\conv\,\dom\,\theta(x,\cdot).
 \e*
For $\theta_1,\theta_2:\Omega\longrightarrow\R$, we say that $\theta_1=\theta_2$, $\muxpw$, if 
 \b*
 \dom_X\theta_1 = \dom_X\theta_2,
 &\mbox{and}&
 \theta_1(X,\cdot)=\theta_2(X,\cdot)
 ~\mbox{on}~\dom_X\theta_1,~\mu-\mbox{a.s.}
 \e*

The main ingredient for our extension is the following.

\begin{Definition}
A measurable function $\theta:\Omega\to\overline{\R}_+$ is a tangent convex function if 
 \b*
 &\theta(x,\cdot)
  ~\mbox{is convex, and}~
  \theta(x,x)=0,
 ~\mbox{for all}~
 x\in\R^d.
 \e*
We denote by $\Theta$ the set of tangent convex functions, and we define
 \b*
 \Theta_\mu
 &:=&
 \big\{\theta\in\L^0(\Omega,\overline{\R}_+):
       \theta = \theta',~\muxpw,
       \mbox{and}~\theta\ge\theta',~\mbox{for some}~
       \theta'\in\Theta   
 \big\}.
 \e*
\end{Definition}

\begin{Definition}
A sequence $(\theta_n)_{n\ge 1}\subset\L^0(\Omega)$ converges $\muxpw$ to some $\theta\in\L^0(\Omega)$ if 
 \b*
 \dom_X\left(\underline\theta_\infty\right) =  \dom_X\theta &\mbox{and}&\theta_n(X,\cdot) \longrightarrow  \theta(X,\cdot),\,\mbox{pointwise on}~\dom_X\theta,
 ~\mu-\mbox{a.s.}
 \e*

{\rm (i)} A subset $\Tc\subset\Theta_\mu$ is $\muxpw$-Fatou closed if $\underline\theta_\infty\in\Tc$ for all $(\theta_n)_{n\ge 1}\subset\Tc$ converging $\muxpw$.

\no{\rm (ii)} The $\muxpw-$Fatou closure of a subset $A\subset\Theta_\mu$ is the smallest $\muxpw-$Fatou closed set containing $A$:
 $$
 \widehat A
 :=
 \bigcap\big\{\Tc\subset\Theta_\mu:~
              A\subset\Tc,
              ~\text{and}~
              \Tc~\mbox{$\muxpw$-Fatou closed}\,
         \big\}.
 $$
\end{Definition}

Recall the definition for $a\ge 0$, of the set $\Cfrak_a:=\big\{f\in\Cfrak:(\nu-\mu)(f)\le a\big\}$, we introduce
$$\widehat{\Tc}(\mu,\nu):=\underset{a\ge 0}{\bigcup}\,\widehat{\Tc}_a,
\text{ where }
\widehat{\Tc}_a:=\widehat{\Tbf(\Cfrak_a)}, ~\mbox{and}~
 \Tbf \big(\Cfrak_a\big)
 :=
 \big\{ \Tbf_p f: f\in\Cfrak_a,p\in\partial f\big\}.
 $$
 
Similar to $\nu\overline{\ominus}\mu$ for $\widetilde{\Tc}(\mu,\nu)$, we now introduce the extended $(\nu-\mu)-$integral:
 \b*
 \nu\widehat{\ominus}\mu[\theta]
 :=
 \inf\big\{a\ge 0 :\theta\in \widehat{\Tc}_a\big\} 
 &\mbox{for}&
 \theta\in\widehat{\Tc}(\mu,\nu).
 \e*

\subsection{Duality result}\label{subsect:proofduality}

As a preparation for the proof of Theorem \ref{thm:dualitymod}, we prove the following Lemma.

\begin{Lemma}\label{lemma:existencehat}
Let $\widehat\theta\in\widehat{\Tc}(\mu,\nu)$, under Assumption \ref{ass:domination}, we may find $\theta\in\widetilde{\Tc}(\mu,\nu)$ such that $\theta\ge\widehat\theta$ and $\nu\overline\ominus{}\mu[\theta]\le\nu\widehat{\ominus}\mu[\widehat\theta]$.
\end{Lemma}
\proof
Let $a> 0$, we consider $\Tc$ the collection of $\widehat\theta\in\Theta_\mu$ such that we may find $\theta\in\widetilde{\Tc}_a$ with $\theta\ge\widehat\theta$. First we have easily $\Tbf(\Cfrak_a)\subset\Tc$, as $\Tbf(\Cfrak_a)\subset\widetilde{\Tc}_a$. Now we consider $(\widehat\theta_n)_{n\ge 1}\subset \Tc$ converging $\muxpw$ to $\underline{\widehat\theta}_\infty$. For each $n\ge 1$, we may find $\theta_n\in\widetilde{\Tc}_a$ such that $\theta_n\ge\widehat\theta_n$ and $\nu\overline\ominus{}\mu[\theta_n]\le a$. Now we may use Assumption \ref{ass:domination}, we may find $\theta\in\widetilde{\Tc}_a$ such that $\theta_n\rightsquigarrow\theta$ by the fact that $\widetilde\Tc_a = a\widetilde\Tc_1$. By the generation properties, $\theta\ge \underline\theta_\infty\ge \underline{\widehat\theta}_\infty$, which implies that $\underline{\widehat\theta}_\infty\in\Tc$. $\Tc$ is $\muxpw-$Fatou closed, and therefore $\widehat\Tc_a\subset\Tc$.

Now let $\widehat\theta\in\widehat{\Tc}(\mu,\nu)$, with $l:=\nu\widehat{\ominus}\mu[\widehat\theta]$. By what we did above, for all $n\ge 1$, we may find $\theta_n\in\widetilde{\Tc}_{l+1/n}$ such that $\widehat\theta\le\theta_n$. We use again Assumption \ref{ass:domination} to get $\theta_n\rightsquigarrow \theta$, by properties of generation, $\theta\ge \underline\theta_\infty\ge\widehat\theta$. By construction, $\nu\overline\ominus{}\mu[\theta]\le l=\nu\widehat{\ominus}\mu[\widehat\theta]$.
\ep\\

\no {\bf Proof of Theorem \ref{thm:dualitymod} }
By Theorem 3.8 in \cite{de2017irreducible}, we may find $(\bar\varphi,\bar\psi,\bar h,\widehat\theta)\in\Leb_+^1(\mu)\x\Leb_+^1(\nu)\x\Leb^0(\R^d,\R^d)\x\widehat{\Tc}(\mu,\nu)$ such that $c\le \bar\varphi\oplus\bar\psi+\bar h^\otimes+\widehat\theta$ on $\{Y\in\aff\,\rf_X\conv\,\dom(\widehat\theta(X,\cdot)+\bar\psi)\}$, furthermore, $\Sbf_{\mu,\nu}(c) = \mu[\bar\varphi]+\nu[\bar\psi]+\nu\widehat{\ominus}\mu[\widehat\theta]$ and $\Sbf_{\mu,\nu}(\widehat\theta)=\nu\widehat{\ominus}\mu[\widehat\theta]<\infty$. By lemma \ref{lemma:existencehat}, we may find $\theta\in\widetilde{\Tc}(\mu,\nu)$ such that $\widehat\theta\le \theta$ and $\nu\overline\ominus{}\mu[\theta]\le\nu\widehat{\ominus}\mu[\widehat\theta]$.

We have that $c\le \bar\varphi\oplus\bar\psi+\bar h^\otimes+\widehat\theta\le \bar\varphi\oplus\bar\psi+\bar h^\otimes+\theta$ on $\{Y\in\aff\,\rf_X\conv\,\dom(\theta(X,\cdot)+\bar\psi)\}$ which is included in $\{Y\in\aff\, \rf_X\conv\,\dom(\widehat\theta(X,\cdot)+\bar\psi)\}$.

As $\widehat\theta\le \theta$, we have $\Sbf_{\mu,\nu}(\theta)\ge \Sbf_{\mu,\nu}(\widehat\theta)=\nu\widehat{\ominus}\mu[\widehat\theta]\ge\nu\overline\ominus{}\mu[\theta]$. From Proposition \ref{prop:numoinsmuqs} (i), we get that $\Sbf_{\mu,\nu}(\theta)  
 =
 \nu\overline\ominus{}\mu[\theta] = \nu\widehat{\ominus}\mu[\widehat\theta]<\infty$. As $\theta\ge \widehat\theta$, we have $(\bar \varphi,\bar \psi,\bar h,\theta)\in\Dc^{\widetilde{mod}}_{\mu,\nu}(c)$. Finally, as $Val(\bar \varphi,\bar \psi,\bar h,\theta)=\mu[\bar\varphi]+\nu[\bar\psi]+\nu\overline{\ominus}\mu{}[\theta]=\mu[\bar\varphi]+\nu[\bar\psi]+\nu\widehat{\ominus}\mu[\widehat\theta]=\Sbf_{\mu,\nu}(c)$, the result is proved.
\ep\\

\no {\bf Proof of Theorem \ref{thm:duality} }
By Theorem \ref{thm:dualitymod}, we may find $(\bar \varphi,\bar \psi,\bar h,\theta)\in\Dc^{\widetilde{mod}}_{\mu,\nu}(c)$ such that $\mu[\bar\varphi]+\nu[\bar\psi]+\nu\widehat{\ominus}\mu[\theta]=\Sbf_{\mu,\nu}(c)$. As Assumption \ref{ass:domination} holds, by Proposition \ref{prop:CctoTc}, we get $f\in\Cfrak_{\mu,\nu}$ and $p\in\partial^{\mu,\nu}f$ such that $\Tbf_pf=\theta$, q.s. Therefore, by definition we have $\nu\overline{\ominus}\mu[f]\le \nu\widehat{\ominus}\mu[\theta]$. Then we denote $\varphi:=\bar\varphi-f$, $\psi:=\bar\psi+f$, and $h:=\bar h -p$. As $\varphi\oplus\psi+h^\otimes= \bar\varphi\oplus\bar\psi+\bar h^\otimes+\theta\ge c$, q.s., (as $Y\in \aff \,\rf_X\conv\,\dom(\widehat\theta(X,\cdot)+\bar\psi)$, q.s.) $\Sbf_{\mu,\nu}(\Tbf_pf) = \nu\widehat{\ominus}\mu[\theta]\ge \nu\overline{\ominus}\mu[f]$. As $\nu\underline{\ominus}\mu[f]:=\Sbf_{\mu,\nu}(\Tbf_pf)\le \nu\overline{\ominus}\mu[f]$ by Proposition \ref{prop:numoinsmuqs} (i), we have $\nu{\ominus}\mu[f] := \nu\underline{\ominus}\mu[f]= \nu\overline{\ominus}\mu[f]$, and therefore $f$ is a $\Mc(\mu,\nu)-$convex moderator for $(\varphi,\psi)$, and as $\mu[\varphi+f] + \mu[\psi-f] + \nu{\ominus}\mu[f] = \Sbf_{\mu,\nu}(c)$, the duality result, and attainment are proved.
\ep

\subsection{Structure of polar sets}\label{subsect:polarproofs}

\no{\bf Proof of Proposition \ref{prop:J}}
\no{\underline{Step 1:}} Let a Borel $N\in\Nc_{\mu,\nu}$ such that $\theta$ is a $N-$tangent convex function. Then $c:=\infty\mathbf{1}_N$ is Borel measurable and non-negative. Notice that $\Sbf_{\mu,\nu}(c)=0$. By Theorem \ref{thm:dualitymod}, we may find $(\varphi_1,\psi_1, h_1,\theta_1)\in\Dc^{\widetilde{mod}}_{\mu,\nu}(c)$ such that $\mu[\varphi_1]+\nu[\psi_1]+\nu\widehat{\ominus}\mu[\theta_1]=\Sbf_{\mu,\nu}(c) = 0$. Then by the pointwise inequality $\infty\mathbf{1}_N\le \varphi_1\oplus\psi_1+h_1^\otimes+\theta_1$ on $\{Y\in\aff\,\rf_X\conv\,D(X)\}$, with $D(X):=\dom\big(\theta_1(X,\cdot)+\psi_1\big)$, (the convention is $0\times\infty = 0$).

By Subsection 6.1 in \cite{de2017irreducible}, we may find $N_\mu'\in\Nc_\mu$, $N_\nu\in\Nc_\nu$, and $\widehat{\theta}\in\widehat{\Tc}(\mu,\nu)$ such that $I(X)\subset D(X)$, $\rf_X\conv(\dom\widehat{\theta}(X,\cdot)\setminus N_\nu) = I(X)$, and $\dom\widehat{\theta}(X,\cdot)\setminus N_\nu\subset \bar J(X)$, on $N_\mu'^c$. By Lemma \ref{lemma:existencehat} we may find $\widehat{\theta}\le\widetilde\theta\in\widetilde{\Tc}(\mu,\nu)$. Up to adding $\mathbf{1}_{N_\mu'}$ to $\phi_1$, $\mathbf{1}_{N_\nu}$ to $\psi_1$, and $\widetilde{\theta}$ to $\theta_1$, we may assume that $\mathbf{1}_{N_\mu'}\le\phi_1$, $\mathbf{1}_{N_\nu}\le\psi_1$, and $\widetilde{\theta}\le\theta_1$. We get that
\b*
N&\subset& \{\varphi_1(X)=\infty\}\cup\{\psi_1(Y)=\infty\}\cup \{Y\notin \dom\theta_1(X,\cdot)\}\cup\{Y\notin\aff\,\rf_X\conv\,D(X)\}\\
&=&\{\varphi_1(X)=\infty\}\cup \big\{Y\notin D(X)\cap\aff\,\rf_X\conv\,D(X) \big\}\\
&=&\{\varphi_1(X)=\infty\}\cup \big\{Y\notin D(X)\cap \aff\,I(X) \big\}.
\e*
We have
\be\label{eq:inclusion_N}
N\subset \dom\big(\theta_1+\varphi_1\oplus\psi_1\big)^c,&\mbox{and}N\subset\{X\in N_\mu\}\cup\{Y\in N_\mu\}\cup\{Y\notin \bar J(X)\}&.
\ee
Notice that as $\mu[\varphi_1]+\nu[\psi_1]=0$, $\{\varphi_1=\infty\}\in\Nc_\mu$ and $\{\psi_1=\infty\}\in\Nc_\nu$. We also have $\nu\overline{\ominus}\mu{}[\theta_1] <\infty$. We may replace $\varphi_1$ by $\infty\mathbf{1}_{\varphi_1=\infty}$, $\psi_1$ by $\infty\mathbf{1}_{\psi_1=\infty}$, and $\theta_1$ by $\infty\mathbf{1}_{\theta_1=\infty}\in\widetilde{\Tc}(\mu,\nu)$, where the fact that $\infty\mathbf{1}_{\theta_1=\infty}\in\widetilde{\Tc}(\mu,\nu)$ stems from the fact that $\frac{1}{n}\theta_1\rightsquigarrow \infty\mathbf{1}_{\theta_1=\infty}\in\widetilde{\Tc}(\mu,\nu)$, proving as well that
\be\label{eq:numuzero}
\nu\overline{\ominus}\mu{}[\infty\mathbf{1}_{\theta_1=\infty}] = 0.
\ee
Thanks to these modifications, $\varphi_1$, $\psi_1$, and $\theta_1$ only take the values $0$ or $\infty$.

\no{\underline{Step 2:}} Now let a Borel set $N_1\in\Nc_{\mu,\nu}$ be such that $\theta_1$ is a $N_1-$tangent convex function. Then similar to what was done for $N$, we may find $(\varphi_2,\psi_2,\theta_2)\in\Leb^1_+(\mu)\x\Leb^1_+(\nu)\x\widetilde{\Tc}(\mu,\nu)$ such that
\b*
N_1&\subset& \dom\big(\theta_2+\varphi_2\oplus\psi_2\big)^c.
\e*
Iterating this process for all $k\ge 2$, we define $(N_{k},\varphi_k,\psi_k,\theta_k)$ for all $k\ge 1$. Now let
\b*
\varphi_\infty:= \sum_{k\ge 1}\varphi_k\in\Leb_+^1(\mu),& \psi_\infty:= \sum_{k\ge 1}\psi_k\in\Leb_+^1(\nu),&
\mbox{and }\theta':= \theta+\sum_{k\ge 1} \theta_k\in\widetilde{\Tc}(\mu,\nu)\ge \theta.
\e*
Let $N_\mu^0:=(\dom\varphi_\infty)^c$, and $N_\nu^0:=(\dom\psi_\infty)^c$. Notice that $\mu[\varphi_\infty]=\nu[\psi_\infty]=0$, and therefore, $(N_\mu^0,N_\nu^0)\in\Nc_{\mu}\x\Nc_{\nu}$. We now fix $(N_\mu^0,N_\nu^0)\subset (N_\mu,N_\nu)\in\Nc_{\mu}\x\Nc_{\nu}$, and denote $\varphi := \infty\mathbf{1}_{N_\mu}$, and $\psi := \infty\mathbf{1}_{N_\nu}$.

Recall that $J(X):=\conv\,\dom\big(\theta'(X,\cdot)+\psi\big)\cap \aff I(X) = \conv\,D_\infty(X)\cap \aff D_\infty(X)$, where we denote $D_\infty(X):= \dom\big(\theta'(X,\cdot)+\psi\big)$
By Proposition 2.1 (ii) in \cite{de2017irreducible}, $\conv\,D_\infty(x)\setminus\rf_x\conv\,D_\infty(x)$ is convex for $x\in\R^d$. Therefore, we may find $u(x)\in (\aff\rf_x\conv\,D_\infty(x)-x)^\perp$ such that $y\in \conv\,D_\infty(x)\setminus\rf_x\conv\,D_\infty(x)$ implies that $u(x)\cdot(y-x)>0$ by the Hahn-Banach theorem, so that
$$J(X) = D_\infty(X)\cap\aff\,\rf_X\conv\,D_\infty(X)=\dom\big((\theta'+\infty u^\otimes)(X,\cdot)+\psi\big),$$
with the convention $\infty-\infty = \infty$. Finally,
\be\label{eq:inclusion_cup_N}
N' &=& \{X\in N_\mu\}\cup\{Y\in N_\nu\}\cup\{Y\notin J(X)\}\nonumber\\
&=& \dom(\varphi\oplus\psi + \infty u+\theta')^c\nonumber\\
&\supset& \cup_{k\ge 1} N_k\cup (\dom\theta_k)^c\cup N\\
&\supset& N.\nonumber
\ee
We proved the inclusion from (ii).

\no{\underline{Step 3:}} Now we prove that $N'^c$ is martingale monotone, which is the end of (ii). Let $\P$ with finite support such that $\P[N'^c] = 1$, and $\P'$ a competitor to $\P$. Let $k\ge 1$, we have $\P[N_k^c] = 1$ by \eqref{eq:inclusion_cup_N}, therefore, as $\theta_k$ is a $N_k-$tangent convex function, $\P'[\theta_k]\le\P[\theta_k]$, therefore, as by \eqref{eq:inclusion_cup_N} we have that $\P[dom\theta_k]=1$, we also have that $\P'[\dom\theta_k]=1$. As this holds for all $k\ge 1$, and for $N$ and the $N-$tangent convex function $\theta$, we have $\P'[\dom\theta']=1$. Now as $\P[\dom\varphi\x\dom\psi]=1$, we clearly have $\P'[\dom\varphi\x\dom\psi]=1$. Recall that by construction, $\dom\big(\theta'+\varphi
\oplus\psi\big)=\big(\infty\mathbf{1}_{\theta' = \infty}+\varphi\oplus\psi\big)^{-1}(0)$, therefore, $\P[\infty\mathbf{1}_{\theta' = \infty}+\varphi\oplus\psi]=\P'[\infty\mathbf{1}_{\theta' = \infty}+\varphi\oplus\psi]=0$. Let $n\ge 1$, $\P[\infty\mathbf{1}_{\theta' = \infty}+n u^\otimes+\varphi\oplus\psi]=\P'[\infty\mathbf{1}_{\theta' = \infty}+n u^\otimes+\varphi\oplus\psi]=0$. As $u^\otimes$ is negative only where the rest of the function is infinite, $\infty\mathbf{1}_{\theta' = \infty}+n u^\otimes+\varphi\oplus\psi\ge 0$ for all $n\ge 1$. Then by monotone convergence theorem, $\P'[\infty\mathbf{1}_{\theta' = \infty}+\infty u^\otimes+\varphi\oplus\psi]=\P[\infty\mathbf{1}_{\theta' = \infty}+\infty u^\otimes+\varphi\oplus\psi] =0$. Therefore, $\P'[N'] = 0$, proving that $N'^c$ is martingale monotone.


\no{\underline{Step 4:}} Now we prove that $J(X)= \conv\big(J(X)\setminus N_\nu\big)\subset \dom\theta(X,\cdot)$, which is the first part of (i).
$$
\dom\big((\theta'+\infty u^\otimes)(X,\cdot)+\psi\big)\subset J(X)\cap \dom\psi_\infty\subset J(X).
$$
Passing to the convex hull, we get $J(X)=\conv\big(J(X)\cap \dom\psi\big)=\conv\big(J(X)\setminus N_\nu\big)$ as $N_\nu=\{\psi=\infty\}$.

\no{\underline{Step 5:}} Now we prove that $J(X)\subset \dom\theta'(X,\cdot)\subset\dom\theta(X,\cdot)$, which is the second part of (i). Let $x\in\R^d$, and $y\in J(x)$. Then $y =\sum_i\lambda_iy_i$, convex combination, with $(y_i)\subset \dom\theta'(x,\cdot)\cap\dom\psi$. Let $\P:=\sum_i\lambda_i\delta_{(x,y_i)}+\delta_{(y,y)}$. Let $k\ge 1$, $\P[N_k^c\cup\{X=Y\}]=1$, $\P[\theta_k]<\infty$, and therefore, as $\P':=\sum_i\lambda_i\delta_{(y,y_i)}+\delta_{(x,y)}$ is a competitor to $\P$, $\P'[\theta_k]\le \P[\theta_k]<\infty$, and $y\in \dom\theta_k(x,\cdot)$. $J(x)\subset \dom\theta_k(x,\cdot)$ for all $k\ge 1$, $J(X)\subset \dom\theta'(X,\cdot)$ on $N_\mu^c$.

\no{\underline{Step 6:}} Now we prove that up to choosing well $I$, and up to a modification of $J$ on a $\mu-$null set, $I\subset J\subset \bar{J}\subset \cl I$, and $J$ is constant on $I(x)$, for all $x\in\R^d$, which is the part concerning $J$ of the end of (iii).

We have that $\{I(x),x\in\R^d\}$ is a partition of $\R^d$, $I\subset \bar J\subset \cl I$, and $\bar J$ is constant on $I(x)$ for all $x$. By looking at the proof of Theorem 2.1 in \cite{de2017irreducible}, we may enlarge the $\mu-$null set $N_\mu^I\in\Nc_\mu$ such that $I = \{X\}$ on $\big(\cup_{x'\notin N_\mu^I} I(x')\big)^c$. We do so by requiring that $N_\mu\subset N_\mu^I$. Now we prove that $J$ is constant on $I(X)$, $\mu-$a.s. Let $x_1,x_2\in \dom\varphi_\infty$, and $y\in \dom(\theta_\infty+\infty u_\infty^\otimes)(x_1,\cdot)\cap \dom\psi_\infty$, then $y-\epsilon (y-x_2)\in I(x_2)$ for $\epsilon>0$ small enough, as $x_2\in\ri I(x_2)=\ri I(x_1)$, and $y\in\cl I(x_1)$, as $J(X)\subset \bar J(X)\subset \cl I(X)$ by \eqref{eq:inclusion_N}. Then we may find $x_1=\sum_i\lambda_iy_i+\lambda y$, convex combination, with $(y_i)\subset \dom(\theta_\infty+\infty u_\infty^\otimes)(x_1,\cdot)\cap \dom\psi_\infty$, and $\lambda>0$. Then let $\P:=\sum_i\lambda_i\delta_{(x_1,y_i)}+\lambda \delta_{(x_1,y)}+\delta_{(x_2,x_2)}$. For all $k\ge 1$, notice that $\P[N_k^c\cup\{X=Y\}] = 1$, as $x_2\in (N_k^c)_{x_2}$. Notice furthermore that $\P[\theta_k]<\infty$, and that $\P':=\sum_i\lambda_i\delta_{(x_2,y_i)}+\lambda \delta_{(x_2,y)}+\delta_{(x_1,x_2)}$ is a competitor to $\P$. Then as $\theta_k$ is a $N_k-$tangent convex function, $\P'[\theta_k]\le\P[\theta_k]<\infty$, and therefore, as $\lambda >0$, $\theta_k(x_2,y)<\infty$. We proved that
$$\dom(\theta_\infty+\infty u_\infty^\otimes)(x_1,\cdot)\cap \dom\psi_\infty\subset \dom\theta_k(x_2,\cdot).$$

Therefore, $\dom(\theta_\infty+\infty u_\infty^\otimes)(x_1,\cdot)\cap \dom\psi_\infty\subset \dom\theta_\infty(x_2,\cdot)$. As the other ingredients of $J$ do not depend on $x$, and as we can exchange $x_1$, and $x_2$ in the previous reasoning,
$$\dom(\theta_\infty+\infty u_\infty^\otimes)(x_1,\cdot)\cap \dom\psi_\infty= \dom(\theta_\infty+\infty u_\infty^\otimes)(x_2,\cdot)\cap \dom\psi_\infty.$$
Taking the convex hull, we get $J(x_1)=J(x_2)$.

\no{\underline{Step 7:}} Now we prove that thanks to the modification of $I$ and $J$, we have that $\Jo$ is constant on all $I(x)$, for $x\in\R^d$, and that $\underline{J}\subset\Jo\subset J$, which is the remaining part of (iii). By its definition, we see that the dependence of $\Jo$ in $x$ stems from a direct dependence in $J(x)$. The map $J$ is constant on each $I(x)$, $x\in\R^d$, whence the same property for $\Jo$. Now for $I(x)\notin I(N_\mu^c)$, all these maps are equal to $\{x\}$, whence the inclusions and the constance.

Now we claim that for $x,x'\in\R^d$ such that $x'\in J(x)$, we have $J(x')\subset J(x)$. This claim will be justified in (iii) of the proof of Remark \ref{rmk:J} above. Now if $x'\in J(x)\setminus N_\mu\subset J(x)$, we have as a consequence that $J(x')\subset J(x)$, and therefore $I(x')\subset J(x)$. We proved that $\Jo\subset J$.

Finally by Proposition 2.4 in \cite{de2017irreducible}, we may find $\widehat{\P}\in\Mc(\mu,\nu)$ such that $\underline{J}(X)\setminus I(X)\subset \{y:\widehat{\P}_X[\{y\}]\}$, on $N_\mu$. Then $\underline{J}\subset \Jo$ on $N_\mu^c$. Otherwise, these maps are again equal to $\{X\}$, whence the result.
\ep\\

\no{\bf Proof of Remark \ref{rmk:J}}
\no{\underline{(i)}} Recall that, with the notations from Proposition \ref{prop:J}, $J(X):=\conv(\dom\theta'(X,\cdot)\setminus N_\nu)\cap\bar J(X)$. $\theta'\in\widetilde{\Tc}(\mu,\nu)$, then $\Sbf_{\mu,\nu}(\theta)<\infty$ and $Y\in \dom\theta'(X,\cdot)$, $\Mc(\mu,\nu)-$q.s. Recall that $Y\in \bar J(X)$, and $Y\notin N_\nu$, $\Mc(\mu,\nu)-$q.s. All these ingredients prove that $Y\in J(X)$, q.s. and $Y\in J(X)\setminus N_\nu$, q.s. The result for $\Jo$ is a consequence of the inclusion \be\label{eq:inclusion_JN_Jo_J}
J\setminus N_\nu\subset\Jo\subset J.
\ee

\no{\underline{(ii)}} Let $x,x'\in\R^d$, we prove that $J(x)\cap J(x') = \aff\big(J(x)\cap J(x')\big)\cap J(x)$. The direct inclusion is trivial, let us prove the indirect inclusion. We first assume that $x,x'\in N_\mu^c$. We claim that
\be\label{eq:capsetminusMnu}
J(x)\cap J(x') = \conv\big(J(x)\cap J(x')\setminus N_{\nu}'\big).
\ee
This claim will be proved in (iii). If $J(x)\cap J(x')= \emptyset$, the assertion is trivial, we assume now that this intersection is non-empty. Let $y_1,...,y_k\in J(x)\cap J(x')\setminus N_{\nu}'$ with $k\ge 1$, spanning $\aff\big(J(x)\cap J(x')\setminus N_{\nu}'\big)$. Let $y\in \aff\big(J(x)\cap J(x')\big)\cap J(x)$, and $y' = \frac{1}{k}\sum_iy_k$. We have $y'\in\ri\,\conv(y_1,...,y_k)$ and $y\in\aff\,\conv(y_1,...,y_k)$, therefore, for $\eps>0$ small enough, $\eps y+(1-\eps)y'\in \ri\,\conv(y_1,...,y_k)\subset J(x)\cap J(x')\subset J(x') = \conv\big(J(x')\setminus N_\nu\big)$ by (i). Then, for $\eps$ small enough, $\eps y+(1-\eps)y' = \sum_i \lambda_i y_i'$, convex combination, with $(y_i)_i\subset J(x')\setminus N_\nu$. Then $\P = \frac12\eps \delta_{(x,y)} +\frac{1}{2k}(1-\eps)\sum_i \delta{(x,y_i)}+\frac12\sum_i\lambda_i\delta{(x',y_i')}$ is concentrated on $N'^c$, and by (iv) we have that its competitor $\P' = \frac12\eps \delta_{(x',y)} +\frac{1}{2k}(1-\eps)\sum_i \delta{(x',y_i)}+\frac12\sum_i\lambda_i\delta{(x,y_i')}$ is also concentrated on $N'^c$. Therefore $y \in J(x')$, and as $y\in J(x)$, we proved the reverse inclusion: $J(x)\cap J(x') = \aff\big(J(x)\cap J(x')\big)\cap J(x)$.

Now if $x,x'\in \cup_{x''\notin N_\mu}I(x'')$, we may find $x_1,x_2\in N_\mu^c$ such that $J(x) = J(x_1)$, and $J(x') = J(x_2)$, whence the result from what precedes. Finally if $x$ or $x'$ is not in $\cup_{x''\notin N_\mu}I(x'')$, If it is $x$, then $I(x) = J(x) = \{x\}$, and if $x\in J(x')$, then the result is $\{x\} = \{x\}$, else it is $\emptyset = \emptyset$. If it is $x'$, then if $x'\in J(x)$, the result is $\{x'\} = \{x'\}$, otherwise, it is again $\emptyset = \emptyset$. In all the cases, the result holds.

Finally we extend this result to $\Jo$. Notice that by \eqref{eq:capsetminusMnu} together with \eqref{eq:inclusion_JN_Jo_J}, we have $\aff\big(J(x)\cap J(x')\big) = \aff\big(\Jo(x)\cap \Jo(x')\big)$. Now consider the equation that we previously proved $J(x)\cap J(x') = \aff\big(J(x)\cap J(x')\big)\cap J(x)$, subtracting $N_\nu\setminus \cup_{x''\notin N_\mu}I(x'')$ and replacing $\aff\big(J(x)\cap J(x')\big)$, we get $\Jo(x)\cap \Jo(x') = \aff\big(\Jo(x)\cap \Jo(x')\big)\cap \Jo(x)$.

\no{\underline{(iii)}} Let $y\in J(x)\cap J(x')$. By (i), $\conv\big(J(x)\setminus N_\nu\big) = J(x)$, and the same holds for $x'$. Then we may find $y_1,...,y_k\in J(x)\setminus N_\nu$ and $y_1',...,y_{k'}'\in J(x')\setminus N_\nu$ with $\sum_i\lambda_i y_i = \sum_i\lambda_i' y_i' = y$, where the $(\lambda_i)$ and $(\lambda_i')$ are non-zero coefficients such that the sums are convex combinations. Now notice that $\P := \frac12\sum_i \lambda_i\delta_{(x,y_i)}+\frac12\sum_i \lambda_i' \delta_{(x',y_i')}$ is supported in $N'^c$. By (iv), its competitor $\P' := \frac12\sum_i \lambda_i\delta_{(x',y_i)}+\frac12\sum_i \lambda_i' \delta_{(x,y_i')}$ is also supported on $N'^c$. Therefore, $y_1,...,y_k,y_1',...,y_{k'}'\in J(x)\cap J(x')\setminus N_\nu$. We proved that $J(x)\cap J(x')\subset\conv\big(J(x)\cap J(x')\setminus N_{\nu}'\big)$, and therefore as the other inclusion is easy, we have $J(x)\cap J(x') = \conv\big(J(x)\cap J(x')\setminus N_{\nu}'\big)$. The extension of this result for $\Jo$ is again a consequence of the inclusion \eqref{eq:inclusion_JN_Jo_J}.

\no{\underline{(iv)}} Now we assume additionally that $I(x')\cap J(x)\neq \emptyset$, let us prove that then $J(x')\subset J(x)$. If $x'\notin \cup_{x''\notin N_\mu}I(x'')$, then $J(x') = \{x'\}$ and the result is trivial. If $x\notin \cup_{x''\notin N_\mu}I(x'')$, then the result is similarly trivial. By constance of $J$ and $I$ on $I(x)$ for all $x$, we may assume now that $x,x'\in N_\mu^c$. Then let $y\in I(x')\cap J(x)\subset \conv\big(J(x')\setminus N_\nu\big)\cap \conv\big(J(x)\setminus N_\nu\big)$. Let $y'\in J(x')\setminus N_\nu$, for $\eps>0$ small enough, $y-\eps(y'-y)\in I(y')$ by the fact that $I(y')$ is open in $\aff J(x')$. Then $y-\eps(y'-y) = \sum_i \lambda_i y_i$, and $y = \sum_i \lambda_i' y_i'$, convex combinations where $(y_i)_i\subset J(x')\setminus N_\nu$, and $(y_i)_i\subset J(x')\setminus N_\nu$. Then $\P := \frac12\frac{\eps}{1+\eps}\delta_{(x,y')}+\frac12\frac{1}{1+\eps}\sum_i\lambda_i\delta_{(x,y_i)}+\frac12\sum_i\lambda_i'\delta_{(x',y_i')}$ is concentrated on $N'^c$, and by (iv), so does its competitor $\P' := \frac12\frac{\eps}{1+\eps}\delta_{(x',y')}+\frac12\frac{1}{1+\eps}\sum_i\lambda_i\delta_{(x',y_i)}+\frac12\sum_i\lambda_i'\delta_{(x,y_i')}$. Then in particular, $y'\in J(x)$. Finally, $J(x')\setminus N_\nu\subset J(x)$, passing to the convex hull, we get that $J(x')\subset J(x)$.

Finally, if $I(x')\cap \Jo(x)\neq\emptyset$, then $I(x')\cap J(x)\neq\emptyset$, and $J(x')\subset J(x)$. Subtracting $N_\nu\setminus \cup_{x''\notin N_\mu}I(x'')$ on both sides, we get $\Jo(x')\subset \Jo(x)$.
\ep\\

\no{\bf Proof of Theorem \ref{thm:polar}}
Let $(N_\mu,N_\nu)\in\Nc_\mu\x\Nc_\nu$, and $J\in{\Jc}(\mu,\nu)$. The "if" part holds as $Y\in J(X)$, $X\notin N_\mu$, and $Y\notin N_\nu$ q.s.

Now, consider an analytic set $N\in\Nc_{\mu,\nu}$. Then $c:=\infty\mathbf{1}_N$ is upper semi-analytic non-negative. Notice that $\Sbf_{\mu,\nu}(c)=0$. By Theorem \ref{thm:dualitymod}, we may find $(\varphi,\psi, h,\theta)\in\Dc^{\widetilde{mod}}_{\mu,\nu}(c)$ such that $\mu[\varphi]+\nu[\psi]+\nu\widehat{\ominus}\mu[\theta]=\Sbf_{\mu,\nu}(c) = 0$. Then by the pointwise inequality $\infty\mathbf{1}_N\le \varphi\otimes\psi+h^\otimes+\theta$, on $\{Y\in\aff\,\rf_X\conv\,D(X)\}$, with $D(X) = \dom\big(\theta(X,\cdot)+\psi\big)$, we get that
\b*
N&\subset& \{\varphi(X)=\infty\}\cup\{\psi(Y)=\infty\}\cup \{Y\notin \dom\theta(X,\cdot)\}\cup\{Y\notin\aff\,\rf_X\conv\,D(X)\}\\
&=&\{\varphi(X)=\infty\}\cup\{\psi(Y)=\infty\}\cup \big\{Y\notin \dom\theta(X,\cdot)\cap\aff\,\rf_X\conv\,D(X) \big\},
\e*
Let $J\in\Jc(\mu,\nu)$ from Proposition \ref{prop:J} for $\theta$, and $N_\nu:=\dom\psi^c\in \Nc_\nu$. We have $J(X)\subset \aff\,\rf_X\conv\,D(X)$ and $J(X)\subset \bar J(X)\subset \dom\theta(X,\cdot)$, $\mu-$a.s. Therefore, we have
\b*
N&\subset& N_0:=\{X\in N_\mu\}\cup\{Y\in N_\nu\}\cup \{Y\notin J(X)\},
\e*
for some $N_\mu\in\Nc_\mu$, and $N_\nu\subset N_\nu\in\Nc_\nu$. By Proposition \ref{prop:J} (i) and (iv), $N_0$ may be chosen canonical up to enlarging $N_\mu$.
\ep

\subsection{Decomposition in irreducible martingale optimal transports}\label{subsect:proofdisintegretion}

In order to prove theorem \ref{thm:desintegration}, we first need to establish the following lemma.

\begin{Lemma}\label{lemma:existdecomp}
Let $\theta\in\widehat{\Tc}(\mu,\nu)$ and $\P\in\Mc(\mu,\nu)$, we may find $\theta'\in \widehat{\Tc}(\mu,\nu)$ such that $\theta\le\theta'$, $\nu\widehat{\ominus}\mu[\theta']\le\nu\widehat{\ominus}\mu[\theta]$, and $\int_{I(\R^d)}\nu^{\P}_I\widehat{\ominus}\mu_I[\theta']\eta(dI)\le\nu\widehat{\ominus}\mu[\theta]$. Furthermore under Assumption \ref{ass:domination}, we may find $f\in\Cfrak_{\mu,\nu}$ and $p\in\partial^{\mu,\nu}f$ such that $\theta\le\Tbf_pf$, q.s., $\nu\overline{\ominus}\mu[f]\le\nu\widehat{\ominus}\mu[\theta]$, and $\int_{I(\R^d)}\nu^{\P}_I\overline{\ominus}\mu_I[f]\eta(dI)\le\nu\widehat{\ominus}\mu[\theta]$.
\end{Lemma}
\proof
Let $a> 0$, we consider $\Tc$ the collection of $\widehat\theta\in\Theta_\mu$ such that we may find $\theta'\in\widehat{\Tc}_a$ with $\theta'\ge\widehat\theta$, $\nu\widehat{\ominus}\mu[\theta']\le a$, and $\int_{I(\R^d)}\nu^{\P}_I\widehat{\ominus}\mu_I[\theta']\eta(dI)\le a$. First we have easily $\Tbf(\Cfrak_a)\subset\Tc$, as $\Tbf(\Cfrak_a)\subset\widetilde{\Tc}_a$, and $\int_{I(\R^d)}\nu^{\P}_I\widehat{\ominus}\mu_I[\theta']\eta(dI)=\int_{I(\R^d)}(\nu^{\P}_I-\mu_I)[\theta']\eta(dI)=(\nu-\mu)[\theta']$, for $\theta'\in\Tbf(\Cfrak_a)$. Now we consider $(\widehat\theta_n)_{n\ge 1}\subset \Tc$ converging $\muxpw$ to $\underline{\widehat\theta}_\infty$. For each $n\ge 1$, we may find $\theta_n\in\widehat{\Tc}_a$ such that $\theta_n\ge\widehat\theta_n$, $\nu\widehat{\ominus}\mu[\theta_n]\le a$, and $\int_{I(\R^d)}\nu^{\P}_I\widehat{\ominus}\mu_I[\theta_n]\eta(dI)\le a$. By the Koml\'os Lemma on $I\longmapsto \nu^{\P}_I\widehat{\ominus}\mu_I[\theta_n]$ under the probability $\eta$ together with Lemma 2.12 in \cite{de2017irreducible}, we may find convex combination coefficients $(\lambda_k^n)_{1\le n\le k}$ such that $\sum_{k=n}^\infty \lambda_k^n \nu^{\P}_I\widehat{\ominus}\mu_I[\theta_k]$ converges $\eta-$a.s. and $\theta_n':=\sum_{k=n}^\infty \lambda_k^n\theta_k$ converges $\muxpw$ to $\theta':= \underline\theta_\infty'$, as $n\longrightarrow\infty$, and moreover $\nu\widehat{\ominus}\mu[\theta']\le a$. As $\theta_n'$ is a convex extraction of $\widehat\theta_n$, we have $\theta':= \underline\theta_\infty'\ge \underline{\widehat\theta}_\infty$. Moreover, by convexity of $\nu^{\P}_I\widehat{\ominus}\mu_I$, we have $\sum_{k=n}^\infty \lambda_k^n \nu^{\P}_I\widehat{\ominus}\mu_I[\theta_k]\ge \nu^{\P}_I\widehat{\ominus}\mu_I[\theta_n']$, and therefore
$$\liminf_{n\to\infty}\sum_{k=n}^\infty \lambda_k^n \nu^{\P}_I\widehat{\ominus}\mu_I[\theta_k]=\limsup_{n\to\infty}\sum_{k=n}^\infty \lambda_k^n \nu^{\P}_I\widehat{\ominus}\mu_I[\theta_k]\ge \limsup_{n\to\infty}\nu^{\P}_I\widehat{\ominus}\mu_I[\theta_n']\ge \nu^{\P}_I\widehat{\ominus}\mu_I[\theta']$$
$\eta-$a.s. Integrating this inequality with respect to $\eta$, and using Fatou's Lemma, we get
$$\int_{I(\R^d)}\nu^{\P}_I\widehat{\ominus}\mu_I[\theta']\eta(dI)\le a.$$
Then $\underline{\widehat\theta}_\infty\in\Tc$. Hence, $\Tc$ is $\muxpw-$Fatou closed, and therefore $\widehat\Tc_a\subset\Tc$.

Now let $\widehat\theta\in\widehat{\Tc}(\mu,\nu)$, with $l:=\nu\widehat{\ominus}\mu[\widehat\theta]$. By the previous step, for all $n\ge 1$, we may find $\theta_n'\in\widehat{\Tc}_{l+1/n}$ with $\int_{I(\R^d)}\nu^{\P}_I\widehat{\ominus}\mu_I[\theta_n']\eta(dI)\le l+1/n$ such that $\widehat\theta\le\theta_n'$. Similar to the proof of Lemma \ref{lemma:existencehat}, we get $\theta'\in \widehat{\Tc}(\mu,\nu)$ such that $\theta\le\theta'$, $\nu\widehat{\ominus}\mu[\theta']\le l$, and $\int_{I(\R^d)}\nu^{\P}_I\widehat{\ominus}\mu_I[\theta']\eta(dI)\le l$, thus proving the result.

We prove the second part of the Lemma similarly, using Assumption \ref{ass:domination} instead of Lemma 2.12 in \cite{de2017irreducible}.
\ep\\

For the proof of next result, we need the following lemma:

\begin{Lemma}\label{lemma:recomposition}
Let $\theta\in\widehat{\Tc}(\mu,\nu)$, $m_X:= \mu[X|I(X)]$, and $f_X(\cdot):= \theta(m_X,\cdot)$. Then we may find a $\mu-$unique measurable $\widehat p(X)\in\aff I(X)-X$ such that for some $N_\mu\in\Nc_\mu$,
\be\label{eq:equalityphat}
\theta=f_X(Y)-f_X(X)-\widehat p(X)\cdot (Y-X),&\mbox{on }\{Y\in\aff\dom_X\theta\}\cap\{X\notin N_\mu\}.&
\ee
\end{Lemma}
\proof
We consider $N_\mu\in\Nc_\mu$ from Proposition 2.10 in \cite{de2017irreducible}, so that for $x_1,x_2\notin N_\mu$, $y_1,y_2\in\R^d$, and $\lambda\in[0,1]$ with $\yb:= \lambda y_1 + (1-\lambda)y_2\in\dom_{x_1}\theta\cap\dom_{x_2}\theta$, we have:
\be\label{eq:x-consistency}
 \lambda \theta(x_1,y_1)+(1-\lambda)\theta(x_1,y_2)-\theta(x_1,\yb)
 =
\lambda \theta(x_2,y_1)+(1-\lambda)\theta(x_2,y_2)-\theta(x_2,\yb)
 \ge 0.&
\ee
By possibly enlarging $N_\mu$, we may suppose in addition that $I(x)\subset\dom_x\theta$ for all $x\in N_\mu^c$. For $x\in N_\mu^c$ and $y\in\dom_x\theta$, we define $H_x(y):= f_x(y)-f_x(x)-\theta(x,y)$. By \eqref{eq:x-consistency}, $H_x$ is affine on $\aff\,\dom_x\theta\cap \dom\theta(x,\cdot)$. Indeed let $y_1\in \aff\,\dom_x\theta\cap \dom\theta(x,\cdot)$, $y_2\in \dom_x\theta$, and $0\le\lambda\le1$, then $\yb:=\lambda y_1+(1-\lambda)y_2\in \dom_x\theta$ and
\b*
H_x(\yb)&=&\theta(m_x,\yb)-\theta(m_x,x)-\theta(x,\yb)\\
&=&\lambda \theta(m_x,y_1)+(1-\lambda)\theta(m_x,y_2)-\lambda \theta(x,y_1)-(1-\lambda)\theta(x,y_2)-\theta(m_x,x)\\
&=&\lambda H_x(y_1)+(1-\lambda) H_x(y_2)
\e*

We notice as well that $H_x(x) = 0$. Then we may find a unique $\widehat p(x)\in\aff I(x)-x$ so that for $y\in\dom_x\theta$, $H_x(y) = \widehat p(x)\cdot(y-x)$. $\widehat p(X)$ is measurable and unique on $N_\mu^c$, and therefore $\mu-$a.e. unique. For $y\in\aff\,\dom_x\theta\cap \dom\theta(x,\cdot)$, it gives the desired equality \eqref{eq:equalityphat}. Now for $y\in \aff\,\dom_x\theta\cap \dom\theta(x,\cdot)^c$, let $0<\lambda< 1$ such that $\yb:=\lambda x+(1-\lambda)y\in\dom_x\theta$. By \eqref{eq:x-consistency},  $\lambda \theta(x,x)+(1-\lambda)\theta(x,y)-\theta(x,\yb)
 =
\lambda f_x(x)+(1-\lambda)f_x(y)-f_x(\yb)$, and therefore $\theta(x,y)$ is finite if and only if $f_x(y)$ is finite. This proves that \eqref{eq:equalityphat} holds for $y\in\aff\,\dom_x\theta$.
\ep\\

\no {\bf Proof of Theorem \ref{thm:desintegration}}
For $\P\in\Mc(\mu,\nu)$, $I_0\in I(\R^d)$, we have by definition of the supremum,

$$\P_{I_0}[c]\leq \Sbf_{\P_{I_0}\circ X^{-1},\P_{I_0}\circ Y^{-1}}(c)= {\bf S}_{\mu_{I_0},\nu^\P_{I_0}}(c),$$

where we denote by $\P_{I}$ a conditional disintegration of $\P$ with respect to the random variable $I$. Now we consider a minimizer for the dual problem $(\bar\varphi,\bar\psi,\bar h,\widehat\theta)\in \Dc^{mod}_{\mu,\nu}(c)$ and $\theta'\in \widehat{\Tc}(\mu,\nu)$ such that $\theta\le\theta'$, $\nu\widehat{\ominus}\mu[\theta']\le\nu\widehat{\ominus}\mu[\theta]$, and $\int_{I(\R^d)}\nu^{\P}_I\widehat{\ominus}\mu_I[\theta']\eta(dI)\le\nu\widehat{\ominus}\mu[\theta]$ from Lemma \ref{lemma:existdecomp}. Recall the notation $m_X:= \mu[X|I(X)] = \P[Y|I(X)]$, by the martingale property, and let $f_X(Y):= \theta'(m_X,Y)$. From Lemma \ref{lemma:recomposition}, we have $\theta'(X,Y) = f_X(Y)-f_X(X)-p_X(X)\cdot(Y-X)$, with $p_X\in\partial f_X(X)$, $\Mc(\mu,\nu)-$q.s. Then let $\varphi:=\bar\varphi - f_X$, $\psi_I(X):=\bar\psi(Y) + f_X(Y)$, $h:=\bar h-p_X$.
$$\mu_I[\bar\varphi]+\nu^\P_I[\bar\psi]+\nu^\P_I\widehat{\ominus}\mu_I[\theta']\ge \mu_I[\varphi]{\oplus}\nu_I^\P[\psi]\geq {\bf I}_{\mu_I,\nu_I^\P}(c) \geq{\bf S}_{\mu_I,\nu^\P_{ I}}(c).$$
Integrating with respect to $\eta$, we get:
\b*
{\bf I}^{mod}_{\mu,\nu}(c) &\ge&\int_{I(\R^d)} \mu_I[\bar\varphi]+\nu^\P_I[\bar\psi]+\nu^\P_I\widehat{\ominus}\mu_I[\theta']\eta(dI) \geq\int_{I(\R^d)} {\bf I}_{\mu_I,\nu^\P_{ I}}(c)\eta(dI)\\
&\ge& \int_{I(\R^d)} {\bf S}_{\mu_I,\nu^\P_{ I}}(c)\eta(dI)\geq \P[c].
\e*
Taking the supremum over $\P$:
$${\bf I}^{mod}_{\mu,\nu}(c)\geq\underset{\P\in\Mc(\mu,\nu)}{\sup}\int_{I(\R^d)} {\bf I}_{\mu_I,\nu^\P_{ I}}(c)\eta(dI)\geq \underset{\P\in\Mc(\mu,\nu)}{\sup}\int_{I(\R^d)} {\bf S}_{\mu_I,\nu^\P_{ I}}(c)\eta(dI)\geq {\bf S}_{\mu,\nu}(c)$$
Then all the inequalities are equalities by the duality Theorem 3.8 in \cite{de2017irreducible}.

We consider $\P^*$ such that $\P^*[c]=\Sbf_{\mu,\nu}(c) = \Ibf_{\mu,\nu}^{mod}(c)$ gives us that there is an optimizer.
\b*
\Sbf_{\mu,\nu}(c)
&=&
\P^*[c]
 =
 \int_{I(\R^d)} \P^*_I[c]\eta(dI)
 \leq
 \int_{I(\R^d)}{\bf S}_{\mu_I,\nu^{\P^*}_{I}}(c)\eta(dI)
 \leq
 \int_{I(\R^d)} {\bf I}^{mod}_{\mu_I,\nu^{\P^*}_{I}}(c)\eta(dI)\\
 &\leq& 
 \int_{I(\R^d)} \mu_I[\bar\varphi]+\nu^{\P^*}_I[\bar\psi]+\nu^{\P^*}_I\widehat{\ominus}\mu_I[\theta']\eta(dI)
 \leq
 \Ibf_{\mu,\nu}^{mod}(c).
 \e*
Then all these inequalities are equalities by duality.

The second part is proved similarly, using the second part of Lemma \ref{lemma:existdecomp}.
\ep

\subsection{Properties of the weakly convex functions}\label{subsect:CctoTc}



The proof of Proposition \ref{prop:CctoTc} is very technical and requires several lemmas as a preparation.

\begin{Lemma}\label{lemma:measselect}
Let $N_\mu\in\Nc_\mu$, we may find $N_\mu\subset N_\mu'\in\Nc_\mu$, and a Borel mapping $\ri\Kcirc\ni K\longmapsto m_K$ such that $m_{I(X)}\in I(X)\setminus N_\mu'$ on $\{X\notin N_\mu'\}$.
\end{Lemma}
\proof
We may approximate $N_\mu^c$ from inside by a countable non-decreasing sequence of compacts $(K_n)_{n\ge 1}$: $\cup_{n\ge 1}K_n\subset N_\mu^c$, and $\mu[\cup_{n\ge 1}K_n] = 1$. Let $N_\mu':=(\cup_{n\ge 1}K_n)^c\in\Nc_\mu$. For $n\in\N$, the mapping $I_n:x\mapsto x+(1-1/n)\big(\cl I(x)-x\big)\cap K_n$ is measurable with closed values. Then we deduce from Theorem 4.1 of the survey on measurable selection \cite{wagner1977survey} that we may find a measurable selection $m^n:\R^d\longrightarrow\R^d$ such that $m^n(x)\in I_n(x)$ for all $x\in\R^d$. Define
 \b*
 m'(x):=m^{n(x)}(x)
 &\mbox{where}&
 n(x):=\inf\{n\ge 1:I_n(x)\neq \emptyset\},
 ~~x\in\R^d,
 \e*
and $m^\infty:=0$. Then for all $x\notin N_\mu'$, we have the inclusion $\emptyset\neq\{x\}\cap N_\mu'^c\subset I(x)\cap N_\mu'^c =\cup_{n\ge 1}I_n(x)$, so that $n(x)<\infty$ and $m'(x)\in I(x)\setminus N_\mu'$. However, we want to find a map from $\Kcirc$ to $\R^d$. Consider again the map $\bar m_I:= \E^\mu[X|I]$. Notice that $\bar m_I\in I$ by the convexity of $I$, and that it is constant on $I(x)$, for all $x\in\R^d$. Then the map $m_I:= m'(\bar m_I)$ satisfies the requirements of the lemma.
\ep\\

We fix a $N-$tangent convex function $\theta\in\widetilde\Tc(\mu,\nu)$. Let $N^0:=\{X\in N_\mu^0\}\cup\{Y\in N_\nu\}\cup\{Y\notin J(X)\}\in\Nc_{\mu,\nu}$, a canonical polar set such that $(N^0)^c\subset N^c\cap \dom\theta$ from Proposition \ref{prop:J}. Consider the map $m_I$ given by Lemma \ref{lemma:measselect} for $N_\mu^0$, let $\Nc_\mu\ni N_\mu\supset N_\mu^0$ such that $m_{I(X)}\in I(X)\setminus N_\mu$ on $\{X\notin N_\mu\}$. By Proposition \ref{prop:J} together with the fact that $N_\mu\supset N_\mu^0$, we may chose the map $I$ so that $N':=\{X\in N_\mu\}\cup\{Y\in N_\nu\}\cup\{Y\notin J(X)\}\in\Nc_{\mu,\nu}$, a canonical polar set such that $N'^c\subset N^c\cap \dom\theta$. For $K\in I(\R^d):=\{I(x):x\in\R^d\}$ we denote $f_K := \theta(m_K,\cdot)$.

\begin{Lemma}\label{lemma:Jo_in_N}
We may find $\Jo\in\Jco(\mu,\nu)$ such that $\{Y\in \Jo(X),X\notin N_\mu\}\subset N^c\cap \dom\theta$, $\conv\Jo = J = \conv\big(J\setminus N_\nu\big)$, and $\conv\big(\Jo(x)\cap\Jo(x')\big) = J(x)\cap J(x') = \conv\big(J(x)\cap J(x')\setminus N_\nu\big)$ for all $x,x'\in\R^d$.
\end{Lemma}
\proof
The map defined by $\Jo(x):= \cup_{x'\in J(x)\setminus N_\mu}I(x')\cup J(x)\setminus N_\nu$ is in $\Jco(\mu,\nu)$. By Proposition \ref{prop:J}, $\Jo\subset J$, therefore $\Jo\subset J = \conv(J\setminus N_\nu)\subset \conv\,\dom\theta(X,\cdot)= \dom\theta(X,\cdot)$ on $N_\mu^c$, whence the inclusion $\{Y\in \Jo(X),X\notin N_\mu\}\subset \dom\theta$.

Now we prove that $\{Y\in \Jo(X),X\notin N_\mu\}\subset N^c$. Recall that $N' = \{Y\in J(X)\setminus N_\nu,X\notin N_\mu\}\subset N^c$. Let $x\notin N_\mu$, and $x'\in J(x)\setminus N_\mu$, then $I(x')\subset N^c_{x'}$. Let $y\in I(x')\subset N'^c_{x'}$, by Proposition \ref{prop:J}, $y\in J(x)\cap J(x') = \conv\big(J(x)\cap J(x')\setminus N_\nu\big)$. Then we may find $y_1,...,y_k\in J(x)\cap J(x')\setminus N_\nu$ such that $y = \sum_i\lambda_iy_i$, convex combination. We also have $y\in N^c_{x'}$, then $\P:=\frac12\sum_i\lambda_i\delta_{(x,y_i)}+\frac12\delta_{x',y}$, and $\P':=\frac12\sum_i\lambda_i\delta_{(x',y_i)}+\frac12\delta_{x,y}$ are competitors such that the only point in their support that may not be in $N^c$ is $(x,y)$, then by Definition \ref{def:Thetamunu} (iii), $(x,y)\in N^c$. We proved that $\{Y\in \Jo(X),X\notin N_\mu\}\subset N^c$.

The other properties are direct consequences of Remark \ref{rmk:J}.
\ep

Let $\Jo\in\Jco(\mu,\nu)$ and $N_\mu\in\Nc_\mu$ from Lemma \ref{lemma:Jo_in_N}.

\begin{Lemma}\label{lemma:theta_is_TpfI}
We have $\theta = \Tbf_{\widehat{p}}f_I(X)$ on $\{Y\in \Jo(X),X\notin N_\mu\}$ for some $\widehat{p}\in\Leb^0(\R^d,\R^d)$, and $\Jo\in\Jco(\mu,\nu)$.
\end{Lemma}
\proof
Let $a_x:=f_{I(x)}-f_{I(x)}(x)-\theta(x,\cdot)$. We claim that $a_x$ is affine on $\Jo(x)$, for all $x\notin N_\mu$, i.e. we may find a measurable map $\widehat{p}$ on $N_\mu^c$ such that, by the above definition of $a_x$ together with the fact that $a_x(x) = 0$,
\b*
\theta = f_{I(X)}(Y)-f_{I(X)}(X)-\widehat p(X)\cdot(Y-X),&\mbox{on }\{Y\in \Jo(X),X\notin N_\mu\}.
\e*
 
  Now we prove the claim. Let $x\notin N_\mu$, and $y,y_1,...,y_k\in \Jo(x)$, for some $k\in\N$, such that $y=\sum_i\lambda_i y_i$, convex combination. Now consider
 \b*
 \P:=\sum_i\delta_{(m_{I(x)},y_i)}+\delta_{x,y}, &\mbox{and}& \P':=\sum_i\delta_{(x,y_i)}+\delta_{m_{I(x)},y}.
 \e*
Notice that $\P$, and $\P'$ are competitors with finite supports, concentrated on $N^c$, by the fact that $m_{I(x)}\notin N_\mu$, together with Lemma \ref{lemma:Jo_in_N}, and the fact that $\Jo$ is constant on $I(x)$ by Proposition \ref{prop:J}. Therefore
 \be\label{eq:convexityinvariance}
 \sum_i \lambda_i\theta(m_{I(x)},y_i)+\theta(x,y)&=&\sum_i\lambda_i\theta(x,y_i)+\theta(m_{I(x)},y),
 \ee
from Definition \ref{def:Thetamunu} (ii). Then the proof that $a_x$ is affine is similar to the proof of Lemma \ref{lemma:recomposition}.

Let $\widehat p(x)$ be a vector in $\nabla \aff I(x)$ representing this linear form. By the fact that $a_x$ is linear and finite on $\Jo(x)$, we have the identity
\be\label{eq:reconstitution}
\theta(x,y) = f_{I(x)}(y)-f_{I(x)}(x)-\widehat p(x)\cdot(y-x),&\mbox{for all }(x,y)\in \{Y\in \Jo(X),X\notin N_\mu\}.
\ee
\ep\\

Recall that we want to find $f:\R^d\longrightarrow\R$, and $p:\R^d\longrightarrow\R^d$ such that $\theta = \Tbf_pf$ on $\{Y\in \Jo(X),X\notin N_\mu\}$. A good candidate for $f$ would be $f_I$, in view of \eqref{eq:reconstitution}. However $f$ defined this way could mismatch at the interface between two components. We now focus on the interface between components. Let $K,K'\in I(\R^d)$, we denote $\interf(K,K'):= \Jo(m_K)\cap \Jo(m_{K'})$ if $m_K,m_{K'}\notin N_\mu$, and $\emptyset$ otherwise.

\begin{Lemma}\label{lemma:existenceTpf}
Let $(A_K)_{K\in I(\R^d)}\subset \Aff(\R^d,\R)$ be such that
\be\label{eq:A_consistency}
A_K(y) - A_{K'}(y) = f_{K'}(y) - f_K(y),&\mbox{for all }y\in \interf(K,K'),&\mbox{and }K,K'\in I(\R^d).
\ee
Then $f(y) := f_K(y) + A_K(y)$ does not depend of the choice of $K$ such that $y\in \Jo(m_K)$, and if we set $p(y):= \widehat{p}(y)+\nabla A_{I(y)}$, we have
\b*
\theta = \Tbf_pf,&\mbox{on}&\{Y\in \Jo(X),X\notin N_\mu\}.
\e*
\end{Lemma}
\proof
Let $K,K'\in I(\R^d)$ such that $y\in \Jo(m_K)\cap \Jo(m_{K'}) = \interf(K,K')$. Then $f_K(y) + A_K(y) - \big(f_{K'}(y) + A_{K'}(y)\big) = 0$ by \eqref{eq:A_consistency}. The first point is proved.

Then $\Tbf_pf = \Tbf_{\widehat{p}+\nabla A_{I}}(f_I+A_I) = \Tbf_{\widehat{p}}f_I+\Tbf_{\nabla A_{I}}A_I = \Tbf_{\widehat{p}}f_I$, where the last equality comes from the fact the $A_I$ is affine in $y$. Then Lemma \ref{lemma:theta_is_TpfI} concludes the proof.
\ep\\

We now use Assumption \ref{ass:domination} (ii) to prove the existence of a family $(A_K)_K$ satisfying the conditions of Lemma \ref{lemma:existenceTpf}. Let $\Cc\subset \Kcirc$, $\Dc\subset \Kcirc$, and $\Rc\subset \Kcirc$ from Assumption \ref{ass:domination} such that $I(X)\in\Cc\cup\Dc\cup\Rc$, $\mu-$a.s. with $\Cc$ well ordered, $\dim(\Dc)\subset\{0,1\}$, and $\cup_{K\neq K'\in\Rc}\big[K\x (\cl K\cap\cl K')\big]\in\Nc_{\mu,\nu}$.

\begin{Lemma}\label{lemma:TKK}
We assume Assumption \ref{ass:domination}, and the existence of $(T_K^{K'})_{K, K'\in \Cc\cup\Rc}\subset \Aff(\R^d,\R)$ such that


\no{\rm (i)} $T_K^{K'}+T_{K'}^{K''} +T^{K}_{K''}=0$, for all $K, K', K''\in \Cc\cup\Rc$;

\no{\rm (ii)} $T_K^{K'}(y) = f_{K'}(y) - f_K(y)$, for all $y\in\interf(K,K')$, $K, K'\in \Cc\cup\Rc$.

Then we may find $(A_K)_{K\in I(\R^d)}$ satisfying the conditions of Lemma \ref{lemma:existenceTpf}.
\end{Lemma}
\proof
We define $A_K$ by for $K\in\Cc\cup\Rc$. If this set is non-empty, we fix $K_0\in\Cc\cup\Rc$. Let $K\in\Cc$, we set $A_K := -T_{K_0}^K$.

Now for $K\in\Dc$, $K$ has at most two end-points, let $x\in\Jo(m_{K})$ be an end-point of $K$. If $x\in \Jo(m_{K'})$ for some $K'\in \Cc\cup\Rc$, then we set $A_K(x) := A_{K'}(x)+f_{K'}(x)-f_{K}(x)$. If $x\in \Jo(m_{K'})$ for some $K'\in \Dc$, then we set $A_K(x) := -f_{K}(x)$. Otherwise, we set $A_K(x):= 0$, and set $A_K$ to be the only affine function on $K$ that has the right values at the endpoints, and has a derivative orthogonal to $K$, which exists as $K$ is at most one-dimensional.

We define $A_K = 0$ for all the remaining $K\in I(\R^d)$.

Now we check that $(A_K)_K$ satisfies the right conditions at the interfaces. Let $K,K'\in I(\R^d)$ such that $\interf(K,K')\neq\emptyset$. If $K\in\Dc$, or $K'\in\Dc$, the value at endpoints has been adapted to get the desired value. Now we treat the remaining case, we assume that $K,K'\in \Cc\cup\Rc$. We have $A_K-A_{K'} = -T_{K_0}^K+T_{K_0}^{K'}$. Property (i) applied to $(K,K,K)$ implies that $T_K^K = 0$, and therefore, (i) applied to $(K_0,K,K)$ gives that $T_{K_0}^K = T^{K_0}_K$. Finally, (i) applied to $(K,K_0,K')$ gives that $A_K-A_{K'} = T_{K}^{K'}$. Finally, by (iii), we get that $A_K-A_{K'} = f_{K'}(y) - f_K(y)$ for all $y\in \interf(K,K')$.
\ep 

\begin{Lemma}\label{lemma:affineoninterf}
Let $K,K'\in I(\R^d)$, we have that $f_{K'} - f_K$ is affine finite on $\interf(K,K')$.
\end{Lemma}
\proof
First, by the fact that $\interf(K,K')\subset \dom\theta(m_K,\cdot)\cap\dom\theta(m_{K'},\cdot)$, $a := f_{K'} - f_K$ is finite on $\interf(K,K')$. Now we prove that this map is affine, let $y_1,...,y_k,y'_1,...,y'_{k'}\in \interf(K,K')$ such that $y = \sum_i\lambda_iy_i = \sum_i\lambda'_iy'_i$, convex combinations. Then $\P := \frac12\sum_i\lambda_i\delta_{(m_{K},y_i)}+\frac12\sum_i\lambda'_i\delta_{(m_{K'},y'_i)}$, and $\P' := \frac12\sum_i\lambda_i\delta_{(m_{K'},y_i)}+\frac12\sum_i\lambda'_i\delta_{(m_{K},y'_i)}$ are competitors that are concentrated on $\{Y\in\Jo(X),X\notin N_\mu\}\subset N^c$ by Lemma \ref{lemma:Jo_in_N}. Therefore, by Definition \ref{def:Thetamunu} (ii) we have $\sum_i\lambda_i\theta(m_{K},y_i)+\sum_i\lambda'_i\theta(m_{K'},y'_i) = \sum_i\lambda_i\theta(m_{K'},y_i)+\sum_i\lambda'_i\theta(m_{K},y'_i)$, which gives
$$\sum_i\lambda_ia(y_i)=\sum_i\lambda'_ia(y'_i).$$
Similar to the proof of Lemma \ref{lemma:recomposition}, we have that $a$ is affine on $\interf(K,K')$.
\ep\\

Let $K, K'\in I(\R^d)$, by the preceding lemma $f_{K'} - f_K$ is affine finite on $\interf(K,K')$. If this set is not empty, let the unique $a_K^{K'}\in\nabla \aff\,\interf(K,K')$ and $b_K^{K'}\in\R$ such that
\b*
f_{K'}(y) - f_K(y) = a_K^{K'}\cdot y+b_K^{K'},&\mbox{for}&y\in\interf(K,K').
\e*
We denote $H_K^{K'}: y\longmapsto  a_K^{K'}\cdot y+b_K^{K'}\in\Aff(\R^d,\R)$. If $\interf(K,K')= \emptyset$, we set $H_K^{K'}:=0$.

\begin{Lemma}\label{lemma:HtoT}
We may find $(T_K^{K'})_{K, K'\in \Cc\cup\Rc}\subset \Aff(\R^d,\R)$ satisfying $(i)$, and $(ii)$ from Lemma \ref{lemma:TKK} if and only if we may find $(\Hb_K^{K'})_{K, K'\in \Cc\cup \Rc}\subset \Aff(\R^d,\R)$ such that $\Hb_K^{K'} = 0$ on $\interf(K,K')$ for all $K,K'\in \Cc\cup \Rc$, and for all triplet $(K_i)_{i = 1,2,3}\in (\Cc\cup\Rc)^3$ such that with the convention $K_4 = K_1$, we have
\be\label{eq:sum_sigma}
\sum_{i = 1}^3 H_{K_i}^{K_{i+1}}+\Hb_{K_i}^{K_{i+1}} = 0.
\ee
\end{Lemma}
\proof
We start with the necessary condition, let $(T_K^{K'})_{K, K'\in \Cc\cup\Rc}\subset \Aff(\R^d,\R)$ satisfying $(i)$, and $(ii)$ from Lemma \ref{lemma:TKK}. Then for $K,K'\in \Cc\cup\Rc$, we introduce $\Hb_K^{K'} := T_K^{K'}-H_K^{K'}$. By (ii), together with the definition of $H_K^{K'}$, we have $\Hb_K^{K'} = 0$ on $\interf(K,K')$. Now let a finite $(K,K',K'')\subset \Cc\cup\Rc$, by (ii) we have
$$H_{K}^{K'}+\Hb_{K}^{K'}+H_{K'}^{K''}+\Hb_{K'}^{K''}+H_{K''}^{K}+\Hb_{K''}^{K} =  T_{K}^{K'}+T_{K'}^{K''}+T_{K''}^{K} = 0.$$

Now we prove the sufficiency. Let $(\Hb_K^{K'})_{K, K'\in \Cc\cup \Rc}\subset \Aff(\R^d,\R)$ such that $\Hb_K^{K'} = 0$ on $\interf(K,K')$ for all $K,K'\in \Cc\cup \Rc$, and for all finite set $\Fc\subset \Cc\cup\Rc$, and all triplet $(K_i)_{i = 1,2,3}\in (\Cc\cup\Rc)^3$ such that we have $\sum_{i = 1}^3 H_{K_i}^{K_{i+1}}+\Hb_{K_i}^{K_{i+1}} = 0$.



Then for $K,K'\in\Cc\cup\Rc$, let $T_K^{K'} := H_K^{K'}+\Hb_K^{K'}$. The property (ii) of $(T_K^{K'})$ follows from the fact that $T_K^{K'}=H_K^{K'}+\Hb_K^{K'} = H_K^{K'} = f_{K'}-f_K$ on $\interf(K,K')$.

Property (i) is a direct consequence of \eqref{eq:sum_sigma} with $(K,K',K'')\in(\Cc\cup\Rc)^3$.
\ep

\begin{Lemma}\label{lemma:H_bar_finite_F}
Let $\Fc\subset I(N_\mu^c)$ finite, we may find $(\Hb_K^{K'})_{K, K'\in \Fc}\subset \Aff(\R^d,\R)$ such that $\Hb_K^{K'} = 0$ on $\interf(K,K')$ for all $K,K'\in \Fc$, and for all triplet $(K_i)_{i = 1,2,3}\in \Fc^3$ such that with the convention $K_4 = K_1$, we have $\sum_{i = 1}^3 H_{K_i}^{K_{i+1}}+\Hb_{K_i}^{K_{i+1}} = 0$.
\end{Lemma}
\proof
Let $p\in\Fc^2$, we denote $H_p := H_{p_1}^{p_2}$, $\interf(p) := \interf(p_1,p_2)$, and the linear map $g_p:A\in\Aff(\R^d,\R)\longmapsto A_{|\aff\,\interf(p)}\in\Aff(\aff\,\interf(p),\R)$. Let the linear map
$$g:(A_p)_{p\in\Fc^2}\in\Aff(\R^d,\R)^{\Fc^2}\longmapsto \big(g_p(A_p)\big)_{p\in\Fc^2}\in\bigtimes_{p\in\Fc^2}\Aff(\aff\,\interf(p),\R),$$
and if we denote $t_{i,j} := (t_i,t_j)\in\Fc^2$ for $t\in\Fc^3$ and $i,j\in\{1,2,3\}$, let the other linear map
$$f:(A_p)_{p\in\Fc^2}\in\Aff(\R^d,\R)^{\Fc^2}\longmapsto\big(A_{t_{1,2}}+A_{t_{2,3}}+A_{t_{3,1}}\big)_{t\in\Fc^3}\in\Aff(\R^d,\R)^{\Fc^3}.$$
Notice that the result may be written in terms of $f$ and $g$ as 
\be\label{eq:fH_in_Im}
f\big((H_p)_{p\in\Fc^2}\big)\in f(\ker g).
\ee
We prove this statement by using the monotonicity principle (ii) of Definition \ref{def:Thetamunu}. Let the canonical basis $(e_j)_{1\le j\le d}$ of $\R^d$, and $e_0:= 0$ so that $(e_j)_{0\le j \le d}$ is an affine basis of $\R^d$, and the scalar product on $\Aff(\R^d,\R)^{\Fc^3}$ defined by $\big\langle(A_t)_{t\in\Fc^3},(A'_t)_{t\in\Fc^3}\big\rangle:=\sum_{t\in\Fc^3,0\le j \le d}A_t(e_j)A'_t(e_j)$. As the dimensions are finite, \eqref{eq:fH_in_Im} is equivalent with the inclusion $f(\ker g)^\perp\subset \left\{f\big((H_p)_{p\in\Fc^2}\big)\right\}^\perp$.

Let $(A_t)_{t\in\Fc^3}\in f(\ker g)^\perp$, we now prove that $(A_t)_{t\in\Fc^3}\in\left\{f\big((H_p)_{p\in\Fc^2}\big)\right\}^\perp$, i.e. that
\b*
\left\langle(A_t)_{t\in\Fc^3},f\big((H_p)_{p\in\Fc^2}\big)\right\rangle &=&\sum_{t\in\Fc^3,0\le j \le d}A_t(e_j)\left(H_{t_{1,2}}(e_j)+H_{t_{2,3}}(e_j)+H_{t_{3,1}}(e_j)\right)\\
&=&0.
\e*
Let $p\in\Fc^2$, $\mathbf{P}_p:=\proj_{\aff\,\interf(p)}$, and $0\le j\le d$. By the fact that $\mathbf{P}_p(e_j)\in\aff\,\interf(p) = \aff\big(J(m_{p_1})\cap J(m_{p_2})\setminus N_\nu\big)$ by Remark \ref{rmk:J}, we may find $(y_{i,j,p})_{1\le i\le d+1}\subset J(m_{p_1})\cap J(m_{p_2})\setminus N_\nu$, and $(\lambda_{i,j,p})_{1\le i\le d+1}\subset \R$ such that $\mathbf{P}_p(e_j) = \sum_{i=1}^{d+1}\lambda_{i,j,p}y_{i,j,p}$, affine combination, and $\sum_{i=1}^{d+1}\lambda_{i,j,p} = 1$. Then with these ingredients we may give the expression of $H_p(e_j)$ as a function of values of $\theta$:
\b*
H_p(e_j) =\sum_{i=1}^{d+1}\lambda_{i,j,p}H_p(y_{i,j,p})&=&\sum_{i=1}^{d+1}\lambda_{i,j,p}\left[\theta(m_{p_2}, y_{i,j,p})-\theta(m_{p_1}, y_{i,j,p})\right]\\
&=& L_p^j[\theta],
\e*
where $L_p^j := \sum_{i=1}^{d+1}\lambda_{i,j,p}\left[\delta_{(m_{p_2}, y_{i,j,p})}-\delta_{(m_{p_1}, y_{i,j,p})}\right]$ is a signed measure with finite support in $\{Y\in J(X)\setminus N_\nu,X\notin N_\mu\}$. We now study the marginals of $L_p^j$: we have obviously from its definition that $L_p^j[Y = y] = 0$ for all $y\in\R^d$. For the X-marginals, $L_p^j[X = m_{p_2}] = -L_p^j[X = m_{p_1}] = \sum_{i=1}^{d+1}\lambda_{i,j,p} = 1$, and $L_p^j[X=x] = 0$ for all other $x\in\R^d$. Finally we look at its conditional barycenter:
\be\label{eq:cond_exp}
L_p^j[Y|X = m_{p_2}] = -L_p^j[Y|X = m_{p_1}] = \sum_{i=1}^{d+1}\lambda_{i,j,p}y_{i,j,p} = \mathbf{P}_p(e_j).
\ee
Now let $t\in\Fc^3$, we denote $L_t^j:= L_{t_{1,2}}^j+L_{t_{2,3}}^j+L_{t_{3,1}}^j$. We still have $L_t^j[Y = y] = 0$ for all $y\in\R^d$ by linearity. Now
\b*
L_t^j[X = t_1] &=& L_{t_{1,2}}^j[X = t_1]+L_{t_{2,3}}^j[X = t_1]+L_{t_{3,1}}^j[X = t_1]\\
&=& -\mathbf{1}_{t_1 = t_1}+\mathbf{1}_{t_2 = t_1}-\mathbf{1}_{t_2 = t_1}+\mathbf{1}_{t_3 = t_1}-\mathbf{1}_{t_3 = t_1}+\mathbf{1}_{t_3 = t_1}\\
&=& 0.
\e*
Similar, $L_t^j[X = t_2] = L_t^j[X = t_3] = 0$, and $L_t^j[X=x] = 0$ for all $x\in\R^d$.

Notice that $\left\langle(A_t)_{t\in\Fc^3},f\big((H_p)_{p\in\Fc^2}\big)\right\rangle = L[\theta]$, with $L:=\sum_{t\in\Fc^3,0\le j \le d}A_t(e_j)L_t^j$. By linearity, we have that
\be\label{eq:marginals}
L[X=x] = L[Y=x] = 0,&\mbox{for all}&x\in\R^d.
\ee
Furthermore, $L$ is supported on $\{Y\in J(X)\setminus N_\nu,X\notin N_\mu\}\subset N^c$ like each $L_p^j$. We claim that $L[Y|X] = 0$, this claim will be justified at the end of this proof. Then we consider the Jordan decomposition $L = L_+-L_-$ with $L_+$ the positive part of $L$ and $L_-$ its negative part. By the fact that $L[\R^d] = 0$, we have the decomposition $L = C(\P_+-\P_-)$, for $C = L_+[R^d] = - L_-[R^d]\ge 0$. Then $\P_+$ and $\P_-$ are two finitely supported probabilities concentrated on $N^c$. By the fact that $L[X] = L[Y] = L[Y|X] = 0$, $\P_+$, and $\P_-$ are furthermore competitors, then by Definition \ref{def:Thetamunu} (ii), $\P_+[\theta] = \P_-[\theta]$, and therefore $\left\langle(A_t)_{t\in\Fc^3},f\big((H_p)_{p\in\Fc^2}\big)\right\rangle = L[\theta] = 0$, which concludes the proof.

It remains to prove the claim that $L[Y|X] = 0$. Recall that $(A_t)_{t\in\Fc^3}\in f(\ker g)^\perp$. Let $K\in\Fc$ and $p\in\Fc^2$ such that $p_1 = K$, and $u\in\R^d$, the map $\xi_p:x\longmapsto u\cdot\big(x - \mathbf{P}_p(x)\big)$ is in $\ker g_p$. For all the other $p'\in\Fc^2$, we set $\xi_{p'} := 0\in\ker g_{p'}$. Then $(\xi_p)_{p\in\Fc^2}\in\ker g$, and therefore $\left\langle(A_t)_{t\in\Fc^3}, f\big((\xi_p)_{p\in\Fc^2}\big)\right\rangle = 0$, we have
\b*
0&=&\sum_{t\in\Fc^3,0\le j \le d}A_t(e_j)\sum_{p=t_{1,2},t_{2,3},t_{3,1}}\mathbf{1}_{p_1 = K}u\cdot\big(e_j - \mathbf{P}_p(e_j)\big)\\
&=& u\cdot \sum_{t\in\Fc^3,0\le j \le d}A_t(e_j)\sum_{p=t_{1,2},t_{2,3},t_{3,1}}\mathbf{1}_{p_1 = K}\big(e_j - \mathbf{P}_p(e_j)\big).
\e*
As this holds for all $u\in\R^d$, we have $\sum_{t\in\Fc^3,0\le j \le d}A_t(e_j)\sum_{p=t_{1,2},t_{2,3},t_{3,1}}\mathbf{1}_{p_1 = K}\big(e_j - \mathbf{P}_p(e_j)\big)=0$. Similarly, we have $\sum_{t\in\Fc^3,0\le j \le d}A_t(e_j)\sum_{p=t_{1,2},t_{2,3},t_{3,1}}\mathbf{1}_{p_2 = K}\big(e_j - \mathbf{P}_p(e_j)\big)=0$. Combining these two equations, and using \eqref{eq:cond_exp} together with the definition of $L$ we get
\b*
L[Y|X = m_K] &=&\sum_{t\in\Fc^3,0\le j \le d}A_t(e_j)\sum_{p=t_{1,2},t_{2,3},t_{3,1}}(\mathbf{1}_{p_2 = K}-\mathbf{1}_{p_1 = K})\mathbf{P}_p(e_j)\\
&=&\sum_{t\in\Fc^3,0\le j \le d}A_t(e_j)\sum_{p=t_{1,2},t_{2,3},t_{3,1}}(\mathbf{1}_{p_2 = K}-\mathbf{1}_{p_1 = K})e_j\\
&=& L[X = m_K]e_j = 0,
\e*
by \eqref{eq:marginals} together with the definition of $L$. We conclude that $L[Y|X= m_K] = 0$, the claim is proved.
\ep 

\begin{Lemma}\label{lemma:globalexistenceH}
Under Assumption \ref{ass:domination}, we may find $(\Hb_K^{K'})_{K, K'\in \Cc\cup\Rc}\subset \Aff(\R^d,\R)$ such that $\Hb_K^{K'} = 0$ on $\interf(K,K')$ for all $K,K'\in \Cc\cup\Rc$, and for all triplet $(K_i)_{i = 1,2,3}\in (\Cc\cup\Rc)^3$ such that with the convention $K_4 = K_1$, we have $\sum_{i = 1}^3 H_{K_i}^{K_{i+1}}+\Hb_{K_i}^{K_{i+1}} = 0$.
\end{Lemma}
\proof
We use the well-order of $\Cc$ from Assumption \ref{ass:domination} to extend the result of Lemma \ref{lemma:H_bar_finite_F} to the possibly infinite number of components.
By the fact that $\Cc$ is well ordered, we have that $\Cc^2$ is also well ordered (we may use for example the lexicographic order based on the well-order of $\Cc$). We shall argue by transfinite induction on $\Cc^2$. For $(K,K')\in\Cc^2$, we denote $\Cc(K,K'):=\{(K_1,K_2)\in\Cc^2:(K_1,K_2)<(K,K')\}$. Finally we fix $\|\cdot\|$, a euclidean norm on the finite dimensional space $\Aff(\R^d,\R)$, and
for $(K,K')\in\Cc^2$, we define an order relation $\preceq_{K,K'}$ on $\Aff(\R^d,\R)^{\Cc(K,K')}$ which is the lexicographical order induced by $(\Cc(K,K'),\le)$, and by the order on affine function $(\Aff(\R^d,\R),\preceq)$, defined by $A\preceq A'$ if $\|A\|\le \|A'\|$. Our induction hypothesis is:

$\Hc(K,K'):$ we may find a unique $(\Hb_{K_1}^{K_2})_{(K_1,K_2)\in\Cc(K)}$ such that:

{\rm (i)} for all finite $\Fc\subset \Cc\cup \Rc$, we may find $(\Ht_{K_1}^{K_2})_{K_1, K_2\in \Fc}\subset \Aff(\R^d,\R)$ such that $\Ht_{K_1}^{K_2} = 0$ on $\interf(K_1,K_2)$ for all $K_1,K_2\in \Fc$, such that for all triplet $(K_i)\in \Fc^3$ we have $\sum_{i = 1}^3 H_{K_i}^{K_{i+1}}+\Ht_{K_i}^{K_{i+1}} = 0$, and finally such that $\Hb_{K_1}^{K_2} = \Ht_{K_1}^{K_2}$ for all $(K_1,K_2)\in \Fc^2\cap \Cc(K,K')$;

{\rm (ii)} for all $(K'',K''')\le (K,K')$, $(\Hb_{K_1}^{K_2})_{K_1,K_2\in\Cc(K'',K''')}$ is the minimal vector satisfying (i) of $\Hc(K'',K''')$, for the order $\preceq_{K'',K'''}$.

Similar to the ordinals, we consider $\Cc^2$ as the upper bound of all the elements it contains, which gives a meaning to $\Hc(\Cc^2)$. The transfinite induction works similarly to a classical structural induction: let $(K_0,K_0')\in\Cc^2$ be the smallest element of $\Cc$, then the fact that $\Hc( K_0,K_0')$ holds, together with the fact that for all $ (K,K')\in\Cc$, we have that $\Hc( K'',K''')$ holding for all $ (K'',K''')< (K,K')$ implies that $\Hc( K,K')$ holds, then the transfinite induction principle implies that $\Hc(\Cc^2)$ holds.

The initialization is a direct consequence of Lemma \ref{lemma:H_bar_finite_F} as $\Cc(K_0,K_0') = \emptyset$. Now let $ (K,K')\in\Cc$, we assume that $\Hc( K',K'')$ holds for all $ (K'',K''')< (K,K')$. Let $ (K_1,K_1')< (K_2,K_2')< (K,K')$. As $\Hc( K_1,K_1')$, and $\Hc( K_2,K_2')$ hold, we may find unique $(\Hb^{1,K''}_{ K'})_{ K',K''\in \Cc(K_1,K_1')}$, and $(\Hb^{2,K''}_{ K'})_{ K',K''\in\Cc(K_2,K_2')}$ satisfying the conditions of the induction hypothesis. The restriction $(\Hb^{2,K''}_{ K'})_{ K',K''\in \Cc(K_1,K_1')}$ satisfies the conditions of $\Hc( K_1,K_1')$ by $\Hc( K_2,K_2')$, and by the fact that for the lexicographic order, if a word is minimal then all its prefixes are minimal as well for the sub-lexicographic orders. Therefore, by uniqueness in $\Hc( K_1,K_1')$, $(\Hb^{1,K''}_{ K'})_{ K',K''\in\Cc( K_1,K_1')}=(\Hb^{2,K''}_{ K'})_{ K',K''\in\Cc( K_1,K_1')}$. For all $ (K'',K''')< (K,K')$ which are not predecessors of $ (K,K')$ (i.e. such that we may find $ (K_{int},K_{int}')\in\Cc^2$ with $ (K'',K''')< (K_{int},K_{int}')< (K,K')$), let $\Hb_{ K'}^{K''}$ be the $ (K',K'')-$th affine function of $(\Hb_{ K_1}^{K_2})_{ K_1,K_2\in\Cc(K_{int},K_{int}')}$ satisfying $\Hc(K_{int},K_{int}')$, which is unique by the preceding reasoning. If $(K,K')$ has no predecessor, then $\Hb:=(\Hb_{ K'}^{K'''})_{ K'',K'''\in\Cc( K,K')}$ solves $\Hc( K,K')$. Now we treat the case in which we may find a predecessor $ (K_{pred},K_{pred}')\in\Cc^2$ to $ (K,K')$. In this case this predecessor is unique because $\Cc^2$ is well ordered. Then we consider $\Hb:=(\Hb_{ K_1}^{K_2})_{ K_1,K_2\in\Cc( K_{pred},K_{pred}')}$ from $\Hc(K_{pred},K_{pred}')$. Now we need to complete $\Hb$ by defining $\Hb^{K_{pred}'}_{K_{pred}}$.






For all finite $\Fc\subset\Cc\cup\Rc$, we define the affine subset $A_\Fc\subset \Aff(\R^d,\R)$ of all $H\in\Aff(\R^d,\R)$ such that $(\Ht_{K_1}^{K_2})_{K_1,K_2\in\Cc(K,K')}$ satisfies (i) of $\Hc(K,K')$, with $\Ht_{K_1}^{K_2} = \Hb_{K_1}^{K_2}$ if $(K_1,K_2)<(K_{pred},K_{pred}')$ and $\Ht_{K_{pred}}^{K_{pred}'}$. By $\Hc(K_{pred},K_{pred}')$ (i) applied to $\Fc\cup\{(K_{pred},K_{pred}')\}$, we have that $A_\Fc$ is non-empty for all $\Fc$. Then the intersection taken on finite sets $A:=\cap_{\Fc\subset\Cc\cup\Rc}A_\Fc$ is also non-empty as we intersect finite dimensional always non-empty affine spaces that have the property $A_{\Fc_1}\cap A_{\Fc_2} = A_{\Fc_1\cup\Fc_2}$.
Then if we chose $H_{K_{pred}}^{K_{pred}'}\in A$, $\Hc(K_{pred},K_{pred}')$ will be verified, except for the minimality.
To have the minimality, we chose the minimal $H\in A$ for the norm $\|\cdot\|$, which is unique as $A$ is affine and the norm is Euclidean.
This uniqueness, together with the uniqueness from the induction hypothesis gives the uniqueness for $\Hc_(K_{pred},K_{pred}')$ by properties of the lexicographic order.
We proved $\Hc_(K_{pred},K_{pred}')$, and therefore $\Hc(\Cc^2)$ holds.


Finally, we need to include $\Rc$ in the indices of $\Hb$. Let the unique $(\Hb_{K}^{K'})_{(K,K')\in\Cc^2}$ from $\Hc(\Cc^2)$. Let $K\in\Rc$, $K'\in\Cc$. Similar to the step in the induction $(K_{pred},K_{pred}')$ to $(K,K')$, we may find a unique $\Hb_{K}^{K'}$ which satisfies the right relations and is minimal for the norm $\|\cdot\|$. As we may do it independently for all $(K,K')\in\Rc\x\Cc$ by the property of $\Rc$ in Assumption \ref{ass:domination}. For $K\in\Cc$ and $K'\in\Rc$, we set $\Hb_{K}^{K'} := -\Hb_{K'}^{K}$. Finally for $K,K'\in\Rc$, if $\Cc = \emptyset$, then we set $\Hb_{K}^{K'} := 0$, else we set $\Hb_{K}^{K'} := \Hb_{K_0}^{K'} - \Hb_{K_0}^{K}$ for some $K_0\in\Cc$. We may prove thanks to $\Hc(\Cc^2)$ that this definition does not depend on the choice of $K_0\in\Cc$, and that $\Hb$ defined this way on $(\Cc\cup\Rc)^2$ satisfies the right conditions.
\ep\\

\no {\bf Proof of Proposition \ref{prop:CctoTc}} The inclusion $\supset$ is obvious from the definition of $\partial^{\mu,\nu}f$. We now prove the reverse inclusion by using Assumption \ref{ass:domination}. Then by Lemma \ref{lemma:globalexistenceH}, we may find $(\Hb_K^{K'})_{K\sim_1 K'\in \Cc\cup \Rc}\subset \Aff(\R^d,\R)$ such that for all finite set $\Fc\subset \Cc\cup\Rc$, and all permutation $\sigma\in\Sc_\Fc$ such that $K\sim_1\sigma(K)$ for all $K\in \Fc$, we have $\sum_{K\in\Fc} H_K^{\sigma(K)}+\Hb_K^{\sigma(K)} = 0$. Then, by Lemma \ref{lemma:HtoT}, we may find $(T_K^{K'})_{K,K'\in I(\R^d):K\sim K'}\subset \Aff(\R^d,\R)$ satisfying $(i)$, $(ii)$, and $(iii)$ from Lemma \ref{lemma:TKK}. Then we may apply Lemma \ref{lemma:TKK}: we may find $(A_K)_{K=I(x), x\in\R^d}\subset \Aff(\R^d,\R)$ such that $A_K(y) - A_{K'}(y) = f_{K'}(y) - f_K(y)$ for all $y\in \interf(K,K')$, and for all $K,K'\in I(\R^d)$. Finally, by Lemma \ref{lemma:existenceTpf}, $f(y) := f_K(y) + A_K(y)$ does not depend of the choice of $K$ such that $y\in \Jo(m_K)$, and if we set $p(y):= \widehat{p}(y)+\nabla A_{I(y)}$, we have
\b*
\theta = \Tbf_pf,&\mbox{on}&\{X\notin N_\mu\}\cap\{Y\in \Jo(X)\}.
\e*
Therefore, $\theta \approx \Tbf_pf$, whence $f\in\Cc_{\mu,\nu}$ and we proved the reverse inclusion.
\ep\\

Now, we prove the convexity of the functions in $\Cfrak_{\mu,\nu}$ on each components.

\no {\bf Proof of Proposition \ref{prop:Mmunuconvex}}
Let $p\in\partial^{\mu,\nu}f$, and $\theta\in\widetilde{\Tc}(\mu,\nu)$ such that $\Tbf_pf=\theta$ on $\{Y\in\Jo(X),X\notin N_\mu\}$ for a $N-$tangent convex function $\theta\in\widetilde{\Tc}(\mu,\nu)$, $N_\mu\in\Nc_\mu$, and $\Jo(\mu,\nu)$. By proposition \ref{prop:J}, we may chose $N_\mu$ and $\Jo$ such that $\{Y\in\Jo(X),X\notin N_\mu\}\subset\dom\theta\cap N^c$. For all $x\notin N_\mu$, and $y\in \Jo(x)$, $f(y) = f(x)+p(x)\cdot(y-x)+\theta(x,y)$, which is clearly convex in $y$ for $x$ fixed. The function $f$ is convex on $\Jo$, $\eta-$a.s.

For all $x\in N_\mu^c$ and $y\in \Jo(x)$, we have $f(y)-f(x)-\proj_{\nabla\aff \Jo}(p)(x)\cdot(y-x) = f(y)-f(x)-p(x)\cdot(y-x)=\theta(x,y)\ge 0$. Then by definition, $\proj_{\nabla\aff \Jo}(p)(x)\in \partial f|_{\Jo}(x)$ for all $x\notin N_\mu$.

For $x\in N_\mu^c$, we define $\widetilde f:= (f \mathbf{1}_{\Jo(x)})_{conv}$ on $J(x)=\conv\big(\Jo(x)\big)$, where the equality comes from Proposition \ref{prop:J} (i) together with the fact that $J\setminus N_\nu\subset \Jo$. We also define $\widetilde f:=f$ on $\cap_{x\in N_\mu^c}\Jo(x)^c\in\Nc_\mu\cap\Nc_\nu = \Nc_{\mu+\nu}$. These definitions are not interfering as if $x'\in J(x)$ then $J(x')\subset J(x)$ by Remark \ref{rmk:J}. Therefore, the convex envelops $(f \mathbf{1}_{\Jo(x)})_{conv}$ and $(f \mathbf{1}_{\Jo(x')})_{conv}$ coincide on $J(x')$.

Then the map $p(X)\cdot(Y-X) = f(Y)-f(X)-\theta(X,Y)$ is Borel measurable on $I(x)\x I(x)$ for all $x\in N_\mu^c$. Let $x\notin N_\mu$, $d_x:=\dim I(x)$, and $\big(y_i\big)_{1\le i\le d_x+1}\in I(x)$, affine basis of $\aff I(x)$. Therefore, $\proj_{\nabla\aff I(x)}\big(p(x')\big) = M^{-1}\Big(p(x')\cdot\big(y_i-y_{d_x+1}\big)\Big)_{1\le i\le d_x}$, with $M:=\big(y_i-y_{d_x+1}\big)_{1\le i\le d_x}$, where everything is expressed in the basis $(y_i-y_{d_x+1})_{1\le i\le d_x}$, is Borel measurable on $I(x)$. Then as it is a subgradient of $f_{|I(x)}$ on $I(x)$ by the fact that $\theta(x,y)=f(y)-f(x)-\proj_{\nabla\aff I(x)}\big(p(x)\big)\cdot (y-x)\ge 0$ for all $x,y\in I(x)$, we have the result.

Finally, notice that $\Tbf_{\widetilde{p}}\widetilde{f} = \Tbf_pf = \theta$ on $\{Y\in\Jo(X),X\notin N_\mu\}$, which proves that $\widetilde{f}\in\Cfrak_{\mu,\nu}$ and $\widetilde{p}\in\partial^{\mu,\nu}\widetilde{f}$.
\ep\\

\no {\bf Proof of Proposition \ref{prop:measurability}}
\no{\rm (i)} Let $(\varphi,\psi,h)\in\widehat\Leb(\mu,\nu)$, and let $f$ be its q.s.-convex moderator, and $p\in\partial^{\mu,\nu}f$. By Proposition \ref{prop:Mmunuconvex}, $f$ is convex and finite on $I$, and $\proj_{\nabla\aff I}(p)\in\partial f_{|I}$, $\eta-$a.s. Then $\psi=(\psi-f)+f$ is Borel measurable on $I$, $\varphi=(\varphi+f)-f$ is Borel measurable on $I$, and $\proj_{\nabla\aff J}(h)=\proj_{\nabla\aff J}(h+p)-\proj_{\nabla\aff J}(p)$ is Borel measurable on $I$, $\eta-$a.s.

\no{\rm (ii)} If one of the conditions in Proposition \ref{prop:emptyboundaries} holds, then condition (iv) holds by Proposition \ref{prop:emptyboundaries}. Then the transfinite induction from the proof of Proposition \ref{prop:CctoTc} becomes a countable induction, thus preserving the measurability. The process of subtracting lines for the one dimensional components is also measurable.
\ep

\subsection{Consequences of the regularity of the cost in x}\label{subsect:regularity}

\no {\bf Proof of Lemma \ref{lemma:regularity}}
We have for all $x,y\in\R^d$, $\varphi(x)+\psi(y)+h(x)\cdot(y-x)\ge c(x,y)$. Then $\varphi(x)\ge\varphi'(x):= -(\psi-c(x,\cdot))_{conv}(x)$. For all $x\in\R^d$, $f_x:=(\psi-c(x,\cdot))_{conv}$ is convex and finite on $D:=\ri\,\conv\,\dom\,\varphi$, let $-h':\R^d\longmapsto\R^d$ be a measurable selection in its subgradient on $D$ (then in $\aff D - x_0$ for some $x_0\in D$). Then for all $y\in\R^d$,
$$-\varphi'(x)-h'(x)\cdot(y-x)\le f_x(y):=(\psi-c(x,\cdot))_{conv}(y)\le \psi(y)-c(x,y).$$
Then $c\le \varphi'\oplus\psi+h'^\otimes$, and therefore, $\P[\varphi'\oplus\psi+h'^\otimes]\ge \P[c]$ is well defined. Subtracting $\P[\varphi\oplus\psi+h^\otimes]<\infty$, we get
$$\mu[\varphi'-\varphi]=\P[(\varphi'-\varphi)(X)+(h'-h)^\otimes]\ge \P[c]-\P[\varphi\oplus\psi+h^\otimes].$$
Finally, taking the supremum over $\P$, we get $\mu[\varphi'-\varphi]\ge \Sbf_{\mu,\nu}(c)-\Sbf_{\mu,\nu}(\varphi\oplus\psi+h^\otimes)=0$. As $\varphi'-\varphi\le 0$, this shows that $\varphi'=\varphi$, $\mu-$a.e. Now
\be\label{eq:phivalue}
f_x(y) = -\inf\left\{\sum_{i=1}^r\lambda_i\big(\psi(y_i)-c(x,y_i)\big):\sum_{i=1}^r\lambda_i y_i = y,\mbox{ and } r\ge 1\right\}
\ee
For $r\ge 1$, and $y=\sum_{i=1}^r\lambda_i y_i$, $x\longmapsto \sum_{i=1}^r\lambda_i\big(\psi(y_i)-c(x,y_i)\big)$ is locally Lipschitz. By taking the infimum, we get that for $x\in D$, $f_{x}(y)$ is uniformly Lipschitz in $x$. Furthermore, $f_x$ is convex on the relative interior of its domain $D$, and therefore locally Lipschitz on it. We claim that for the convex function $f_x$, the Lipschitz constant on a compact $K\subset D$ is bounded by $\frac{\max_{K'}f_x-\min_Kf_x}{\delta}$, where $\delta = \inf_{(x,y)\in K\x K'}|x-y|$, for any compact $K'\subset D$ such that $K\subset\ri\,K'$ (cf proof of Theorem 9.3 in \cite{de2017irreducible}). Then if we fix $K$ and $K'$, the Lipschitz constant of $f_x$ is dominated on $K$ as $x\longmapsto(\max_K'f_x,\min_Kf_x)$ is Locally Lipschitz. Then for $K\subset D$ compact, we may find $L$, and $L'$, Lipschitz constants for both variables. Finally, for $x_1,x_2\in B$,
$$|\varphi'(x_1)-\varphi'(x_2)|\le |f_{x_1}(x_1)-f_{x_1}(x_2)|+|f_{x_1}(x_2)-f_{x_2}(x_2)|\le (L+L')|x_1-x_2|.$$
In the proof of Theorem 9.3 in \cite{de2017irreducible}, the bound $\frac{\max_K'f_x-\min_Kf_x}{\delta}$ is in fact a bound for the subgradients of $f_x$. As $-h'$ is a subgradient of $f_x$ in $x$, its component in $\aff D - x_0$ (for some $x_0\in D$) is bounded in $K$.

\section{Verification of Assumptions \ref{ass:domination}}\label{sect:assumption}
\setcounter{equation}{0}

\subsection{Marginals for which the assumption holds}\label{subsect:verifass}

In preparation to prove Proposition \ref{prop:emptyboundaries}, we first need to prove two lemmas.

\begin{Lemma}\label{lemma:dominproba}
Assume that there exists $\Q\in\Pc(\Omega)$ such that
\be\label{eq:convergence_implication}
(\theta_n)_{n\ge 1}\subset\widetilde{\Tc}_1,~\mbox{converges}~\Mc(\mu,\nu)-\mbox{q.s.}~\mbox{whenever}~(\theta_n)_{n\ge 1},~\mbox{converges}~\Q-\mbox{a.s.}
\ee
Then for all $(\theta_n)_{n\ge 1}\subset \widetilde\Tc_1$, we may find $\theta\in\widetilde\Tc_1$ such that $\theta_n\rightsquigarrow \theta$.
\end{Lemma}
\proof
Let $\Q\in\Pc(\Omega)$ satisfying \eqref{eq:convergence_implication}. Let $\Q':=\frac12\Q+\frac12\mu(dx)\otimes \sum_{n\ge 1}2^{-n}\delta_{f_n(x)}(dy)$, where $(f_n)_{n\ge 1}\subset\Leb^0(\R^d,\R^d)$ is chosen such that $\{f_n(x):n\ge 1\} \subset \aff I(x)$ is dense in $\aff I(x)$ for all $x\in\R^d$ (see Step 2 in the proof of Proposition 2.7 in \cite{de2017irreducible}). Then by Koml\'os lemma, we may find $\widehat\theta_n\in\conv(\theta_k,k\ge n)$ such that $\widehat\theta_n$ converges $\Q'-$a.s. Therefore, $\widehat\theta_n$ converges q.s. to $\theta:=\underline{\widehat\theta}_\infty$. As $\widehat\theta_n\in\conv(\theta_k,k\ge n)$, we have the inequality $\underline{\widehat\theta}_\infty\ge \underline{\theta}_\infty$. We also have by Fatou's lemma $\P[\underline{\widehat\theta}_\infty]\le\liminf_{n\to\infty}\P[\widehat\theta_n]\le\limsup_{n\to\infty}\P[\theta_n]$, for all $\P\in\Pc(\Omega)$. Finally we need to prove that $\theta\in\Theta_{\mu,\nu}$. For $n\ge 1$, let $N_n\in\Nc_{\mu,\nu}$ be the set from Definition \ref{def:Thetamunu} for $\theta_n$, and let $N_{cvg}\in\Nc_{\mu,\nu}$ be the set where $\widehat\theta_n$ does not converge. We set $ N:= \cup_{n\ge 1} N_n\cup N_{cvg}\in\Nc_{\mu,\nu}$. As $\theta_n(X,X) = 0$ for all $n\ge 1$, we have obviously $\{X=Y\}\subset N_{cvg}^c$, and $\{X=Y\}\subset N^c$. By convexity of $\theta_n(x,\cdot)$, the $\mu(dx)\otimes \sum_{n\ge 1}2^{-n}\delta_{f_n(x)}(dy)-$convergence implies pointwise convergence of $\theta(X,\cdot)$ on $I(X)$, $\mu-$a.s. as in the case of $\muxpw-$convergence. Then $\theta(x,\cdot)$ is convex on $ N_x^c$ by passing to the limit, $I(X)\subset N_X^c$, $\mu-$a.s. By Lemma 6.1 in \cite{de2017irreducible}, we may chose $N_\mu\in\Nc_\mu$ so that if $\Nc_\mu\ni N_\mu'\supset N_\mu$, then $\{Y\in I(X)\}\cap\{X\in N_\mu'^c\}$ is a Borel set, and therefore, the function $\mathbf{1}_{\{Y\in I(X)\}\cap\{X\in N_\mu'^c\}}\theta_{\infty}$ is Borel and Definition \ref{def:Thetamunu} (iv) holds.

For $\P$ with finite support on $N^c$, and $\P'$ competitor to $\P$, $\P[\theta]=\lim_{n\to\infty}\P[\widehat\theta_n]$, and $\P'[\theta]\le\liminf_{n\to\infty}\P'[\widehat\theta_n]$ by Fatou's Lemma. As for all $n\ge 1$, $\P[\widehat\theta_n]\ge \P'[\widehat\theta_n]$, we get the inequality $\P[\theta]\ge\P'[\theta]$. Furthermore, if we suppose to the contrary that $\{\omega\}:=\supp\, \P'\cap N$ is a singleton, $\omega\notin N_n$ for all $n\ge 1$ by Definition \ref{def:Thetamunu} (iii). Then for all $n\ge 1$, $\P[\theta_n]=\P'[\theta_n]$, and $\P'[\omega]\theta_n(\omega)=\P[\theta_n]-\P'[\theta_n\mathbf{1}_{\Omega\setminus\{\omega\}}]$. Then as the term on the right of this equality converges, $\theta_n(\omega)$ converges as well, and $\omega\in N^c$. We got the contradiction, (iii) of Definition \ref{def:Thetamunu} holds.
\ep

\begin{Lemma}\label{lemma:nocharge}
Assume that $\nu$ is dominated by the Lebesgue measure. Then $Y\notin \partial I(X)$ whenever $\dim I(X)\ge d-1$, $\Mc(\mu,\nu)-$q.s.
\end{Lemma}
\proof
First the components of dimension $d$ are at most countable, and their boundary is Lebesgue negligible as they are convex. Then, if we enumerate the countable $d-$dimensional components $(I_k)_{k\ge 1}$, we have $Y\notin \cup_{k\ge 1}\partial I_k$, $\nu-$a.s. and therefore $\Mc(\mu,\nu)-$q.s.

Now we deal with the $(d-1)-$dimensional components. $I$ is a Borel map, and therefore by Lusin theorem (see Theorem 1.14 in \cite{evans2015measure}), for all $\epsilon>0$, we may find $K_\epsilon\subset \{\dim I(X) = d-1\}$ with $\mu[K_\epsilon]\ge \mu[\dim I(X) = d-1]-\epsilon$, on which $I$ is continuous. We may also assume that $K_\epsilon$ is compact. Then for all $x\in K_\epsilon$ such that $\dim I(x)=d-1$, $I(x)$ contains a closed $d-1-$dimensional ball $B_x:=I(x)\cap B_{r_x}(x)$ for some $r_x>0$. As $I$ is continuous on $K_\epsilon$, we may find $\epsilon_x>0$ such that for $x'\in B_{\epsilon_x}(x)$, $B_x\subset\proj_{\aff I(x)}\big(I(x')\big)$, and such that the angle between the normals of $I(x)$ and $I(x')$ is smaller than $\eta :=\pi/4<\pi/2$. We denote $l_x$ the line from $x$, normal to $I(x)$. The balls $B_{\epsilon_x}(x)$ cover $K_\epsilon$, then by the compactness of $K_\epsilon$, we may consider $x_1,...,x_k\in K_\epsilon$ for $k\ge 1$ such that $K_\epsilon\subset \cup_{i=1}^k B_{\epsilon_{x_i}}(x_i)$. Let $1\le i\le k$, by Lemma C.1. in \cite{ghoussoub2015structure}, we may find a bi-Lipschitz flattening map $F:\cup_{x'\in A_i}I(x')\longrightarrow \R^d = \aff I(x_i)\x l_{x_i}$, where $A_i:=B_{\epsilon_{x_i}}(x_i)\cap l_{x_i}$, such that for all $x'\in A_i$ and all $(v,w)\in I(x')$, $F(v,w) = (v,x')$. Notice that for all $x'\in B_{\epsilon_{x_i}}(x_i)$, $I(x')\cap A_i\neq\emptyset$. Then for all $x'\in A_i$, $F\big(I(x')\big)\subset \aff I(x_i)\x\{x'\}$. Now, let $\lambda$ be the Lebesgue measure. By the Fubini theorem, $\lambda\big[F\big(\cup_{x'\in A_i}\partial I(x')\big)\big] = \int_{l_x}\mathbf{1}_{x'\in A_i}\lambda_{x'}\big[F\big(\partial I(x')\big)]dx'$. By the facts that $F$ is bi-Lipschitz, $\partial I(x')$ is Lebesgue-negligible in $\aff I(x')$, and $\lambda_{x'}$ is a $d-1-$dimensional Lebesgue measure, we have $\lambda_{x'}\big[F\big(\partial I(x')\big)\big]=0$, $\mathbf{1}_{x'\in A_i}dx'-$a.e. Therefore, $\lambda\big[F\big(\cup_{x'\in A_i}\partial I(x')\big)\big]=0$, and as $F$ is bi-Lipschitz, $\lambda[\cup_{x'\in A_i}\partial I(x')]=0$. Then summing up on all the $1\le i\le k$ and by the fact that $\nu$ is dominated by the Lebesgue measure, we get $\nu[\cup_{x\in K_{\epsilon}}\partial I(x)] = 0,$ so that for all $\P\in\Mc(\mu,\nu)$, we have
$$
\P[Y\in\partial I(X),\dim I(X)=d-1] \le \P[X\notin K_\epsilon,\dim I(X)=d-1]+\P[Y\in\cup_{x\in K_\epsilon}\partial I(x)]\le \epsilon.
$$
As this holds for all $\epsilon>0$ and for all $\P\in\Mc(\mu,\nu)$, the lemma is proved.
\ep\\



\no {\bf Proof of Proposition \ref{prop:emptyboundaries}}
Let us first prove the equivalence from (i). First for $\P\in\Mc(\mu,\nu)$. As $Y\in I(X)$, $\P$-a.s., we have $I(X)=I(Y)$, $\P-$a.s., and therefore, for all $A\in \Bc(\Kc)$,
\b*
\nu\circ I^{-1}[A] = \P[I(Y)\in A] = \P[I(X)\in A] = \mu[I(X)\in A] = \mu\circ I^{-1}[A]
\e*

Conversely, suppose that $\mu\circ I^{-1} = \nu\circ I^{-1}$. We will prove by backward induction on $0\leq k\leq d+1$ that $Y\in I(X)$, $\Mc(\mu,\nu)$-q.s., conditionally to $\dim I(X) \ge k$. For $k=d+1$ this is trivial because the dimension is lower than $d$. Now for $k\in\N$ we suppose that the property is true for $k'>k$. Then conditionally to $\dim I(x) = k$, we have that $Y\in \cl I(X)$, q.s. Then for $\P\in\Mc(\mu,\nu)$,
$$\P[\dim I(Y) = k] = \P[Y\in I(X)\text{ and }\dim I(X) = k]+\P[Y\in \partial I(X)\text{ and }\dim I(X) > k]$$
By the induction hypothesis, $\P[Y\in \partial I(X)\text{ and }\dim I(X) > k]=0$. (i) gives that $\P[\dim I(Y) = k] = \P[\dim I(X) = k]$. Then
$$\P[\dim I(X) = k]=\P[Y\in I(X)\text{ and }\dim I(X) = k],$$
implying that $\P-a.s.$, $\dim I(X) = k\implies Y\in I(X)$. As holds true for all $\P\in\Mc(\mu,\nu)$, combined with the induction hypothesis, we proved the result at rank $k$. By induction, $Y\in I(X)$, q.s. The equivalence is proved.

It remains to show that (iv) is implied by all the other conditions. If (i) holds, $\cup_{x\in\R^d} I(x)\x\partial I(x)\in\Nc_{\mu,\nu}$ then (iv) holds with $\Cc=\Dc=\emptyset$, and $\Rc:= I(\R^d)$. If (ii) holds, as $I(\R^d)$ is a partition of $\R^d$, there can be at most countably many components with full dimension. Therefore (iv) holds with $\Cc:=I(\{\dim I = d\})$, and $\Dc:=\{\dim I \le 1\}$.

Now we suppose (iii), by Lemma \ref{lemma:nocharge}, $Y\notin \partial I(X)$ if $\dim I(X)=d-1$, $\Mc(\mu,\nu)-$q.s. Then we just set $\Dc:=\{\dim I \le 1\}$, $\Cc:=\{\dim I = d\}$, and $\Rc:=\{\dim I = d-1\}$. Now we prove the claim.

We suppose that (iv) holds. The second part of the proposition follows from the fact that a countable set can be well ordered. Now let us deal with the first part. According to Lemma \ref{lemma:dominproba}, we just need to find a probability measure $\Q$ that implies the quasi-sure convergence of functions in $\widetilde{\Tc}_1$. This is possible thanks to the convexity of these functions in the second variable: the interior of the components can be dealt with $\mu(dx)\otimes \sum_{n\ge 1}2^{-n}\delta_{f_n(x)}(dy)$, where $(f_n)_{n\ge 1}\subset\Leb^0(\R^d,\R^d)$ is chosen such that $\{f_n(x):n\ge 1\} \subset \aff I(x)$ is dense in $\aff I(x)$ for all $x\in\R^d$ (see the proof of Lemma \ref{lemma:dominproba}).

For the boundaries, the measure $\mu\otimes\nu$ will deal with the countable components of $\Cc$. Indeed, let $K\in\Cc$ such that $\eta[K]>0$. Let $(\theta_n)_n\subset \widetilde{\Tc}_1$, converging $\mu(dx)\otimes \sum_{n\ge 1}2^{-n}\delta_{f_n(x)}(dy)+\mu\otimes\nu-$a.s. to some function $\theta$. We already have that $\theta_n(x,\cdot)\longrightarrow\theta(x,\cdot)$ on $K$ for all $x\in N_\mu^c\cap K$, for some $N_\mu\in\Nc_\mu$ by the previous step. For all $n\ge 1$, let $N_n\in\Nc_{\mu,\nu}$ be such that $\theta_n$ is a $N_n-$tangent convex function. By \eqref{eq:polar} and by possibly enlarging the $\mu-$null set $N_\mu$, we may assume that we may find $(N_\nu,\theta)\in\Nc_\nu\x\widehat{\Tc}(\mu,\nu)$ such that $N_\mu^c\x N_\nu^c\cap\{Y\in J_\theta(X)\}\subset N^c:=(\cup_{n\ge 1}N_n)^c$, and that $N_\mu^c\x \{Y\in I(X)\}\subset N^c$. Then for $x,x'\in N_\mu^c\cap K$, $x_0\in K$, and $y\in J_\theta(x)\setminus (K\cup N_\nu)$, let the probability measures
\b*
4\P:=\delta_{x,x_0}+\delta_{x,y}+2\delta_{x',y'},&\mbox{and}&
4\P':=\delta_{x',x_0}+\delta_{x',y}+2\delta_{x,y'}
\e*
with $y' := \frac12(y + x_0)$. Let $n\ge 1$, notice that $\P$ and $\P'$ are competitors and concentrated on $N_n$, then by $\theta_n-$martingale monotonicity of $N_n$, we have
$$\theta_n(x,x_0)+\theta_n(x,y)+2\theta_n(x',y') = \theta_n(x',x_0)+\theta_n(x',y)+2\theta_n(x,y').$$
We re-order the terms
\be\label{eq:convinvariance}
\theta_n(x,y)-2\theta_n(x,y')+\theta_n(x,x_0) = \theta_n(x',y)-2\theta_n(x',y')+\theta_n(x',x_0).
\ee
Then $\theta_n(x,y)-2\theta_n(x,y')+\theta_n(x,x_0)$ does not depend on the choice of $x\in K\cap N_\mu^c$. As we assumed that $\theta_n$ converges $\mu(dx)\otimes \sum_{n\ge 1}2^{-n}\delta_{f_n(x)}(dy)+\mu\otimes\nu-$a.s. by possibly enlarging $N_\mu$, without loss of generality, we may assume that for all $x\in N_\mu^c$, $\theta_n(x,\cdot)$ converges pointwise to $\theta$ on $I(x)$, and $\theta_n(x,Y)$ converges $\nu-$a.s. Let $x'\in N_\mu^c\cap K$, up to enlarging $N_\nu$, we may assume that $\theta_n(x',y)$ converges to $\theta(x',y)$ for all $y\in N_\nu^c$. Then if $x,y\in (K\cap N_\mu^c)\x N_\nu^c$, and $x\in K$, identity \eqref{eq:convinvariance} implies that $\theta_n(x,y)$ converges, as all the other terms have a limit, and $\theta(x,y')$ and $\theta(x',y')$ are finite. Now for $\P\in\Mc(\mu,\nu)$, $\P[(K\cap N_\mu^c)\x N_\nu^c] =\eta[K]$. Then $\theta_n$ converges $\P-$a.s. on $K\x\R^d$. This holds for all $K\in\Cc$, and $\P\in\Mc(\mu,\nu)$.

For the $1$-dimensional components of $\Dc$, if we call $a(x)$ and $b(x)$ their (measurably selected) endpoints, the measure $\mu(dx)\otimes \frac{\delta_{a(x)}+\delta_{b(x)}}{2}$ will fit. Finally, in the case of the components in $\Rc$, for all probability $\P\in\Mc(\mu,\nu)$, $\P_x$  does not send mass to $\partial K$ for $\mu-$a.e. $x\in K\in\Rc$ by assumption. We take
$$\Q(dx,dy) := \mu(dx) \sum_{n\ge 1}2^{-n}\delta_{f_n(x)}(dy) + \mu(dx)\nu(dy) + \mu(dx) \frac{\delta_{a(x)}+\delta_{b(x)}}{2}(dy).$$
the convergence of $\theta_n$, $\Q-$a.s. implies its convergence $\Mc(\mu,\nu)-$q.s. Assumption \ref{ass:domination} holds.
\ep\\

\no{\bf Proof of Remark \ref{rmk:invariancenu}}
The fact the $\nu_I^\P$ is independent of $\P\in\Mc(\mu,\nu)$ for $d=1$ is proved by Beiglb\"ock \& Juillet \cite{beiglboeck2016problem}.

Now we assume that (i) in Proposition \ref{prop:emptyboundaries} holds. If $Y\in I(X)$, $\Mc(\mu,\nu)-$q.s., then by symmetry as $\{I(x):x\in\R^d\}$ is a partition of $\R^d$, we have $X\in I(Y)$, $\nu-$a.s. Then similar to $\mu_I$, $\nu_I := \nu^\P_I:=\P\circ (Y|X\in I)^{-1}$ does not depend on the choice of $\P\in\Mc(\mu,\nu)$.

Now in the case of (iii) in Proposition \ref{prop:emptyboundaries}, let $\nu_I:=\nu\circ I^{-1}$. On $\{\dim I(X) \ge d-1\}$, $Y\notin \partial I(X)$, q.s. by Lemma \ref{lemma:nocharge}, so that for all $\P\in\Mc(\mu,\nu)$, $\nu_I^\P = \nu_I$ on $\{\dim I(X) \ge d-1\}$. Now on $\{\dim I(X) =0\}$, $\mu_I = \nu_I^\P$ is also independent of $\P$. Finally, on $\{\dim I(X) =1\}$, by the fact that there is not mass coming from higher dimensional components, we have $\nu_I^\P=\nu_I+\lambda_1(I)\delta_{a(I)}+\lambda_2(I)\delta_{b(I)}$, where $\lambda_1(I),\lambda_2(I)\ge 0$, and $a(I),b(I)$ are measurable selections of the boundary of $I$. Then $\mu_I-\nu_I = \lambda_1(I)+\lambda_2(I)$, and $\mu_I[X]-\nu_I[Y]=\lambda_1(I)a(I)+\lambda_2(I)b(I)$. Therefore, $\lambda_1$ and $\lambda_2$ depend only on $\mu_I$ and $\nu_I$, therefore, $\nu_I^\P$ does not depend on the choice of $\P$.
\ep\\

\no{\bf Proof of Remark \ref{rmk:diffusion}}
We consider $\tau$ the stopping time, and write $\Q$ the probability measure associated with the diffusion. We claim that the components $\csupp\P_{X_0}\subset I(X_0)$, $\mu-$a.s. have dimension $d$, $\mu$-a.s, where $\P\in\Mc(\mu,\nu)$ is the joint law of $(X_0,X_\tau)$. Then (iii) of Proposition \ref{prop:emptyboundaries} holds, which proves the remark.

Now we prove the claim. Let $p>0$. For $x\in\R^d$, we consider $\tau_x$, the stopping time $\tau$ conditional to $X_0=x$, and $\sigma^x_t$, which is $\sigma_t$ conditional to $X_0=x$. Now we fix $x\in\R^d$. As $\sigma_0$ has rank $d$, $\Vert\sigma^x_0\Vert:=\inf_{|u|=1}|u^t\sigma^x_0|>0$, a.s. Then we may find $\alpha >0$ such that $\Q[\Vert\sigma^x_0\Vert\le \alpha]\le p$. Similarly, we consider $\delta>0$ small enough so that
\be\label{eq:stoppingbound1}
\Q[\tau <\delta]\le p.
\ee
Finally, by the fact that $\sigma^x_t$ is right-continuous in $0$, a.s, we may lower $\delta>0$ so that $\Q\left[\sup_{t\le \delta}|\sigma^x_t-\sigma^x_0|^2> \beta\right]\le  p$ for some $\beta>0$ that we will fix later. Now we use these ingredients to prove that $(X_t)$ "spreads out in all directions" for $t$ close to $0$. Let $u\in\R^d$ with $|u|=1$ and $\lambda >0$,
\be\label{eq:brownianineq} 
\Q[u\cdot\sigma^x_0W_\delta\ge \lambda\alpha\sqrt{\delta}]\ge \Q[ v\cdot W_1\ge \lambda]-p \ge \frac12 - 2p,
\ee
with $v = u\cdot\sigma^x_0/|u\cdot\sigma^x_0|$, for $\lambda$ small enough, independent of $\alpha$ and $\delta$. Now recall that $\Q\left[\sup_{t\le \delta}|\sigma^x_t-\sigma^x_0|^2> \beta\right]\le p$. As a consequence, the stopping time $\tilde{\tau} =\inf\{t\ge 0 :|\sigma^x_t-\sigma^x_0|^2\ge \beta\}$ satisfies
\be\label{eq:stoppingbound2}
\Q[\tilde{\tau}< \delta]\le p.
\ee
Now, stopping $X_t$, we get, conditionally to $X=x$: $\E^\Q[(\int_0^{\delta\wedge\tilde{\tau}} (\sigma_t^x-\sigma_0^x)dW_t)^2]\le \delta\beta$ by It\^o isometry, and therefore, by the Markov inequality, $\Q\left[\left|\int_0^{\delta\wedge\tilde{\tau}}(\sigma_t^x-\sigma_0^x)dW_t\right|\ge\alpha\lambda\sqrt{\delta}/2\right]\le \frac{ 4\delta\beta}{\alpha^2\lambda^2\delta}$. Then if we chose $\beta = p\frac{\alpha^2\lambda^2}{4}$ (not depending on $\delta$), we finally get that
\be\label{eq:probalessthanp}
\Q\left[\left|\int_0^{\delta\wedge\tilde{\tau}}(\sigma_t^x-\sigma_0^x)dW_t\right|\ge\alpha\lambda\sqrt{\delta}/2\right]\le p.
\ee
Therefore $\Q[(X_{t\wedge\tau}-x)\cdot u\ge \alpha\lambda\sqrt{\delta}/2|X=x]$ is greater than
\b*
\Q\left[\sigma^x_0 W_{t\wedge\tau}\cdot u\ge \alpha\lambda\sqrt{\delta}\mbox{, and }\left|\int_0^{\delta}(\sigma_t^x-\sigma_0^x)dW_t\right|\le\alpha\lambda\sqrt{\delta}/2\mbox{, and }\tilde{\tau}\ge\delta\mbox{, and } \tau\ge\delta|X=x\right]\\
\ge \Q[u\cdot\sigma^x_0W_\delta\ge \lambda\alpha\sqrt{\delta}]-3p\ge \frac12 - 5p,
\e*
by \eqref{eq:stoppingbound1}, \eqref{eq:brownianineq}, \eqref{eq:stoppingbound2}, and \eqref{eq:probalessthanp}. Then by setting $p = \frac{1}{12}$, for all $u$ of norm $1$, we get
\be\label{eq:goindirection}
\Q[(X_{t\wedge\tau}-x)\cdot u\ge \alpha_0|X_0=x] \ge p_0,
\ee
with $\alpha_0 := \alpha\lambda\sqrt{\delta}/2>0$, and $p_0 := \frac{1}{12}>0.$

We use \eqref{eq:goindirection} to prove that $\csupp\P_x$ is $d$ dimensional. Indeed, we suppose for contradiction that $\csupp\P_x\subset H$, where $H$ is a hyperplane. $H$ contains $0$, as it contains $\csupp\P_x$. Let $u$ be a unit normal vector to $H$, by \eqref{eq:goindirection}, we have $\Q[(X_{t\wedge\tau}-x)\cdot u\ge \alpha_0|X=x]\ge p_0$. Then by the martingale property (the volatility is bounded) combined with the boundedness of $\tau$, we have $\E^\Q[X_{\tau}|\Fc_{t\wedge\tau}] = X_{t\wedge\tau}$. Therefore, $\P_x[Y\cdot u \ge \alpha_0/2] = \Q[X_{\tau}\cdot u \ge \alpha_0/2|X=x]>0$, which contradicts the inclusion of the support of $\P_x$ in $H$.
\ep

\subsection{Medial limits}\label{subsect:verifassaxiom}

Medial limits, introduced by Mokobodzki \cite{mokobodzki1967ultrafiltres} (see also Meyer \cite{meyer1973limites}), are powerful instruments. It is an operator from the set of real bounded sequences $l^\infty$ to $\R$ satisfying the following properties:

\begin{Definition}\label{def:medial}
A linear operator $\med : l^\infty\to\R$ is a medial limit if\\
{\rm (i)} $\med$ is nonnegative: if $u\ge 0$ then $\med(u)\ge 0$.\\
{\rm (ii)} $\med$ is invariant by translation: if $\Tc$ is the translation operator ($\Tc:(u_n)_n\mapsto (u_{n+1})_n$) then $\med(\Tc u) = \med(u)$.\\
{\rm (iii)} $\med((1)_n) = 1$.\\
{\rm (iv)} $\med$ is universally measurable on the unit ball $[0,1]^\N$.\\
{\rm (v)} $\med$ is measure linear: for any sequence of Borel-measurable functions $f_n:[0,1]\to [0,1]$, if we write $f:=\med((f_n)_n)$ (defined pointwise), then for any Borel measure $\lambda$ on $[0,1]$, $f$ is $\lambda$-measurable and
$$\int f d\lambda = \int \med(f_n)d\lambda = \med\left(\int f_n d\lambda\right).$$
\end{Definition}

We can extend any medial limit $\med$ to $\R_+^\N$ by setting $\med(u):=\sup_{N\in\N}\med((u_n\wedge N)_n)$. It keeps the same properties, except {\rm (v)} which becomes a kind of Fatou's Lemma: for any sequence of Borel-measurable functions $f_n:[0,1]\to \R_+$, then for any Borel measure $\lambda$ on $[0,1]$,
\be\label{eq:Fatoumed}
\int \med(f_n)d\lambda &\leq& \med\left(\int f_n d\lambda\right).
\ee
The existence of medial limits is implied by Martin's axiom. Notice that Martin's axiom is implied by the continuum hypothesis (See Chapter I of Volume 5 of \cite{fremlin2000measure}). Kurt G\"odel \cite{godel1947cantor} provides 6 paradoxes implied by the continuum hypothesis, Martin's axiom implies only 3 of these paradoxes. All these axioms are undecidable either under ZF and under ZFC, indeed Paul Larson \cite{larson2009filter} proved that if ZFC is consistent, then ZFC+"there exists no medial limits" is also consistent (Corollary 3.3 in \cite{larson2009filter}). See \cite{limmed54562} for a complete survey.\\

\no {\bf Proof of Proposition \ref{prop:axioms}}
Axiom of choice on $\R$ implies that $\R$ can be well-ordered, which proves that Assumption \ref{ass:domination} (ii) holds. Now let us prove the first part. For $(\theta_n)_{n\ge 1}\subset\widetilde{\Tc}_1$, we denote $\theta:=\med(\theta_n)$. The Proposition is proved if we show that $\theta_n\rightsquigarrow\theta$. $\theta=\med(\theta_n)\ge \underline{\theta}_\infty$ by linearity of a medial limit together with Definition \ref{def:medial} (i) and (ii). Let $\P\in\Pc(\Omega)$, $\P[\theta]\le \med(\P[\theta_n])\le\limsup_{n\to\infty}\P[\theta_n]$ by \eqref{eq:Fatoumed}. Finally the linearity combined with Definition \ref{def:medial} (i) give that $\theta\in\Theta_{\mu,\nu}$, as it is a property of comparison of linear combinations of values of $\theta$, $\theta$ is a $\emptyset-$tangent convex function. Finally, we prove that we may have (iv) in Definition \ref{def:Thetamunu}. Up to assuming that we applied the Koml\'os Lemma to $(\theta_n)_{n\ge 1}$ (which only reduces the superior limits and increase the inferior limits, thus preserving the previous properties) under the probability $\mu(dx)\otimes \sum_{n\ge 1}2^{-n}\delta_{f_n(x)}(dy)$, where $(f_n)_{n\ge 1}\subset\Leb^0(\R^d,\R^d)$ is chosen such that $\{f_n(x):n\ge 1\} \subset \aff I(x)$ is dense in $\aff I(x)$ for all $x\in\R^d$ as in the proof of Lemma \ref{lemma:dominproba}, we may assume without loss of generality that $(\theta_n)$ converges pointwise on $\{X\in N_\mu'^c\}\cap\{Y\in I(X)\}$. Then let $N_\mu^n\in\Nc_\mu$ be from Definition \ref{def:Thetamunu} (iv) for $\theta_n$. Let $N_\mu = \cup_{n\ge 1}N_\mu^n\cup N_\mu'$. Let $A:=\{X\in N_\mu^c\}\cap\{Y\in I(X)\}$, $\mathbf{1}_{A}\theta$ is Borel measurable as the pointwise limit of Borel measurable functions $\mathbf{1}_{A}\theta_n$, as the medial limit coincides with the real limit when convergence holds.
\ep

\bibliographystyle{plain}
\bibliography{mabib}

\end{document}